\documentclass[a4paper,10pt]{article}
\makeatletter

\usepackage[utf8]{inputenc}
\usepackage[T1]{fontenc}
\usepackage[czech,english]{babel}

\usepackage{geometry}
\usepackage{ifpdf}
\ifpdf
        \@ifpackageloaded{hyperref}{}{\usepackage[pdftex]{hyperref}}
	\@ifpackageloaded{xcolor}{}{\usepackage[pdftex,svgnames,table]{xcolor}}
        \@ifpackageloaded{graphicx}{}{\usepackage[pdftex]{graphicx}}
	\usepackage{epstopdf}
\else
        \@ifpackageloaded{hyperref}{}{\usepackage[hypertex]{hyperref}}
	\@ifpackageloaded{xcolor}{}{\usepackage[dvips,svgnames,table]{xcolor}}
        \@ifpackageloaded{graphicx}{}{\usepackage[dvips]{graphicx}}
\fi
\graphicspath{{fig/}}

\hypersetup{
    breaklinks=true,
    colorlinks=true,
    linkcolor=myblue,
    citecolor=myblue,
    urlcolor=myblue,
    unicode,
    naturalnames
}
\usepackage{xurl} %
\usepackage{amssymb,amsfonts}
\usepackage{lmodern}
\usepackage{fancyhdr}

\usepackage{amsmath,amsthm,amscd,upref}
\usepackage{mathtools}
\mathtoolsset{showonlyrefs=true} %

\usepackage{paralist} %
\usepackage{tikz}

\usepackage[round,sort]{natbib}
\let\oldcite=\cite
\let\oldcitep=\citep

\renewcommand\cite[1]{\hypersetup{linkcolor=Navy}\hyperlink{#1}{\oldcite{#1}}}
\renewcommand{\citep}[1]{\hypersetup{linkcolor=Navy}\hyperlink{#1}{\oldcitep*{#1}}}
\renewcommand{\citet}[1]{\hypersetup{linkcolor=Navy}\hyperlink{#1}{\oldcite*{#1}}}

\newcommand{\doi}[1]{DOI~\href{\detokenize{https://doi.org/#1}}{\detokenize{#1}}}
\newcommand{\zblnumber}[1]{Zbl~\href{\detokenize{https://zbmath.org/?q=an:#1}}{\detokenize{#1}}}
\newcommand{\mrnumber}[1]{\href{\detokenize{https://www.ams.org/mathscinet-getitem?mr=#1}}{\detokenize{MR#1}}}

\newcommand{\msc}[1]{\href{https://mathscinet.ams.org/mathscinet/msc/msc2020.html?s=#1}{#1}}
\newcommand{\acm}[1]{\href{https://www.acm.org/publications/computing-classification-system/1998/#1}{\uppercase{#1}}}
\newcommand{\jel}[1]{\href{https://www.aeaweb.org/jel/guide/jel.php}{#1}}%

\usepackage{listings}
\colorlet{myblue}{Navy}
\definecolor{favorange}{cmyk}{0.02,0.22,1.0,0.15}
\definecolor{lightyellow}{cmyk}{0,0,0.1,0}

\definecolor{myfullred}{RGB}{244,63,43}
\definecolor{myhalfred}{RGB}{158,22,7}
\definecolor{myfullorange}{RGB}{255,128,0}
\definecolor{myhalforange}{RGB}{147,73,0}
\definecolor{mygreen}{RGB}{70,180,5}
\definecolor{mylilas}{RGB}{170,55,241}

\colorlet{pagebgcolor}{white}
\colorlet{textcolor}{black}
\colorlet{titlebgcolor}{pagebgcolor}

\usepackage{xspace}

\usepackage{soul}
\sethlcolor{LimeGreen}

\soulregister\cite7%
\soulregister\citep7%
\soulregister\eqref7%
\soulregister\pageref7%
\soulregister\ref7%

\def\revised#1{Revised~#1\par}

\@ifundefined{ccode}{%
\newcommand{\ccode}[2]{\par
        \vspace*{8pt}
        {{\leftskip18pt\rightskip\leftskip
        \noindent{\it #1}\/: #2\par}}\par}
}{}
\@ifundefined{keywords}{\newcommand{\keywords}[1]{\ccode{Keywords}{#1}}}{}
\@ifundefined{email}{\newcommand{\email}[1]{\href{mailto:#1}{#1}}}{}

\usepackage{xifthen}
\newboolean{default}
\newcommand{\case}{}
\newcommand{\default}{}
\newenvironment{switch}[1]{%
    \setboolean{default}{true}
    \renewcommand{\case}[2]{\ifthenelse{\equal{#1}{##1}}{%
        \setboolean{default}{false}##2}{}}%
    \renewcommand{\default}[1]{\ifthenelse{\boolean{default}}{##1}{}}
}{}

\@ifclassloaded{article}{%

\geometry{a4paper,top=20mm,left=20mm,right=20mm,bottom=20mm} %
\setlength{\headheight}{18.0pt} %
\rhead[]{\fancyplain{}{\slshape\rightmark}}
\chead{}
\lhead[\fancyplain{}{\slshape\leftmark}]{}
\cfoot{\fancyplain{\bfseries --~\thepage~--}{\bfseries --~\thepage~--}}

\usepackage{authblk}
\renewcommand\section{\@startsection {section}{1}{\z@}%
                                   {-3.5ex \@plus -1ex \@minus -.2ex}%
                                   {2.3ex \@plus.2ex}%
                                   {\color{myblue}\normalfont\Large\bfseries }}

}{}%

\@ifclassloaded{beamer}{%
\usetheme{default}
\setbeamertemplate{navigation symbols}{}
\setbeamertemplate{sidebar left}{\rlap{\hskip1mm\insertlogo}}
\setbeamertemplate{sidebar right}{}
\newcommand{\myfootlinefull}
{%
  \leavevmode%
  \hbox{\begin{beamercolorbox}[wd=\paperwidth,ht=2.5ex,dp=1.125ex,leftskip=3mm,rightskip=3mm]{author in head/foot}%
  \usebeamerfont{author in head/foot}\insertshortdate\hfill\insertshortauthor:\:
  \usebeamerfont{title in head/foot}\insertshorttitle\hfill \insertframenumber~/~\inserttotalframenumber
  \end{beamercolorbox}}%
  \vskip0pt%
}

\setbeamertemplate{footline}{\myfootlinefull}
\setbeamercolor{background canvas}{bg=pagebgcolor}
\setbeamercolor{alerted text}{fg=Navy}
\setbeamercolor{structure}{fg=Navy}
\setbeamercolor{frametitle}{fg=textcolor}
\setbeamercolor{author in head/foot}{fg=textcolor}
\setbeamertemplate{frametitle}[default][right,rightskip=1mm]
\addtobeamertemplate{frametitle}{}{\noindent\makebox[\linewidth]{\textcolor{favorange}{\rule{\paperwidth}{1pt}}}}
\setbeamerfont{title}{shape=\scshape}

\setbeamertemplate{bibliography item}{}
\setbeamersize{text margin left=5mm,text margin right=5mm}

\addtobeamertemplate{theorem begin}{%
    \setbeamercolor{block title}{fg=black,bg=favorange}%
    \setbeamercolor{block body}{fg=black,bg=lightyellow}%
}{}

\logo{%
\begin{tabular}{p{26mm}p{30mm}}%
\vspace{0mm}\href{http://www.fav.zcu.cz/en}{\includegraphics[height=13mm,trim=5mm 5mm 7mm 5mm,clip]{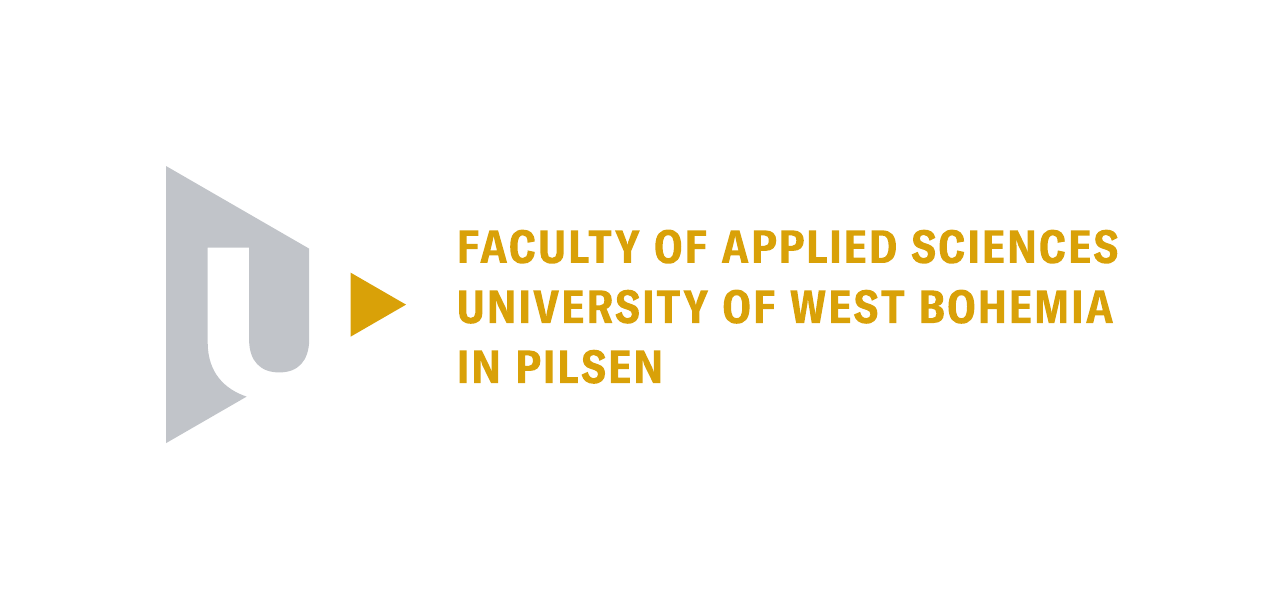}}&%
\vspace{1.5mm}\href{http://www.kma.zcu.cz/en}{\includegraphics[height=13mm,trim=0mm 0mm 135mm 0mm,clip]{logo-kma_en}}
\end{tabular}
}
\titlegraphic{%
\begin{tabular}{p{26mm}p{30mm}p{70mm}}%
\vspace{0mm}\href{http://www.fav.zcu.cz/en}{\includegraphics[height=13mm,trim=5mm 5mm 7mm 5mm,clip]{logo-fav_en}}&%
\vspace{1.5mm}\href{http://www.kma.zcu.cz/en}{\includegraphics[height=13mm,trim=0mm 0mm 135mm 0mm,clip]{logo-kma_en}}&%
\vspace{0mm}\hfill\href{http://www.ntis.zcu.cz/en}{\includegraphics[height=13mm,trim=0mm 0mm 0mm 0mm,clip]{logo-ntis}}
\end{tabular}
}

}{} %

\renewcommand{\d}{\,\mathrm{d}} %

\newcommand{\e}{\mathrm{e}}     %

\DeclareMathOperator{\logp}{log1p}
\newcommand{\NaN}{\mathrm{NaN}}

\newcommand{\N}{\mathbf{N}}	%
\newcommand{\R}{\mathbf{R}}	%

\addto\captionsczech{%
  \def\theoremname{V\v{e}ta}%
  \def\definitionname{Definice}%
  \def\examplename{P\v{r}\'{\i}klad}%
  \def\remarkname{Pozn\'{a}mka}%
  \def\propositionname{Tvrzen\'{\i}}%
  \def\lemmaname{Lemma}%
  \def\corollaryname{D\r{u}sledek}%
}%
\addto\captionsenglish{%
  \def\theoremname{Theorem}%
  \def\definitionname{Definition}%
  \def\examplename{Example}%
  \def\remarkname{Remark}%
  \def\propositionname{Proposition}%
  \def\lemmaname{Lemma}%
  \def\corollaryname{Corollary}%
}%

\theoremstyle{definition}

\theoremstyle{remark}
       
\theoremstyle{plain}

\newcommand{\nocontentsline}[3]{}%
\makeatother

\newcommand{\alert}[1]{\textcolor{myblue}{#1}}
\newcommand{\jpTitle}{Challenges in automatic differentiation and numerical integration in physics-informed neural networks modelling}

\newcommand{\jpKeywords}{variable precision arithmetic; automatic differentiation; numerical integration; adaptive quadrature; physics-informed neural networks}
\newcommand{\jpMSC}{\msc{65D30}; \msc{65G30}; \msc{68T07}}
\newcommand{\jpACM}{\acm{g.1.0}; \acm{g.1.4}; \acm{i.2.6}; \acm{j.2}}
\newcommand{\jpJEL}{\jel{C63}}
\newcommand{\jpDateRevised}{25 July, 2025}

\newcommand{\jpDate}{}%

\author[1]{Josef Dan\v{e}k}%
\author[1]{Jan Posp\'{\i}\v{s}il\thanks{Corresponding author, \email{honik@kma.zcu.cz}}}
\affil[1]{Department of Mathematics, Faculty of Applied Sciences, \authorcr University of West Bohemia in Pilsen, Univerzitn\'{\i} 2732/8, 301 00 Plze\v{n}, Czech Republic,\vspace*{3pt}}%

\title{\textcolor{Navy}{\textsc{\jpTitle}}}
\date{\jpDate}

\begin{document}

\maketitle

\begin{center}
\revised{\jpDateRevised}
\end{center}

\begin{abstract}
In this paper, we numerically examine the precision challenges that emerge in automatic differentiation and numerical integration in various tasks now tackled by physics-informed neural networks (PINNs). Specifically, we illustrate how ill-posed problems or inaccurately computed functions can cause serious precision issues in differentiation and integration. A major difficulty lies in detecting these problems. A simple large-scale view of the function or good-looking loss functions or convergence results may not reveal any potential errors, and the resulting outcomes are often mistakenly considered correct. To address this, it is critical to determine whether standard double-precision arithmetic suffices or if higher precision is necessary. Three problematic use-cases for solving differential equations using PINNs are analysed in detail. For the case requiring numerical integration, we also evaluate several numerical quadrature methods and suggest particular numerical analysis steps to choose the most suitable method.

 \end{abstract}

\ccode{MSC classification}{\jpMSC}
\ccode{ACM classification}{\jpACM}
\ccode{JEL classification}{\jpJEL}

\setcounter{tocdepth}{2}
\tableofcontents
\clearpage

\section{Introduction}\label{sec:introduction}

Floating-point numbers are widely used in physics-related numerical calculations because they are easy to understand in this context. They have become the standard method for representing real numbers in computers. The IEEE-754 Standard, established in 1985, has significantly contributed to the widespread adoption of floating-point arithmetic by providing clear and practical guidelines for its implementation. To recall the history of floating-point arithmetic as well as its specification mandated by the IEEE-754 Standard we refer to the recent survey by \cite{Boldo2023floating}.

However, in many mathematical models we can observe that IEEE-754 Standard 32-bit or 64-bit floating-point format (denoted \texttt{binary32} or \texttt{binary64} respectively) is not always fully sufficient, see for example works by \citet{Bailey05g}, who studied especially applications in physics \citep{Bailey12,Bailey15}. Computations that require higher than \texttt{double} precision for robust and exact decision making were introduced in \cite{Pal04}. For a brief overview of numerical quadratures and their failures we refer the reader to the Section~4 in \cite{DanekPospisil20ijcm}. 

Solving differential equations of physical phenomena using machine learning of neural networks (NNs) has emerged recently as a new field of scientific machine learning. Physics-informed neural networks (PINNs) are a type of universal function approximators that can embed the knowledge of any physical laws directly to the learning process \citep{Raissi2019pinns}. PINNs modelling tasks involve among others an automatic differentiation \citep{Baydin2018AD} and a numerical calculation of integrals that usually depend on several model parameters. Optimization tasks consist of large number of integral evaluations with high precision and low computational time requirements. However, for some model parameters, many numerical quadrature algorithms fail to meet these requirements. We can observe an enormous increase in function evaluations, serious precision problems and a significant increase of computational time.

To overcome the problems caused by the floating-point arithmetic limits, high-precision or variable-precision arithmetic is a rapidly growing part of scientific computing environments. In the above mentioned papers we can find among others currently available software packages for high-precision floating-point arithmetic. For example in MATLAB, there exists a possibility of defining variables and perform numerical calculations in variable-precision arithmetic using the \texttt{vpa} bundle that is part of the Symbolic Math Toolbox (according to the documentation, \texttt{vpa} was introduced already before R2006a). Variable in the name suggests that a user can set the number of significant digits arbitrary, by default it is 32 significant digits. Numerical integration using variable precision was introduced later in MATLAB R2016b, however, it can be used only for functions that can be easily converted to symbolic expressions.

In recent years we can identify a systematic research in the area of computer arithmetic that are regularly presented for example at the IEEE International Symposia on Computer Arithmetic (ARITH) or newly at International Conferences for Next Generation Arithmetic (CoNGA). One area of research focus on correctly rounded evaluation of elementary and special functions, see for example a survey by \cite{Brisebarre2024crsurvey}. Some techniques are based on the error free transformations of the elementary floating point operations (additions, subtraction, multiplication and division) such as \texttt{TwoSum} algorithm by \cite{Knuth97taocp2ed3} that requires 6 floating point operations (flop) or \texttt{TwoProd} algorithm by \cite{Dekker71} that require 17 flop. These two algorithms are for example used in the compensated Horner's scheme to evaluate accurately the polynomials \citep{Langlois2007how}. 

In terms of accurate higher precision numerical quadratures, the research is not at all systematic and the majority papers focus only on solutions tailored to specific problems. For papers with a more general impact we refer to the results about Gauss--Legendre quadrature by \cite{Fousse2007accurate}, Newton--Cotes quadratures by \cite{Fousse2007multiple} or Gauss--Tur{\'a}n quadrature by \cite{Gautschi2014high}. A comparison of three high-precision quadrature schemes was performed on several testing examples by \cite{Bailey05a}, namely the Gaussian quadrature (\texttt{quadgs}), an error function based quadrature scheme (\texttt{quaderf}), and a $\tanh$--$\sinh$ quadrature (\texttt{quadts}) were analysed and compared.

In our original conference paper \cite{DanekPospisil15tcp} we show that even a very simple quadratic function can produce serious numerical problems if evaluated incorrectly. In a subsequent paper \cite{DanekPospisil20ijcm} we analysed more complicated integrands that occur in inverse Fourier transforms integrals. Although the motivation comes from the option pricing problems in mathematical finance, we can observe similar misbehaviour in many multi-physics and multi-scale models as well. In this paper we show that in some cases it is actually necessary to use the high or variable precision arithmetic (\texttt{vpa}) in order to get correct integration results that are crucial for the training of neural networks. The aim of this paper is to show how similar problems arise in PINNs modelling and how to handle the precision within desired bounds.

The structure of the paper is as follows. In Section \ref{sec2:inaccurate_functions}, we introduce the numerical precision challenges and give examples of inaccurately evaluated functions and their derivatives and integrals. Both a simple quartic polynomial example as well as a rather complicated function from the mathematical finance are presented. In Section \ref{sec3:quadratures}, we compare different numerical quadratures applied to studied integrals and suggest a usage the variable precision arithmetic in problematic cases. In the second part of the paper we study three problematic use-cases, in particular one linear and one non-linear ordinary differential equation and one linear partial differential equation are introduced in Section~\ref{sec:models} and solve them using PINNs in Section~\ref{sec:pinns}. 
We conclude in Section~\ref{sec:conclusion}.

\section{Inaccurately evaluated functions and their derivatives}\label{sec2:inaccurate_functions}

IEEE-754 Standard specifies \texttt{double} precision binary floating-point format as \texttt{binary64} having 1 sign bit, 11 bits for exponent and 53 bits for significand precision (1 bit is implicit and 52 bits are explicitly stored). Machine epsilon\footnote{For different definitions see for example \url{https://en.wikipedia.org/wiki/Machine_epsilon}.} (or machine precision) is an upper bound on the relative approximation error due to rounding in floating point number systems. In \texttt{binary64} this value (in particular the interval machine epsilon\footnote{It is the largest relative interval between two nearest numbers in finite-precision, i.e. $\epsilon = b^{-(p-1)}$, where $b$ is the base and $p$ is the precision.}) is $\epsilon = 2^{-52} \approx \texttt{2.22e-16}$. Among physical constants\footnote{See for example \url{https://en.wikipedia.org/wiki/List_of_physical_constants}.} we can find those that are many degrees smaller than this value and many being on the other hand extremely huge. Therefore, one can easily imagine many arithmetic operations with these constants in \texttt{double} precision can become problematic. 

In the past we witnessed already several examples when erroneous navigation input to the rocket's guidance system (not only in Ariane 5) lead to the placement of satellites to the wrong orbit. Many other particular examples not only from physics can be found for example in the papers by \cite{Bailey05g,Bailey12,Bailey15}. 

For the purposes of presented manuscript we introduce several numerical precision challenges in differentiation and integration of functions. In Subsection \ref{ss:phix} we describe in detail how a simple polynomial of fourth order can be inaccurately evaluated very easily. Accurate evaluation of functions is crucial not only for numerical quadratures, that are in detail studied in Section \ref{sec3:quadratures} below, but also in differentiation. Problems occur not only in numerical differentiation, but in particular also in automatic differentiation (see Section \ref{ss:diff}). In Subsection \ref{ss:fk} we present a more complex example coming from mathematical finance that was for the first time introduced by \cite{DanekPospisil20ijcm}. 

\subsection{Numerical precision challenges}\label{ss:list}

In the following, we list several important categories of problems that occur when evaluating functions and consequently in differentiation and integration. For each category we provide particular examples with short explanation.

\smallskip\noindent Subtractive-cancellation (loss of significance):
\begin{itemize}
\item $f(x) = 1-\dfrac{\sin x}{x}$, $x\in[10^{-12},10^{-2}]$, first 8–10 digits cancel, only round-off left;
\item $f(x) = \dfrac{1-\cos x}{x^2}$, $x\in[0,10^{-6}]$, both numerator and denominator are of the same order and close to zero, only round-off left;
\item $f(x) = \log(1+x)-x$, $x\in[0,10^{-9}]$, if $\logp$ function is not used, difference loses approximately $x/2$ digits.
\end{itemize}

\smallskip\noindent Extreme dynamic-range (overflow/underflow or huge ratio of scales):
\begin{itemize}
\item $f(x) = \e^{100x} - \e^{100(x+\varepsilon)}$, $\varepsilon\approx 10^{-8}$, $x\gtrsim 0.7$, each term overflows to $\infty$; difference becomes $\NaN$;
\item $f(x) = \exp(-1/x^2)$, $x\in[0,10^{-8}]$, underflows to 0, tail area completely missing;
\item $f(x) = \dfrac{1}{1+\e^{-1000x}}$, $|x|>0.02$, evaluates only to 0 or 1.
\end{itemize}

\smallskip\noindent Highly oscillatory or near-singular behaviour:
\begin{itemize}
\item $f(x) = \dfrac{\sin(k\cdot x)}{x}$, $k\approx 10^6$, $0<|x|<10^{-8}$, round-off dominates,
\item $f(x) = \sin(1/x)$, $x\in[0,10^{-2}]$, infinite oscillations as $x\to 0$, 
\item Bessel function of the first kind $f(x) = J_0(k\cdot x)$, $k\approx 10^6$, $x\in[0,1]$, root spacing approximately $\pi/k$ causes problems to numerical quadratures.
\end{itemize}

Apart from the above categories, one has to be careful also in other scenarios such as integrable algebraic or logarithmic singularities, piece-wise discontinuities (sharp corners), random-looking ``spiky'' data (e.g. noise-amplified derivatives in PINNs), very high dimension (already $d \gg 3$) or for example chaotic or stiff ODE-generated solutions (see also Sections \ref{ss:4:linode} and \ref{ss:4:nonlinode} below).

\subsection{Problematic quartic function example}\label{ss:phix} %

Let us now consider a simple continuous bounded real function $\phi(x)$ of one real variable 
\begin{align}
\phi(x)&=(b^2(x)-a^2(x))/\delta, \label{e:phi_delta} 
\intertext{where} 
a(x)&=x+\delta^2, \\
b(x)&=\sqrt{a^2(x)+\delta(x-1)^2(x+1)^2}.
\end{align}
Function $\phi(x)$ depends \emph{internally} on some real parameter $\delta>0$. This particular example has its simplified form 
\begin{equation}\label{e:phi_exact}
\phi(x)=(x-1)^2(x+1)^2
\end{equation}
that even does not depend on $\delta$. It is worth to mention that such simplification is usually not possible like in the example in Section~\ref{ss:fk}. We denote the MATLAB implementations of both forms as \texttt{phi\_double()} (that implements \eqref{e:phi_delta}) and \texttt{phi\_exact()} (that implements \eqref{e:phi_exact}), see Figure~\ref{f:phi_double+exact}.

\begin{figure}[ht!]
\newlength\jpColsep
\setlength\jpColsep{5mm}
\noindent\begin{minipage}[t]{\dimexpr0.5\textwidth-0.5\jpColsep\relax}
{%
\begin{lstlisting}
function value=phi_double(x,opt)
    delta = opt.delta;
    a     = x + delta.^2;
    b     = (a.^2+delta.*(x-1).^2.*(x+1).^2).^.5;
    value = (b.^2-a.^2)./delta; 
end
\end{lstlisting}
}
\end{minipage}\hfill
\begin{minipage}[t]{\dimexpr0.5\textwidth-0.5\jpColsep\relax}
{%
\begin{lstlisting}
function value=phi_exact(x,opt)
    % delta = opt.delta;
    value = (x-1).^2.*(x+1).^2;
end 
\end{lstlisting} 
}
\end{minipage}
\caption{MATLAB implementation of functions \texttt{phi\_double()} and \texttt{phi\_exact()}.}\label{f:phi_double+exact}
\end{figure}

We further introduce the so called variable precision arithmetic (\texttt{vpa}), also called high or arbitrary precision arithmetic. In MATLAB, \texttt{vpa(x,d)} uses variable precision arithmetic to evaluate each element of the symbolic input \texttt{x} to at least \texttt{d} significant digits, where \texttt{d} is the value of the \texttt{digits} function, whose default value is 32. MATLAB implementation \texttt{phi\_vpa()} of our function in variable precision arithmetic is listed in Figure~\ref{f:phi_vpa}. Please note that the high precision value is at the end converted back to \texttt{double}.

\begin{figure}[ht!]
\begin{lstlisting}
function value=phi_vpa(x,opt)
    delta = opt.delta;
    a     = vpa(x + delta.^2);
    b     = vpa((a.^2+delta.*(x-1).^2.*(x+1).^2).^.5);
    value = vpa((b.^2-a.^2)./delta);
    value = double(value);
end
\end{lstlisting} 
\caption{MATLAB implementation of function \texttt{phi\_vpa()}.}\label{f:phi_vpa}
\end{figure}

\begin{figure}[ht!]
\includegraphics[width=0.49\textwidth]{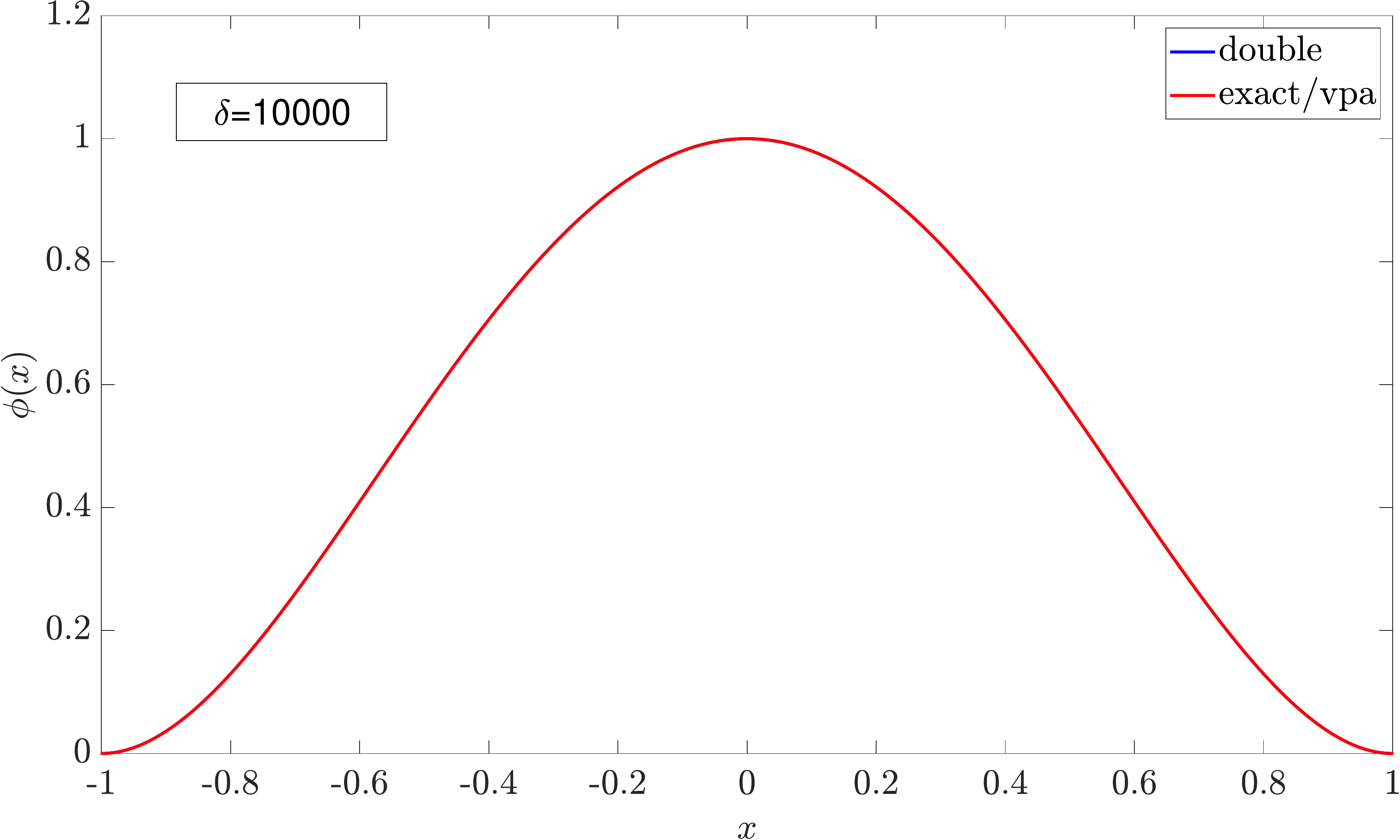}
\includegraphics[width=0.49\textwidth]{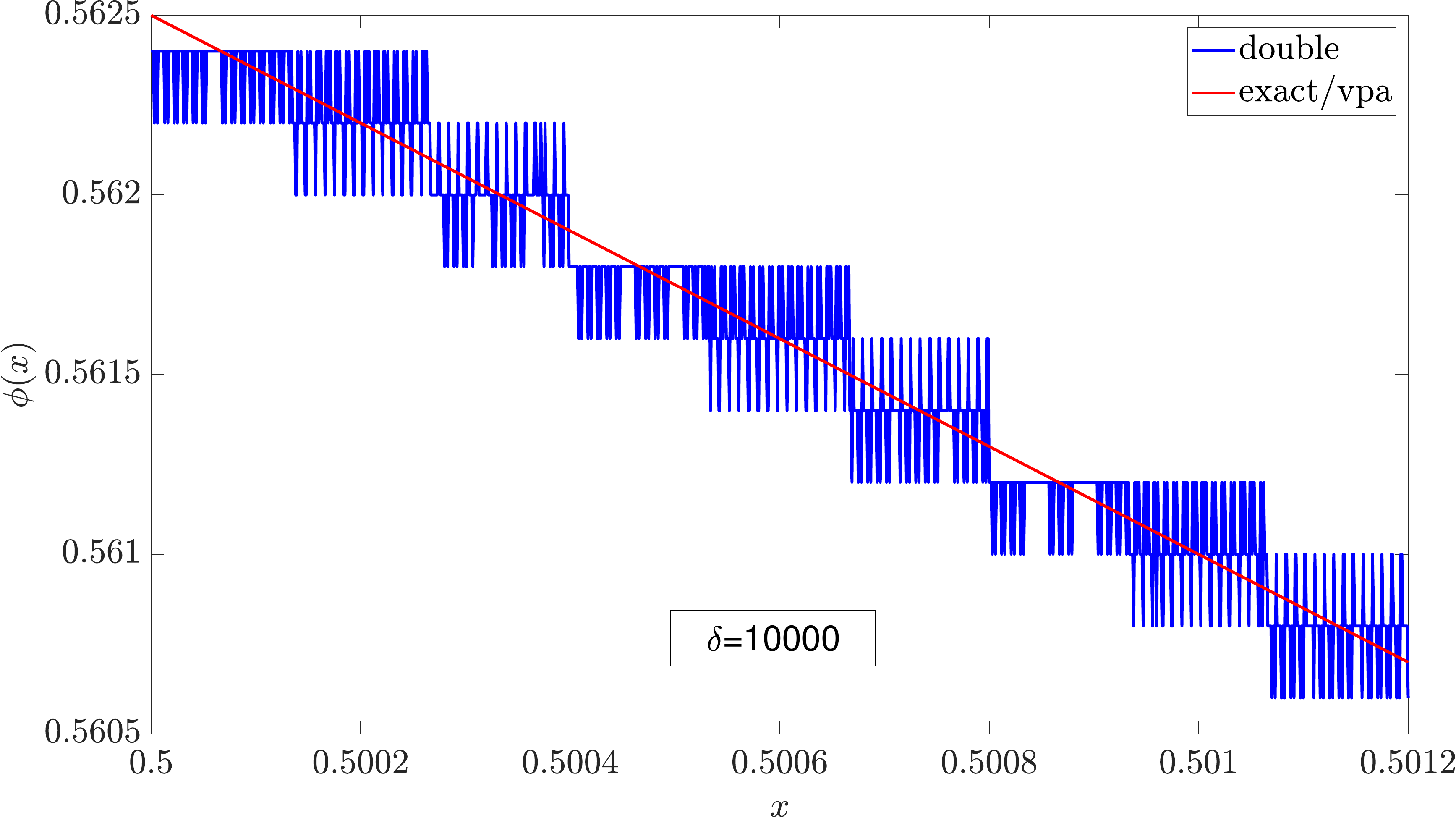}\\
\includegraphics[width=0.49\textwidth]{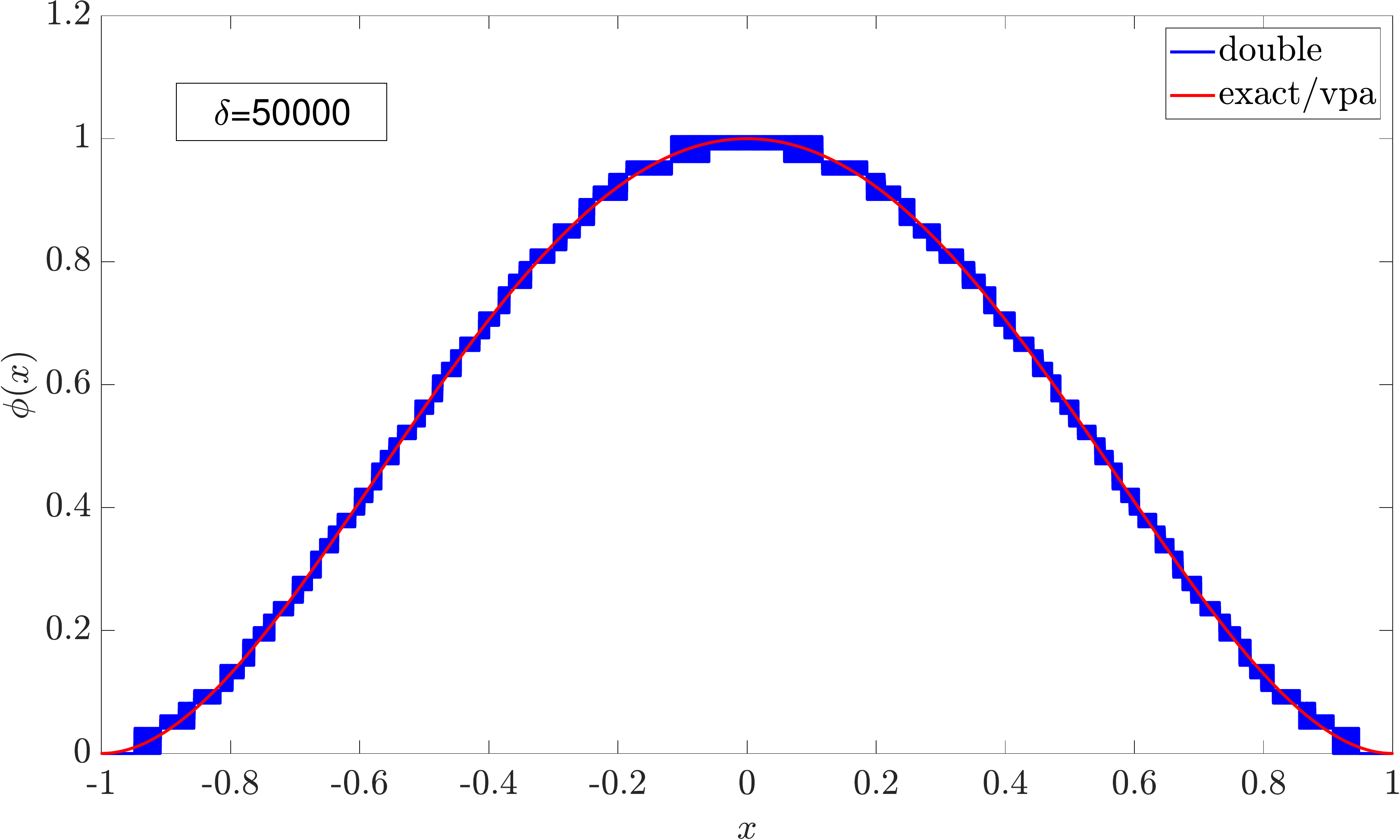}
\includegraphics[width=0.49\textwidth]{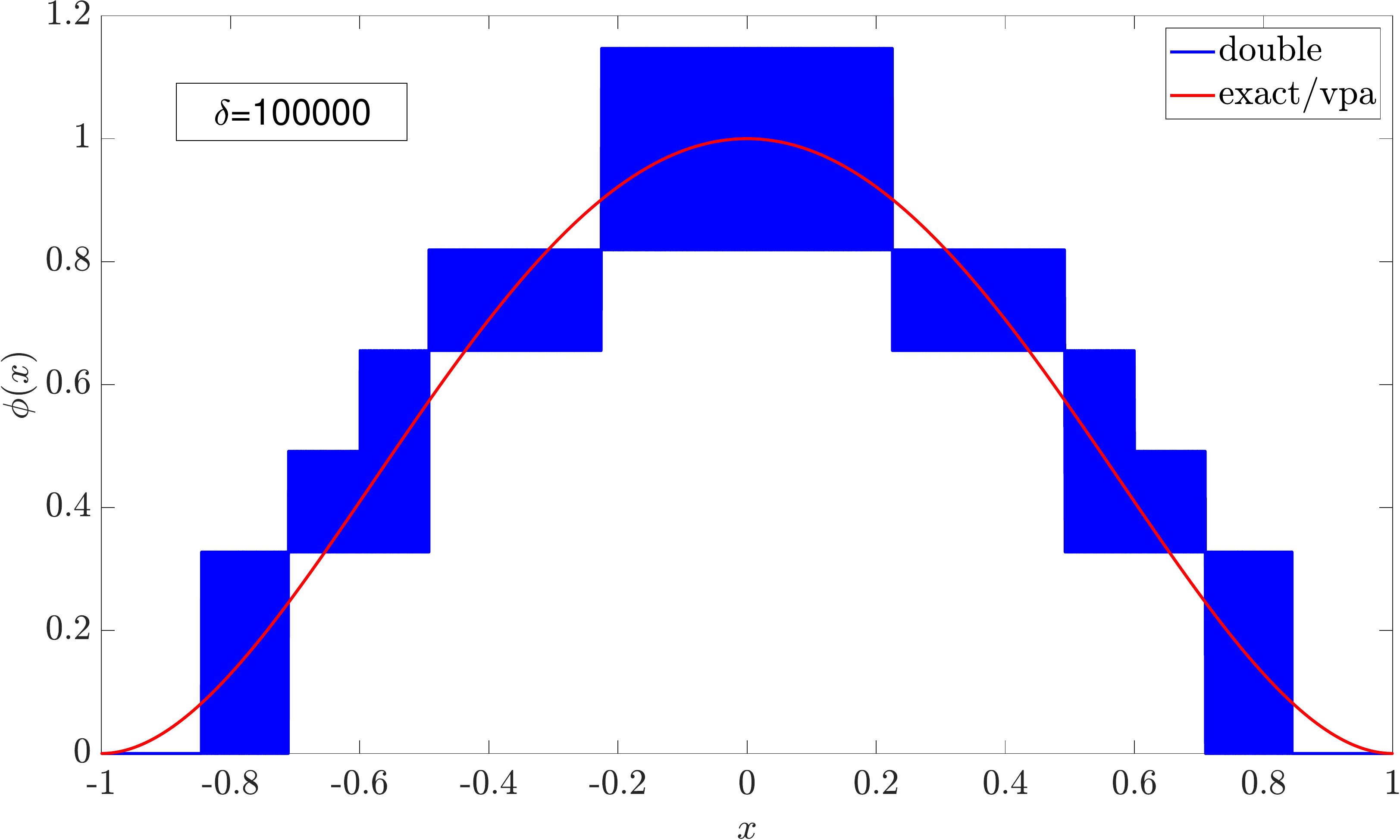}
\includegraphics[width=0.49\textwidth]{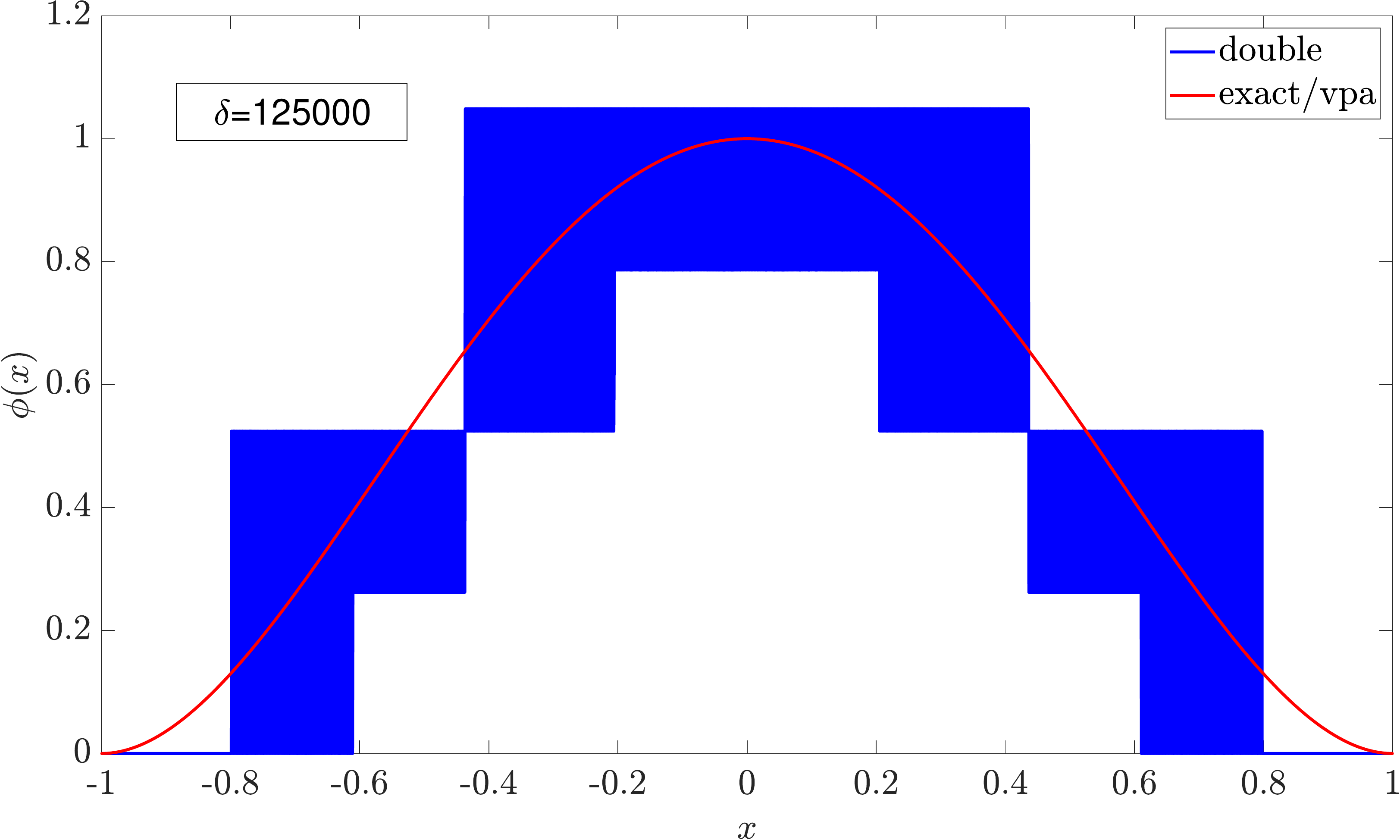}
\includegraphics[width=0.49\textwidth]{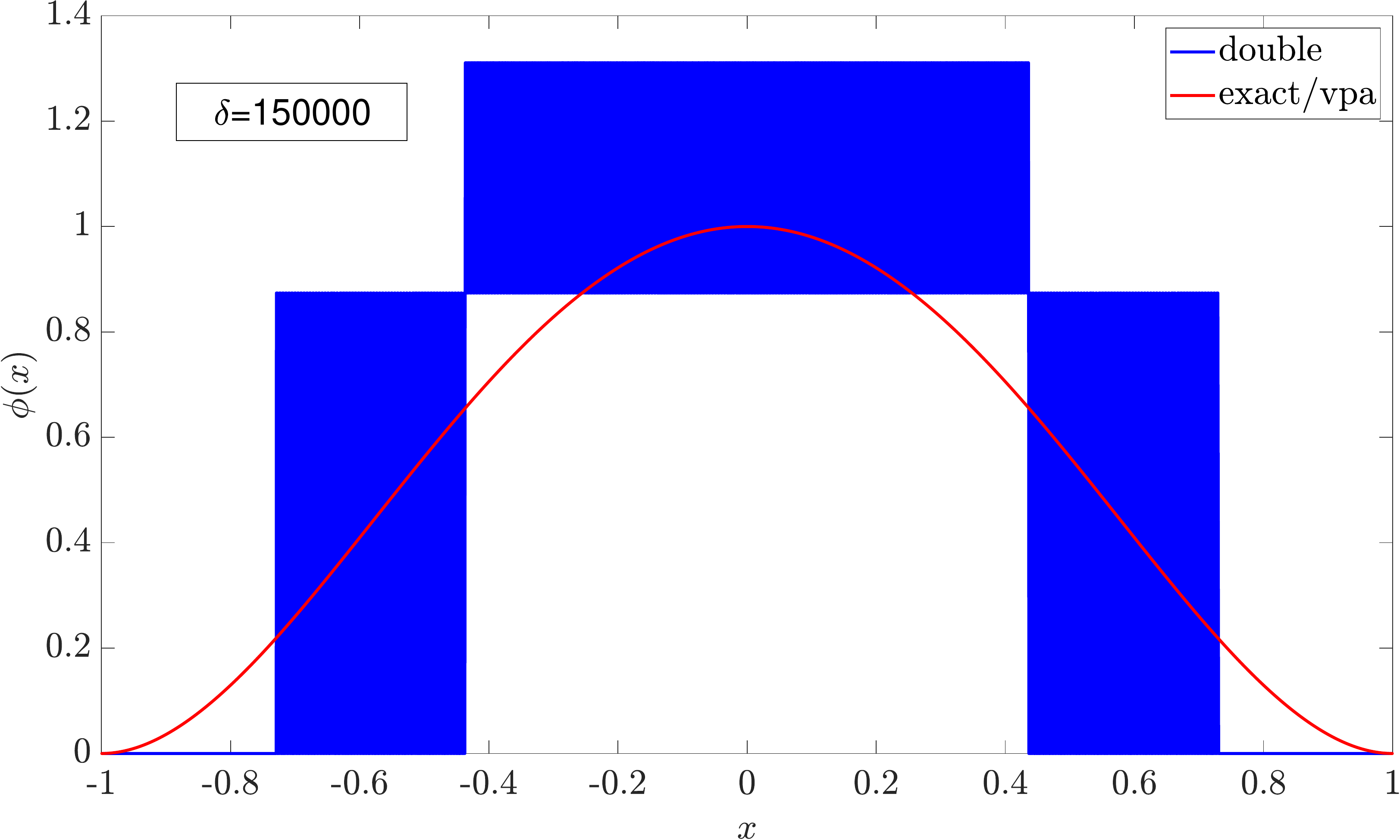}
\includegraphics[width=0.49\textwidth]{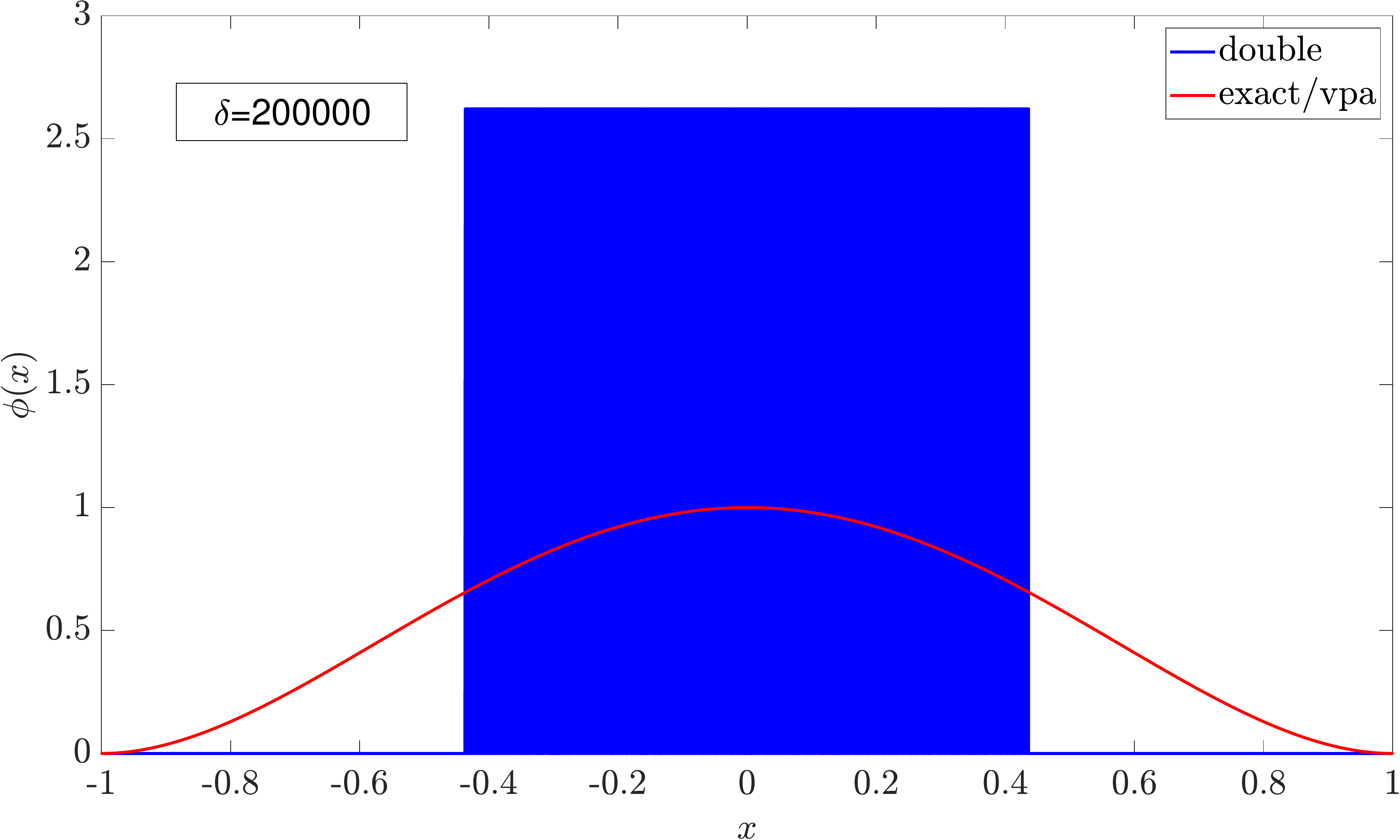}
\includegraphics[width=0.49\textwidth]{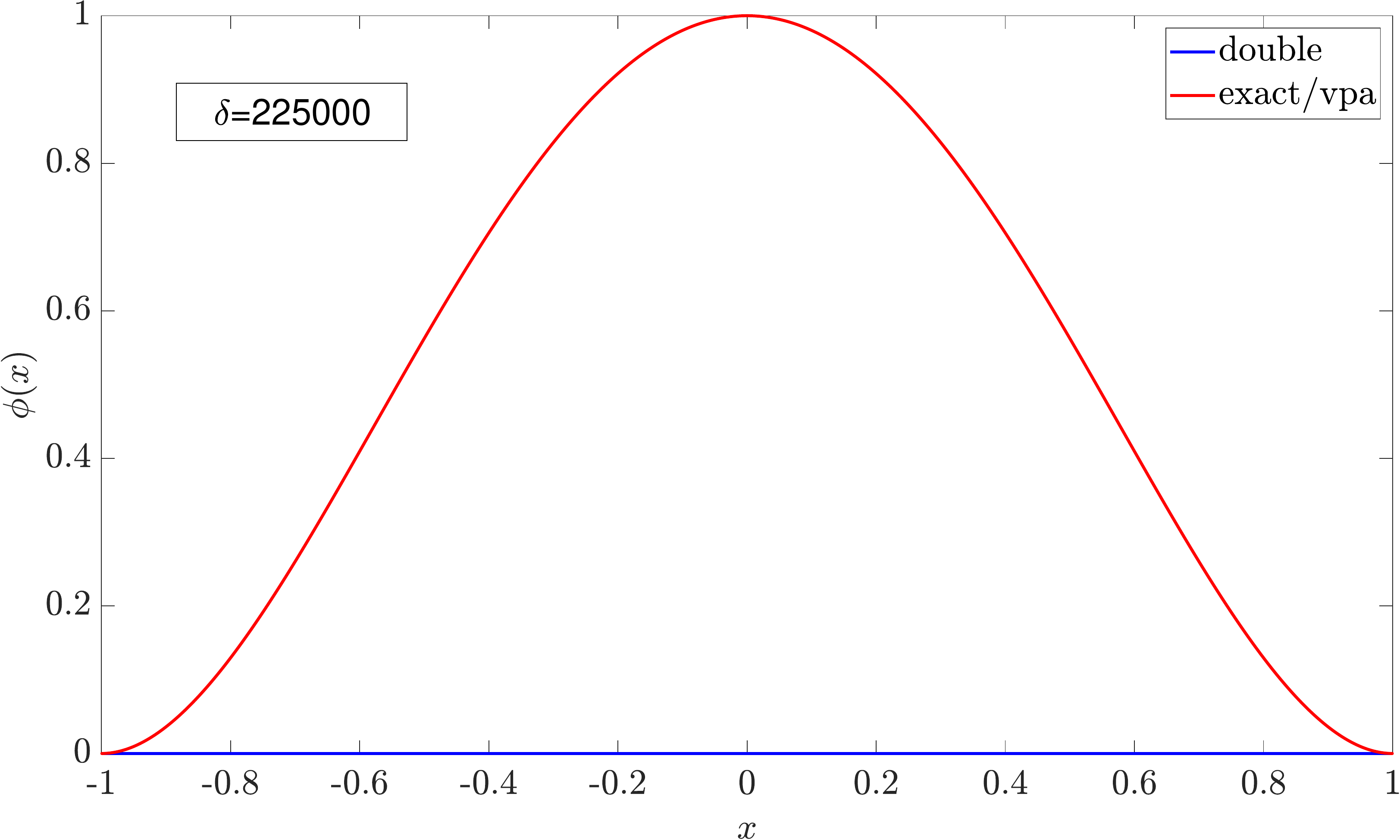}
\caption{Detailed views to the inaccurately evaluated function \texttt{phi\_double()} (\textcolor{blue}{blue}) in comparison to the accurate \texttt{phi\_vpa()} (\textcolor{red}{red}) for different values of $\delta$.}\label{f:phi_zoom}
\end{figure}

Figure~\ref{f:phi_zoom} demonstrates the inaccuracy in evaluating the function values for different values of the parameter $\delta$. Whereas a global view to the function $\phi(x)$ for values of delta less than ca $10000$ (top left picture) does not show any discrepancies, a detailed zoom (top right picture) indicates the differences in evaluating the function values in \texttt{double} (\textcolor{blue}{blue}) and \texttt{vpa} (\textcolor{red}{red}) arithmetic. It is worth to remind that the differences between \texttt{exact} and \texttt{vpa} (converted back to \texttt{double}) values are within the machine precision. Increasing values of $\delta$ further increases the inaccuracies of the \texttt{phi\_double()}. For $\delta$ equal to $225000$ (bottom right picture) or greater \texttt{phi\_double()} returns zero for all arguments.

\subsection{Numerical and automatic differentiation}\label{ss:diff} 

\begin{figure}[ht!]
\includegraphics[width=0.49\textwidth]{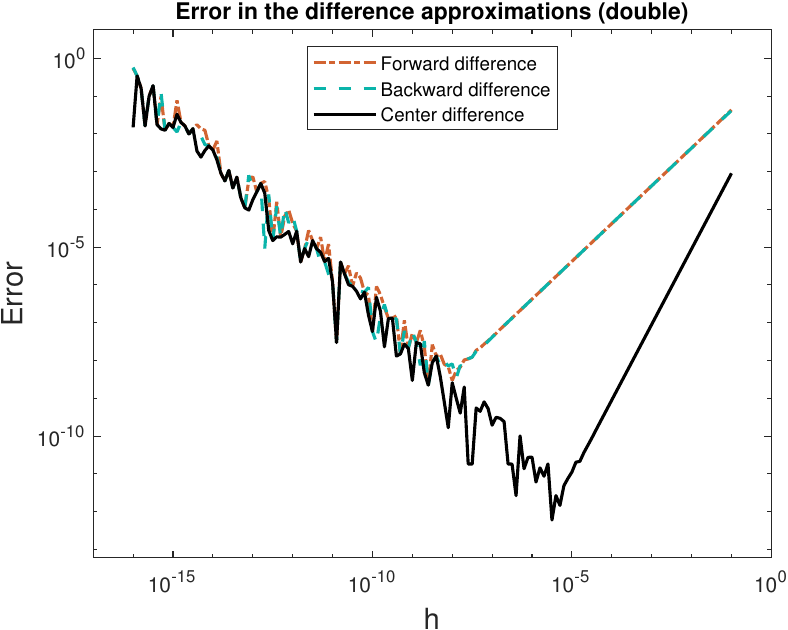}
\includegraphics[width=0.49\textwidth]{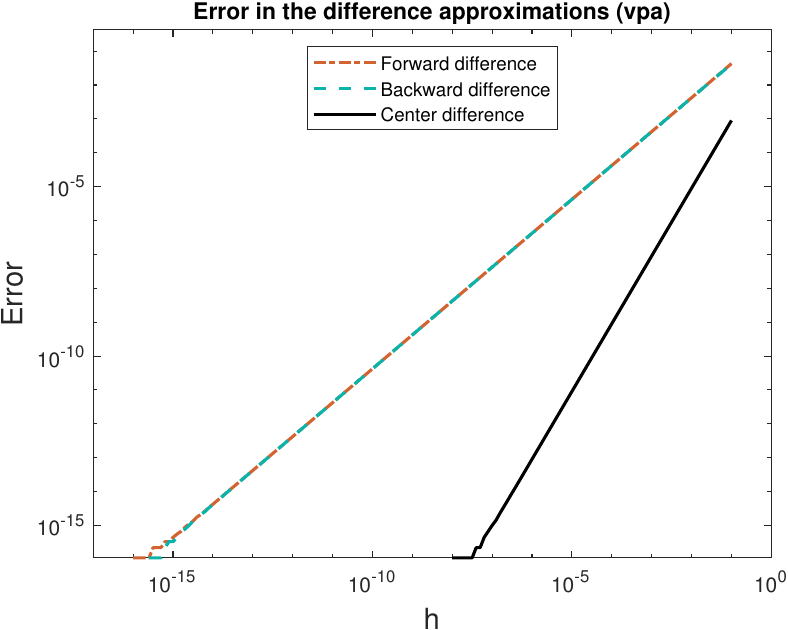}
\caption{Error in the difference approximation in \texttt{double} and \texttt{vpa} arithmetic.}\label{f:conditionality}
\end{figure}

Numerical approximations of derivatives are in general ill-conditioned and unstable especially due to the truncation error (often called error of the method, it is the error of approximation) and round-off error. In Figure~\ref{f:conditionality} we can see the error in the forward, backward and central difference approximating the derivative of function $\sin(x)$ at the point $x_0=1$ for different stepsize $h$. Standard \texttt{double} arithmetic is used on the left and \texttt{vpa} on the right. In both pictures on the right we can see that error of the method is dominant, the slope of the lines (there is a $\log-\log$ scale) corresponds to the order of the truncation error ($O(h)$ for the forward and backward difference schemes and $O(h^2)$ for the central difference). On the other hand, in \texttt{double} arithmetic we observe that with decreasing stepsize $h$ the round-off increases and becomes dominant.

Automatic differentiation (AD) is a computational technique for evaluating the derivatives of functions specified by computer programs, by systematically applying the chain rule to elementary operations and functions as the code executes. Unlike numerical differentiation, AD produces exact derivatives up to machine precision without step‑size tuning. And whereas symbolic differentiation manipulates mathematical expressions to obtain closed‑form derivative formulas, AD works directly on the program’s operations, neither requiring expression simplification nor generating unwieldy symbolic representations. By interleaving derivative computations with the original function evaluation, AD enables efficient gradient and Jacobian calculation for a wide class of functions. Although AD frameworks generate algebraically exact derivatives, they still execute those formulas in finite‑precision floating‑point \texttt{double} arithmetic. Therefore, the same numerical pathologies that occur in ordinary numerical codes re‑appear in the forward or reverse passes of AD. It is worth to mention that back propagation used for training neural networks is a special case of reverse accumulation pass \citep{Baydin2018AD}.

\subsection{Inverse Fourier transform integrand example}\label{ss:fk} %

Let us now consider a motivating example coming from mathematical finance that was for the first time studied by \cite{DanekPospisil20ijcm}. It can be shown \citep{BaustianMrazekPospisilSobotka17asmb} that in the so called approximate fractional stochastic volatility jump-diffusion (AFSVJD) model the price of a European call option with strike price $K$ and time to maturity $\tau$ can be written using the inverse complex Fourier transforms integral as

\begin{equation}
V~= S~- Ke^{-r\tau}\frac{1}{\pi}\int\limits_{0+i/2}^{+\infty+i/2} \underbrace{e^{-ikX}\frac{\hat{F}(k,v,\tau)}{k^2-ik}\phi(-k)}_{f(k)} dk,\label{e:price}
\end{equation}
with $X = \ln\frac{S}{K} + r\tau$, fundamental transform $\hat{F}(k,v,\tau) = \exp (C(k,\tau) + D(k,\tau) v)$ and
\begin{align*}
C(k,\tau) &= \kappa\theta \underbrace{\left(Y\tau - \frac{2}{B^2}\right)}_{C_1(k,\tau)} \underbrace{\ln\left(\frac{1-ge^{-d\tau}}{1-g}\right)}_{C_2(k,\tau)},\quad
D(k,\tau) = Y \frac{1-e^{-d\tau}}{1-ge^{-d\tau}},\quad Y = -\frac{k^2-ik}{b+d},\\
g &= \frac{b-d}{b+d},\quad d = \sqrt{b^2+B^2(k^2-ik)},\quad b = \kappa + ik\rho B,\quad B = \varepsilon^{H-1/2}\sigma, \\
\phi(k) &= \exp\left\{ -i \lambda \beta k~\tau + \lambda \tau \biggl[\hat{\varphi}(k)-1] \biggr] \right\}, \quad
\hat{\varphi}(k) = \exp\left\{ i\mu_J k~-\frac12 \sigma_J^2 k^2 \right\}, \\
\beta &= \exp\{ \mu_J + \sigma_J^2/2 \}-1.
\end{align*}
In this paper we omit all the mathematical finance details and focus strictly on the integrand $f(k)$ and its numerical integral later in Subsection \ref{ss:intfk}.

\begin{figure}[ht!]
\includegraphics[width=0.9\textwidth]{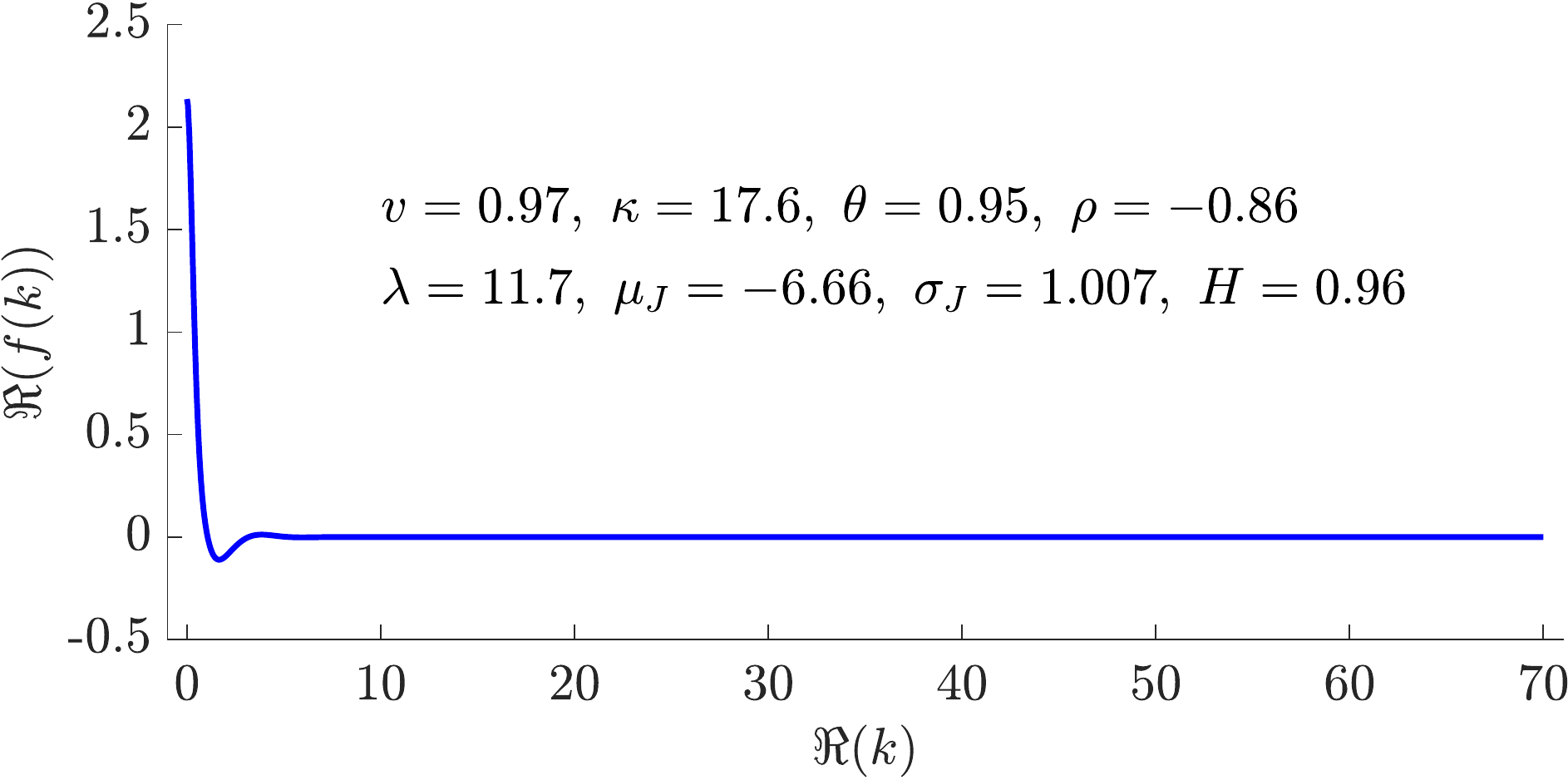}
\includegraphics[width=0.49\textwidth]{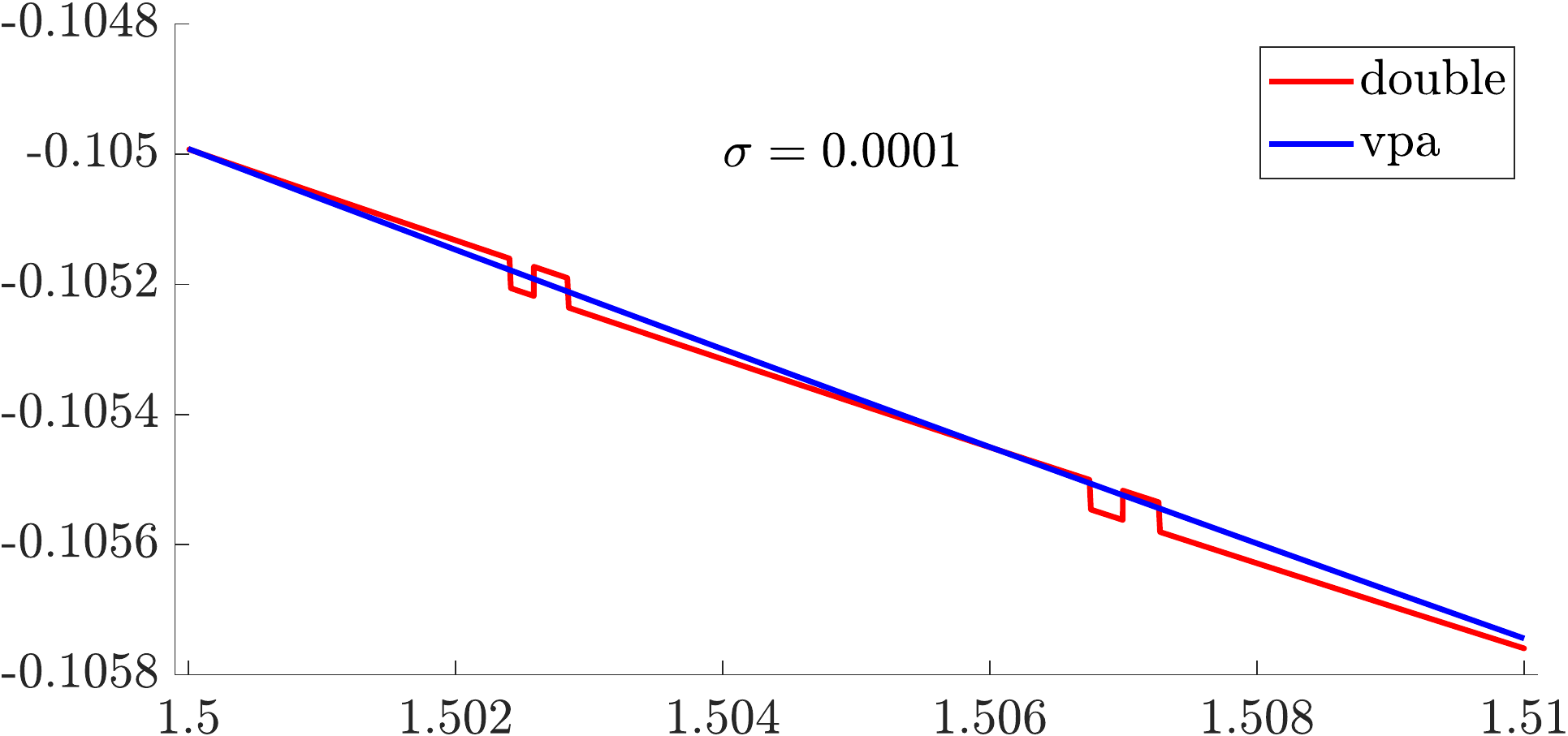}
\includegraphics[width=0.49\textwidth]{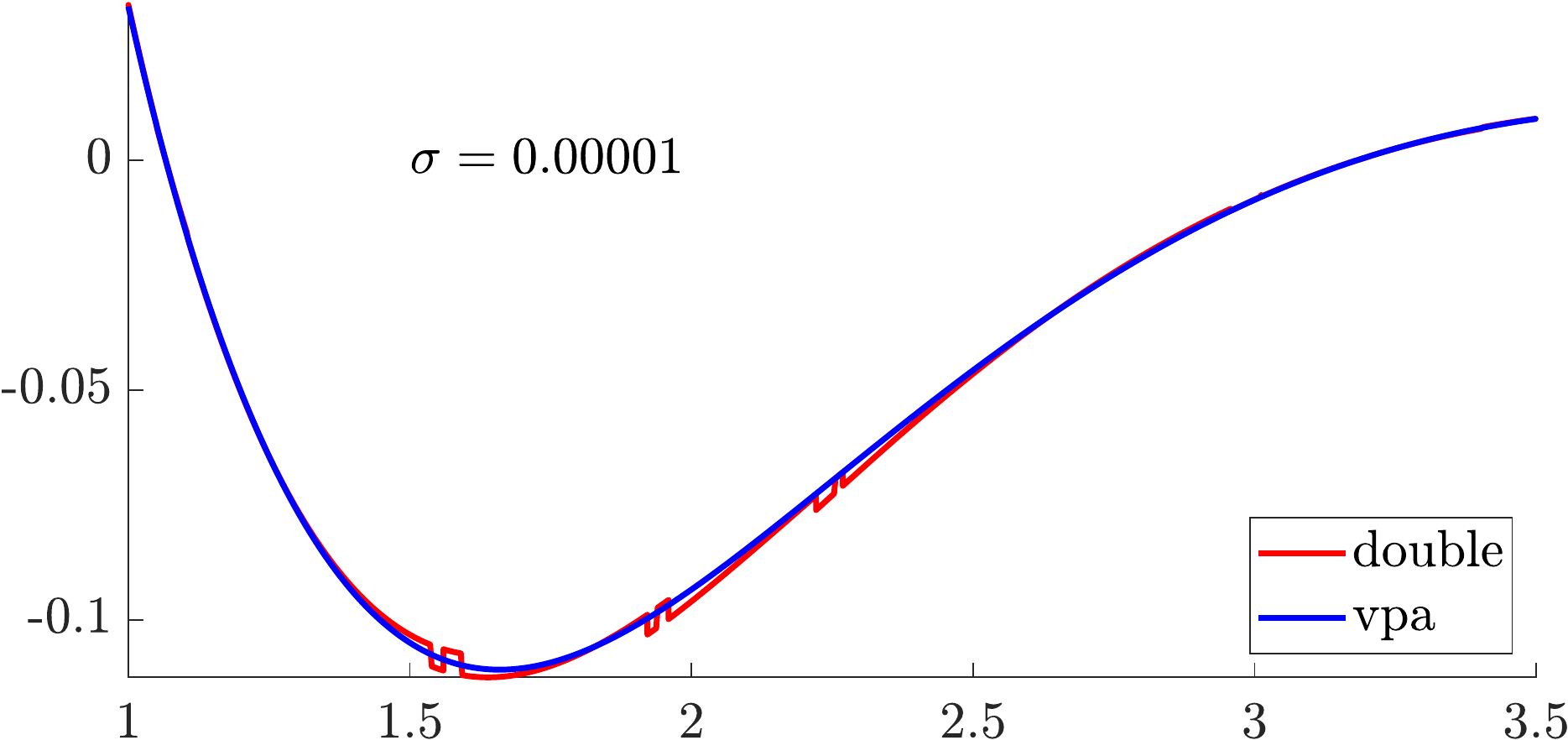}\\
\includegraphics[width=0.49\textwidth]{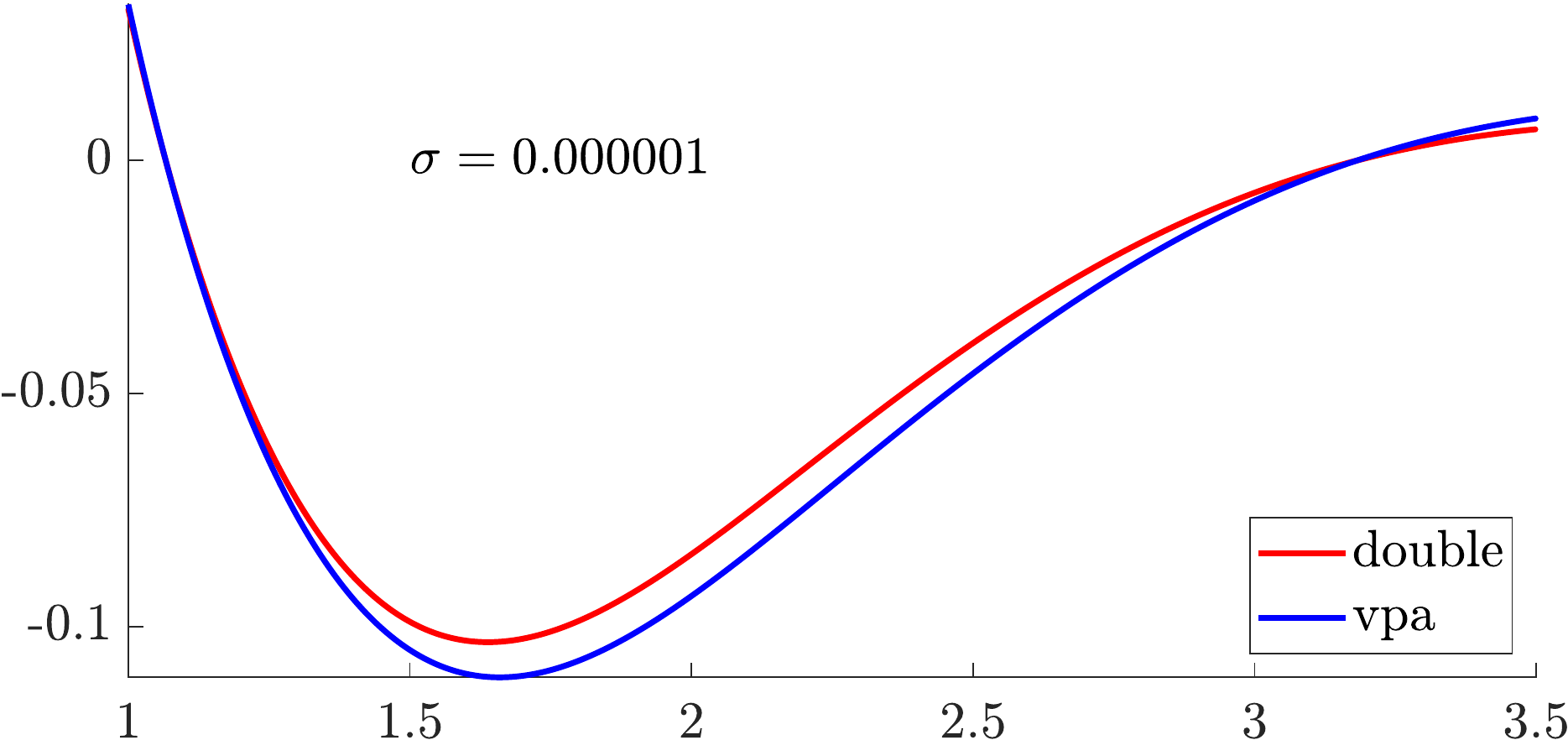}
\includegraphics[width=0.49\textwidth]{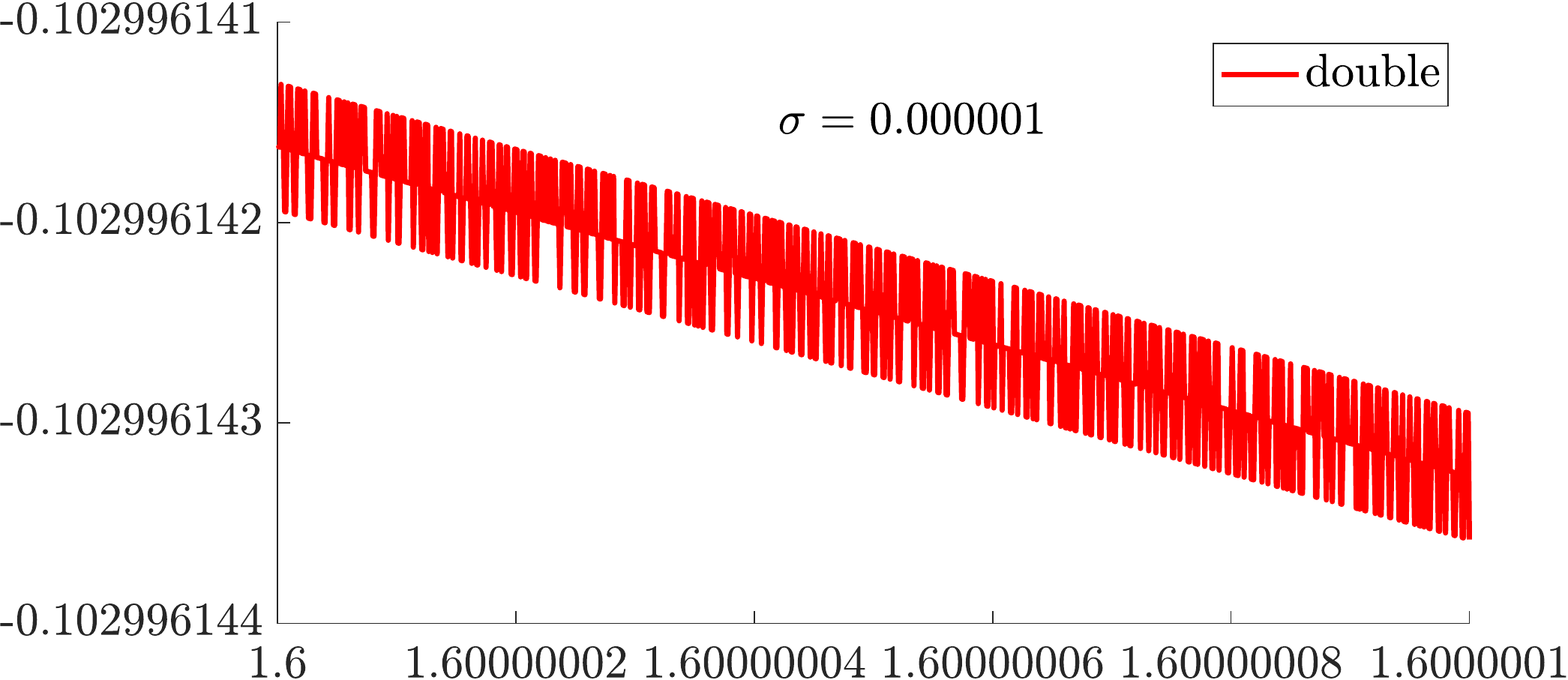}
\caption{Function $f(k)$ for some values of $\sigma$.}\label{f:finance}
\end{figure}

In Figure \ref{f:finance} we consider market data for European call options to the FTSE 100 index dated 8 January 2014, see \citep{PospisilSobotka16amf}, namely $\tau = 0.120548, K=6250, r=0.009, S=6721.8$ and model parameters as listed directly in the picture. Similarly as \cite{DanekPospisil20ijcm} we demonstrate the numerical misbehaviour by changing only the values of $\sigma$ and also the zoom depth. Whereas the inaccurately enumerated values are in \texttt{double} precision arithmetic (\textcolor{red}{red}), the smooth values are in \texttt{vpa} (\textcolor{blue}{blue}). In detailed zoom we can see that for greater values of $\sigma$, inaccuracies are more of local character, and for decreasing $\sigma$ we observe not only oscillating imprecision, but also a severe numerical bias, see bottom  pictures for $\sigma=0.000001$. Bottom right picture with even more detailed zoom of the \texttt{double} evaluated integrand reveals inaccurate numerical oscillations that are not visible on the less detailed zoom on the left.

 \clearpage
\section{Numerical quadratures and their failures}\label{sec3:quadratures}

The aim of this section is to show how numerical quadratures \alert{fail} to numerically evaluate the integral of functions introduced in the previous section within a prescribed precision. Moreover, if adaptive quadratures are used, we observe an enormous increase of function evaluations and hence the \alert{increase of computational time}. This increase could be one of the first signals to start doing the numerical analysis properly, both for the quadrature used as well as checking the integrand. 

\bigskip
Tested quadratures (existing functions in MATLAB):
\renewcommand*\descriptionlabel[1]{\hspace\labelsep\normalfont\textcolor{myblue}{\texttt{#1}}:}
\begin{description}
\item[integral] adaptive Gauss quadrature, in particular GaussKronrod(7,15),
\item[quad] adpative Simpson quadrature,
\item[quadl] adaptive Lobatto quadrature,
\item[trapz] trapezoidal numerical integration,
\item[tanh-sinh] tanh-sinh quadrature,
\item[vpaintegral] adaptive Gauss quadrature, again GaussKronrod(7,15), but in variable precision arithmetic, or more precisely in symbolic operations, i.e. the integrand must be convertible to a symbolic expression. 
\end{description}

Adaptive numerical quadratures either fail or in more recent versions of MATLAB at least warn you that the approximate bound on error is huge or that the maximum number of function evaluations was reached. In the tables below we will denote these two types of warnings using numbered superscripts as 

\smallskip\noindent\textcolor{red}{\texttt{%
Warning$^1$: Reached the limit on the maximum number of intervals in use. Approximate bound on error is $\dots$. The integral may not exist, or it may be difficult to approximate numerically to the requested accuracy. \\[3pt]
Warning$^2$: Maximum function count exceeded; singularity likely.
}}

\smallskip Failure of the adaptiveness of the general Gauss-Kronrod($n,2n+1$) quadrature on even a piece-wise constant function is described in Example 4.1 in \cite{DanekPospisil20ijcm} where authors suggest how to set the adaptive refinement in order to satisfy the local error estimate to be within the prescribed relative error, but even if this local error bound is zero, the value of the calculated integral does not have to be precise at all. 

Further, by \textcolor{myblue}{\texttt{integral\_vpa}} we denote the method in which the integrand (values of the integrated function) is evaluated in the variable precision arithmetic with the result being converted back to \texttt{double}. Similarly we denote \textcolor{myblue}{\texttt{quadl\_vpa}} and \textcolor{myblue}{\texttt{trapz\_vpa}}. 

For all quadratures we set the same error tolerances, absolute tolerance \texttt{AbsTol=1e-12} and relative tolerance \texttt{RelTol=1e-8} with only one exception, namely for \textcolor{myblue}{\texttt{vpaintegral}} default values \texttt{AbsTol=1e-10} and \texttt{RelTol=1e-6} were used.
 
Following the same structure as in the previous section, in Subsection \ref{ss:intfk} we will first numerically integrate the function $f(k)$ coming from mathematical finance (introduced in Subsection \ref{ss:fk}) and in Subsection \ref{ss:intphix} the numerical integral of function $\phi(x)$ (introduced in Subsection \ref{ss:phix}). Since this function $\phi(x)$ will serve as the main example later in the paper, in Subsection \ref{ss:intphixgnx} we also show how inaccuracy of the quadratures changes when the function is multiplied by the Fourier basis functions.

\subsection{Inverse Fourier transform integration}\label{ss:intfk} %

Let us now evaluate the integral in \eqref{e:price}, in particular
\begin{equation}\label{e:int_f}
I = \int\limits_{0+i/2}^{L+i/2} f(k)\d k,
\end{equation}
where $L=100$ is chosen for demonstration purposes only. Due to the complicated form of function $f(k)$ that depends on many parameters, see \eqref{e:price}, the exact result is not known. In Table~\ref{t:ex-finance} we can see comparison of different quadratures applied to \eqref{e:int_f} that is depicted in Figure~\ref{f:finance}. In the last column, Fevals denote the number of function (integrand) evaluations, i.e. the number of adaptive refinements. As we can see, for all values of $\sigma$, \textcolor{myblue}{\texttt{integral\_vpa}}, \textcolor{myblue}{\texttt{quadl\_vpa}} and \textcolor{myblue}{\texttt{trapz\_vpa}} gives results comparable within the machine precision and computational times. The value of \textcolor{myblue}{\texttt{integral\_vpa}} is therefore taken as the reference value and the value Error in the third column is the difference to this value. If \texttt{vpa} is not used, all quadratures are naturally much faster, but less precise. It is worth to emphasize that even if calculating the integrand in \texttt{vpa} is time consuming, the computational time can be comparable to the case when the integrand is calculated in \texttt{double} arithmetic only, which is caused by a massive increase of adaptive refinements caused by the inaccurate evaluation and hence the increase of Fevals that are in the order of hundreds of thousands. 

To understand the results of \textcolor{myblue}{\texttt{vpaintegral}} one has to realize that its argument is a function that is inaccurately evaluated in \texttt{double} arithmetic and the function tries to convert it to the symbolic expression that cannot be further simplified. Results obtained by \textcolor{myblue}{\texttt{vpaintegral}} in comparison to \textcolor{myblue}{\texttt{integral\_vpa}} are less precise and more time consuming, in particular precision worsens and computational time increases with decreasing $\sigma$. As we can see in Table~\ref{t:ex-finance} for $\sigma=0.00001$, the computational time of \textcolor{myblue}{\texttt{vpaintegral}} is more than four times higher than for \textcolor{myblue}{\texttt{integral\_vpa}} even for weaker tolerances that are described above.

\begin{table}[ht!]
\begin{center}
\caption{Comparison of different quadratures applied to \eqref{e:int_f} for some values of parameter $\sigma$.}\label{t:ex-finance}
\begin{tabular}{lcccc}
\hline 
Method & Value & Error & Time [s] & Fevals \\ 
\hline 
\hline 
{\it case} $\sigma=0.001$ &   &   &   &   \\ 
\hline 
\textcolor{myblue}{\texttt{integral}} & 0.7768149482534377 & 0.0000001607693588 & 2.825 & 271200 \\ 
\textcolor{myblue}{\texttt{integral\_vpa}} & 0.7768147874840788 & 0.0000000000000000 & 7.486 & 390 \\ 
\textcolor{myblue}{\texttt{vpaintegral}} & 0.7768147879892994 & 0.0000000005052205 & 0.389 & 1 \\ 
\textcolor{myblue}{\texttt{quad}} & 0.7769560723152419 & 0.0001412848311630 & 0.131 & 5010 \\ 
\textcolor{myblue}{\texttt{quadl}} & 0.9420883388070267 & 0.1652735513229479 & 0.082 & 2024 \\ 
\textcolor{myblue}{\texttt{quadl\_vpa}} & 0.7768147874840784 & -0.0000000000000004 & 6.979 & 236 \\ 
\textcolor{myblue}{\texttt{trapz-0.01}} & 0.7768149982921251 & 0.0000002108080462 & 0.014 & 1 \\ 
\textcolor{myblue}{\texttt{trapz-0.001}} & 0.7768149922690513 & 0.0000002047849724 & 0.025 & 1 \\ 
\textcolor{myblue}{\texttt{trapz\_vpa-0.01}} & 0.7768147874840787 & -0.0000000000000001 & 12.131 & 1 \\ 
\hline 
\hline 
{\it case} $\sigma=0.0001$ &   &   &   &   \\ 
\hline 
\textcolor{myblue}{\texttt{integral}} & 0.7768395484341176 & 0.0000247608438022 & 6.508 & 584700 \\ 
\textcolor{myblue}{\texttt{integral\_vpa}} & 0.7768147875903154 & 0.0000000000000000 & 8.523 & 390 \\ 
\textcolor{myblue}{\texttt{vpaintegral}} & 0.7768146694066477 & -0.0000001181836677 & 4.652 & 1 \\ 
\textcolor{myblue}{\texttt{quad}} & 0.7769532420467165 & 0.0001384544564010 & 0.146 & 5016 \\ 
\textcolor{myblue}{\texttt{quadl}} & 0.9352610191414621 & 0.1584462315511467 & 0.102 & 2030 \\ 
\textcolor{myblue}{\texttt{quadl\_vpa}} & 0.7768147875903146 & -0.0000000000000008 & 8.152 & 236 \\ 
\textcolor{myblue}{\texttt{trapz-0.01}} & 0.7768424929317990 & 0.0000277053414836 & 0.010 & 1 \\ 
\textcolor{myblue}{\texttt{trapz-0.001}} & 0.7768394964992882 & 0.0000247089089728 & 0.030 & 1 \\ 
\textcolor{myblue}{\texttt{trapz\_vpa-0.01}} & 0.7768147875903151 & -0.0000000000000003 & 12.740 & 1 \\ 
\hline 
\hline 
{\it case} $\sigma=0.00001$ &   &   &   &   \\ 
\hline 
\textcolor{myblue}{\texttt{integral}} & 0.7682354116124637 & -0.0085793759884747 & 1.460 & 94680 \\ 
\textcolor{myblue}{\texttt{integral\_vpa}} & 0.7768147876009384 & 0.0000000000000000 & 8.577 & 390 \\ 
\textcolor{myblue}{\texttt{vpaintegral}} & 0.7768968651904418 & 0.0000820775895034 & 36.579 & 1 \\ 
\textcolor{myblue}{\texttt{quad}} & 0.7682368656720007 & -0.0085779219289377 & 0.129 & 4392 \\ 
\textcolor{myblue}{\texttt{quadl}} & 0.7681158192130139 & -0.0086989683879245 & 0.113 & 2024 \\ 
\textcolor{myblue}{\texttt{quadl\_vpa}} & 0.7768147876009377 & -0.0000000000000007 & 8.133 & 236 \\ 
\textcolor{myblue}{\texttt{trapz-0.01}} & 0.7682814251749993 & -0.0085333624259391 & 0.013 & 1 \\ 
\textcolor{myblue}{\texttt{trapz-0.001}} & 0.7682369706788630 & -0.0085778169220754 & 0.030 & 1 \\ 
\textcolor{myblue}{\texttt{trapz\_vpa-0.01}} & 0.7768147876009375 & -0.0000000000000009 & 11.891 & 1 \\ 
\hline 
\end{tabular} \end{center}
\end{table}

\subsection{Integration of the problematic quartic function}\label{ss:intphix} %

The goal is to evaluate 

\begin{equation}\label{e:int_phi}
J = \int\limits_{-1}^{1} \alert{\phi(x)}\d x. 
\end{equation}

Although the exact result ($16/15\approx 1.066666666666667$) of the simplified \texttt{phi\_exact()} is known, many numerical quadratures \alert{fail} to numerically evaluate the integral of \texttt{phi\_double()} within a prescribed precision. In contrast to the previous example, when the further simplification of the integrand was not possible and the integration results were not satisfactory, here the automatic conversion of \texttt{phi\_double()} to symbolic expression performed by \textcolor{myblue}{\texttt{vpaintegral}} gives exactly the \texttt{phi\_exact()} and hence the exact integration result. Further we compare the integration results of quadratures \textcolor{myblue}{\texttt{integral}}, \textcolor{myblue}{\texttt{quad}}, \textcolor{myblue}{\texttt{quadl}}, \textcolor{myblue}{\texttt{trapz}} and \textcolor{myblue}{\texttt{tanh-sinh}} applied to integrands \texttt{phi\_double()}, \texttt{phi\_vpa()} and \texttt{phi\_exact()}. In Table \ref{t:int_phi} we can see the results of numerical integration for different quadratures and different values of $\delta$. As we can see all quadratures applied to \texttt{phi\_vpa()} or \texttt{phi\_exact()} gives the exact results, what differs is the computational time, \textcolor{myblue}{\texttt{integral}} and \textcolor{myblue}{\texttt{quad}} being the fastest. When the quadratures were applied to \texttt{phi\_double()} their results were problematic, worsening with increasing $\delta$. For higher $\delta$, the inaccuracy of the integrand evaluation is higher and hence the Fevals in the adaptive quadrature algorithm. In the last column we indicate the appearance of the above mentioned warnings. It is worth to mention that one has to be careful with low computational times that can sometimes occur for inaccurately evaluated integrand. The case $\delta=250000$ gives a completely wrong result $0.000000000000000$, because the integrand is already completely wrongly evaluated, see also Figure~\ref{f:phi_zoom}. Similarly as in the previous example, the massive increase of adaptive refinements is caused by the inaccurately evaluated integrand which leads to a huge increase of Fevals and consequently also the computational times, see for example the case $\delta=100000$ for \textcolor{myblue}{\texttt{quad}} and \textcolor{myblue}{\texttt{quadl}}.

\begin{table}[ht!]
\begin{center}
\caption{Comparison of different quadratures applied to \eqref{e:int_phi} for some values of parameter $\delta$.}\label{t:int_phi}
{\small%
\begin{tabular}{lclcccc}
\hline 
Integral Rule & $\delta$ & Integrand Rule & Value & Time [s] & Fevals & Warning \\ 
\hline 
\hline 
\textcolor{myblue}{\texttt{integral}} & 1000 & \textcolor{myblue}{\texttt{phi\_double}} & 1.066666662383775 & 0.180 & 12 & No \\ 
  &   & \textcolor{myblue}{\texttt{phi\_vpa}} & 1.066666666666666 & 0.069 & 1 & No \\ 
  &   & \textcolor{myblue}{\texttt{phi\_exact}} & 1.066666666666666 & 0.035 & 1 & No \\ 
\hline 
\textcolor{myblue}{\texttt{quad}} & 1000 & \textcolor{myblue}{\texttt{phi\_double}} & 1.066666641747667 & 0.068 & 14 & No \\ 
  &   & \textcolor{myblue}{\texttt{phi\_vpa}} & 1.066666666666667 & 0.113 & 14 & No \\ 
  &   & \textcolor{myblue}{\texttt{phi\_exact}} & 1.066666666666667 & 0.063 & 14 & No \\ 
\hline 
\textcolor{myblue}{\texttt{quadl}} & 1000 & \textcolor{myblue}{\texttt{phi\_double}} & 1.066666543018575 & 0.040 & 2 & No \\ 
  &   & \textcolor{myblue}{\texttt{phi\_vpa}} & 1.066666666666667 & 0.054 & 2 & No \\ 
  &   & \textcolor{myblue}{\texttt{phi\_exact}} & 1.066666666666667 & 0.040 & 2 & No \\ 
\hline 
\textcolor{myblue}{\texttt{trapz}} & 1000 & \textcolor{myblue}{\texttt{phi\_double}} & 1.066666666650391 & 0.038 & 1 & No \\ 
  &   & \textcolor{myblue}{\texttt{phi\_vpa}} & 1.066666666666667 & 2.346 & 1 & No \\ 
  &   & \textcolor{myblue}{\texttt{phi\_exact}} & 1.066666666666667 & 0.032 & 1 & No \\ 
\hline 
\textcolor{myblue}{\texttt{tanh-sinh}} & 1000 & \textcolor{myblue}{\texttt{phi\_double}} & 1.066666609903423 & 0.034 & 1 & No \\ 
  &   & \textcolor{myblue}{\texttt{phi\_vpa}} & 1.066666666666667 & 0.051 & 1 & No \\ 
  &   & \textcolor{myblue}{\texttt{phi\_exact}} & 1.066666666666667 & 0.060 & 1 & No \\ 
\hline 
\hline 
\textcolor{myblue}{\texttt{integral}} & 10000 & \textcolor{myblue}{\texttt{phi\_double}} & 1.066667217711548 & 0.084 & 11 & Yes$^1$ \\ 
  &   & \textcolor{myblue}{\texttt{phi\_vpa}} & 1.066666666666666 & 0.073 & 1 & No \\ 
  &   & \textcolor{myblue}{\texttt{phi\_exact}} & 1.066666666666666 & 0.046 & 1 & No \\ 
\hline 
\textcolor{myblue}{\texttt{quad}} & 10000 & \textcolor{myblue}{\texttt{phi\_double}} & 1.066721252633333 & 0.083 & 22 & No \\ 
  &   & \textcolor{myblue}{\texttt{phi\_vpa}} & 1.066666666666666 & 0.176 & 14 & No \\ 
  &   & \textcolor{myblue}{\texttt{phi\_exact}} & 1.066666666666667 & 0.073 & 14 & No \\ 
\hline 
\textcolor{myblue}{\texttt{quadl}} & 10000 & \textcolor{myblue}{\texttt{phi\_double}} & 1.066649867548471 & 0.350 & 128 & No \\ 
  &   & \textcolor{myblue}{\texttt{phi\_vpa}} & 1.066666666666667 & 0.057 & 2 & No \\ 
  &   & \textcolor{myblue}{\texttt{phi\_exact}} & 1.066666666666667 & 0.041 & 2 & No \\ 
\hline 
\textcolor{myblue}{\texttt{trapz}} & 10000 & \textcolor{myblue}{\texttt{phi\_double}} & 1.066667380000000 & 0.032 & 1 & No \\ 
  &   & \textcolor{myblue}{\texttt{phi\_vpa}} & 1.066666666666667 & 2.544 & 1 & No \\ 
  &   & \textcolor{myblue}{\texttt{phi\_exact}} & 1.066666666666667 & 0.033 & 1 & No \\ 
\hline 
\textcolor{myblue}{\texttt{tanh-sinh}} & 10000 & \textcolor{myblue}{\texttt{phi\_double}} & 1.066565320860938 & 0.031 & 1 & No \\ 
  &   & \textcolor{myblue}{\texttt{phi\_vpa}} & 1.066666666666667 & 0.060 & 1 & No \\ 
  &   & \textcolor{myblue}{\texttt{phi\_exact}} & 1.066666666666667 & 0.038 & 1 & No \\ 
\hline 
\hline 
\textcolor{myblue}{\texttt{integral}} & 100000 & \textcolor{myblue}{\texttt{phi\_double}} & 1.073588382586769 & 0.125 & 12 & Yes$^1$ \\ 
  &   & \textcolor{myblue}{\texttt{phi\_vpa}} & 1.066666666666667 & 0.141 & 1 & No \\ 
  &   & \textcolor{myblue}{\texttt{phi\_exact}} & 1.066666666666666 & 0.091 & 1 & No \\ 
\hline 
\textcolor{myblue}{\texttt{quad}} & 100000 & \textcolor{myblue}{\texttt{phi\_double}} & 1.075099781351111 & 16.240 & 5006 & Yes$^2$ \\ 
  &   & \textcolor{myblue}{\texttt{phi\_vpa}} & 1.066666666666667 & 0.117 & 14 & No \\ 
  &   & \textcolor{myblue}{\texttt{phi\_exact}} & 1.066666666666667 & 0.072 & 14 & No \\ 
\hline 
\textcolor{myblue}{\texttt{quadl}} & 100000 & \textcolor{myblue}{\texttt{phi\_double}} & 1.075775476723900 & 4.734 & 2024 & Yes$^2$ \\ 
  &   & \textcolor{myblue}{\texttt{phi\_vpa}} & 1.066666666666667 & 0.062 & 2 & No \\ 
  &   & \textcolor{myblue}{\texttt{phi\_exact}} & 1.066666666666667 & 0.038 & 2 & No \\ 
\hline 
\textcolor{myblue}{\texttt{trapz}} & 100000 & \textcolor{myblue}{\texttt{phi\_double}} & 1.071808512000003 & 0.034 & 1 & No \\ 
  &   & \textcolor{myblue}{\texttt{phi\_vpa}} & 1.066666666666667 & 2.355 & 1 & No \\ 
  &   & \textcolor{myblue}{\texttt{phi\_exact}} & 1.066666666666667 & 0.062 & 1 & No \\ 
\hline 
\textcolor{myblue}{\texttt{tanh-sinh}} & 100000 & \textcolor{myblue}{\texttt{phi\_double}} & 1.062458555214656 & 0.037 & 1 & No \\ 
  &   & \textcolor{myblue}{\texttt{phi\_vpa}} & 1.066666666666667 & 0.053 & 1 & No \\ 
  &   & \textcolor{myblue}{\texttt{phi\_exact}} & 1.066666666666667 & 0.032 & 1 & No \\ 
\hline 
\hline 
\textcolor{myblue}{\texttt{integral}} & 250000 & \textcolor{myblue}{\texttt{phi\_double}} & 0.000000000000000 & 0.046 & 1 & No \\ 
  &   & \textcolor{myblue}{\texttt{phi\_vpa}} & 1.066666666666666 & 0.076 & 1 & No \\ 
  &   & \textcolor{myblue}{\texttt{phi\_exact}} & 1.066666666666666 & 0.045 & 1 & No \\ 
\hline 
\textcolor{myblue}{\texttt{quad}} & 250000 & \textcolor{myblue}{\texttt{phi\_double}} & 0.000000000000000 & 0.044 & 4 & No \\ 
  &   & \textcolor{myblue}{\texttt{phi\_vpa}} & 1.066666666666667 & 0.118 & 14 & No \\ 
  &   & \textcolor{myblue}{\texttt{phi\_exact}} & 1.066666666666667 & 0.066 & 14 & No \\ 
\hline 
\textcolor{myblue}{\texttt{quadl}} & 250000 & \textcolor{myblue}{\texttt{phi\_double}} & 0.000000000000000 & 0.041 & 2 & No \\ 
  &   & \textcolor{myblue}{\texttt{phi\_vpa}} & 1.066666666666667 & 0.059 & 2 & No \\ 
  &   & \textcolor{myblue}{\texttt{phi\_exact}} & 1.066666666666667 & 0.038 & 2 & No \\ 
\hline 
\textcolor{myblue}{\texttt{trapz}} & 250000 & \textcolor{myblue}{\texttt{phi\_double}} & 0.000000000000000 & 0.037 & 1 & No \\ 
  &   & \textcolor{myblue}{\texttt{phi\_vpa}} & 1.066666666666667 & 2.374 & 1 & No \\ 
  &   & \textcolor{myblue}{\texttt{phi\_exact}} & 1.066666666666667 & 0.032 & 1 & No \\ 
\hline 
\textcolor{myblue}{\texttt{tanh-sinh}} & 250000 & \textcolor{myblue}{\texttt{phi\_double}} & 0.000000000000000 & 0.031 & 1 & No \\ 
  &   & \textcolor{myblue}{\texttt{phi\_vpa}} & 1.066666666666667 & 0.052 & 1 & No \\ 
  &   & \textcolor{myblue}{\texttt{phi\_exact}} & 1.066666666666667 & 0.035 & 1 & No \\ 
\hline 
\hline 
\end{tabular} %
}
\end{center}
\end{table}
\subsection{Calculation of Fourier coefficients}\label{ss:intphixgnx} %

Let us now consider basis functions 
\begin{equation}\label{e:g_n}
g_n(x) = \sin\left(\pi(2n-1)(x+1)/2\right),\quad n\in\N, 
\end{equation}
that are even, orthonormal at $[-1,1]$ and satisfy $g_n(-1)=g_n(1)=0$. Alternatively, an equivalent formula 
$$g_n(x) = (-1)^{n+1}\cos\left(\pi(2n-1)x/2\right),\quad n\in\N,$$
can be considered.

The goal now is to evaluate the Fourier coefficients 
\begin{equation}\label{e:A_n}
A_n = \int\limits_{-1}^{1} \alert{\phi(x)}g_n(x)\d x,\quad n\in\N. 
\end{equation}

\begin{figure}[!ht]
\includegraphics[width=0.49\textwidth]{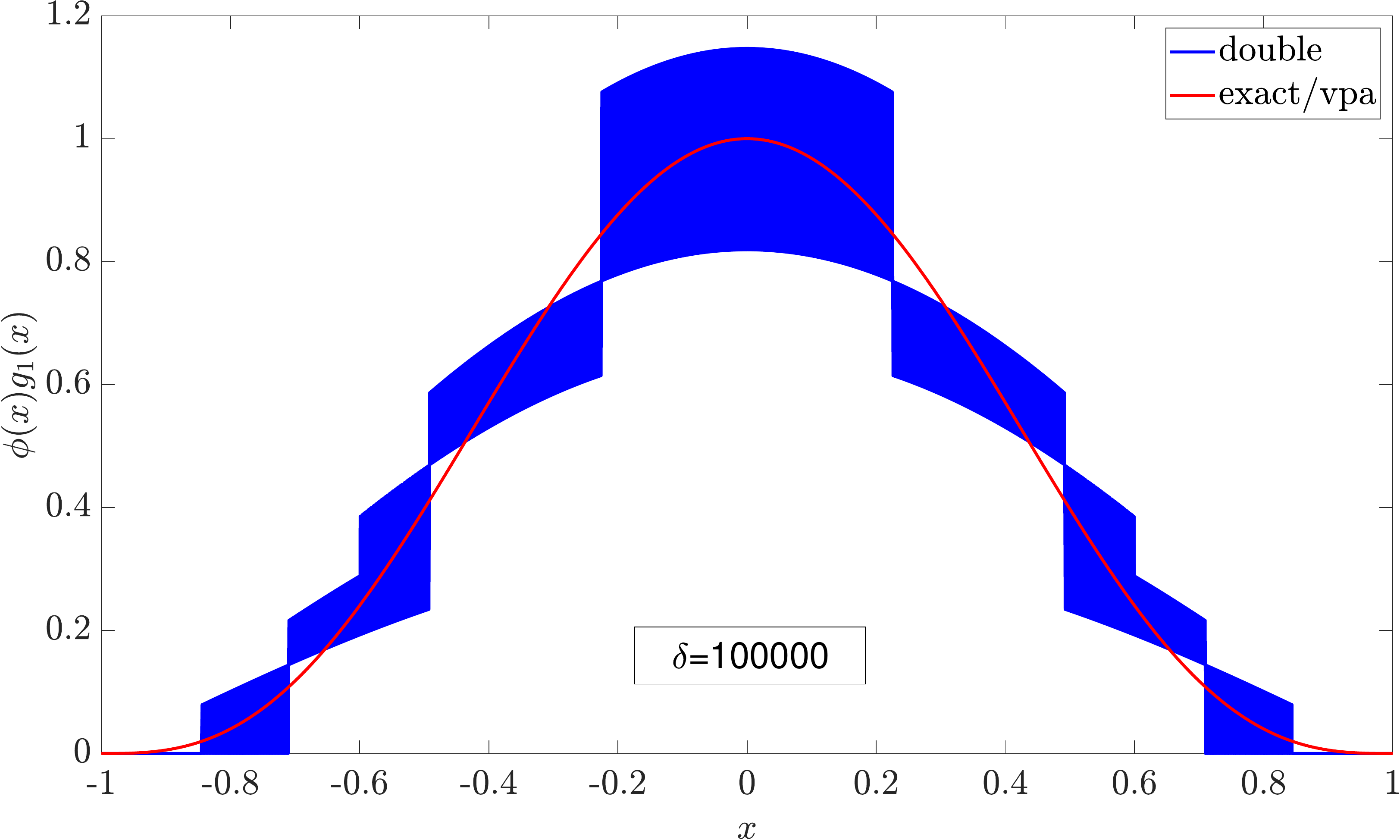}
\includegraphics[width=0.49\textwidth]{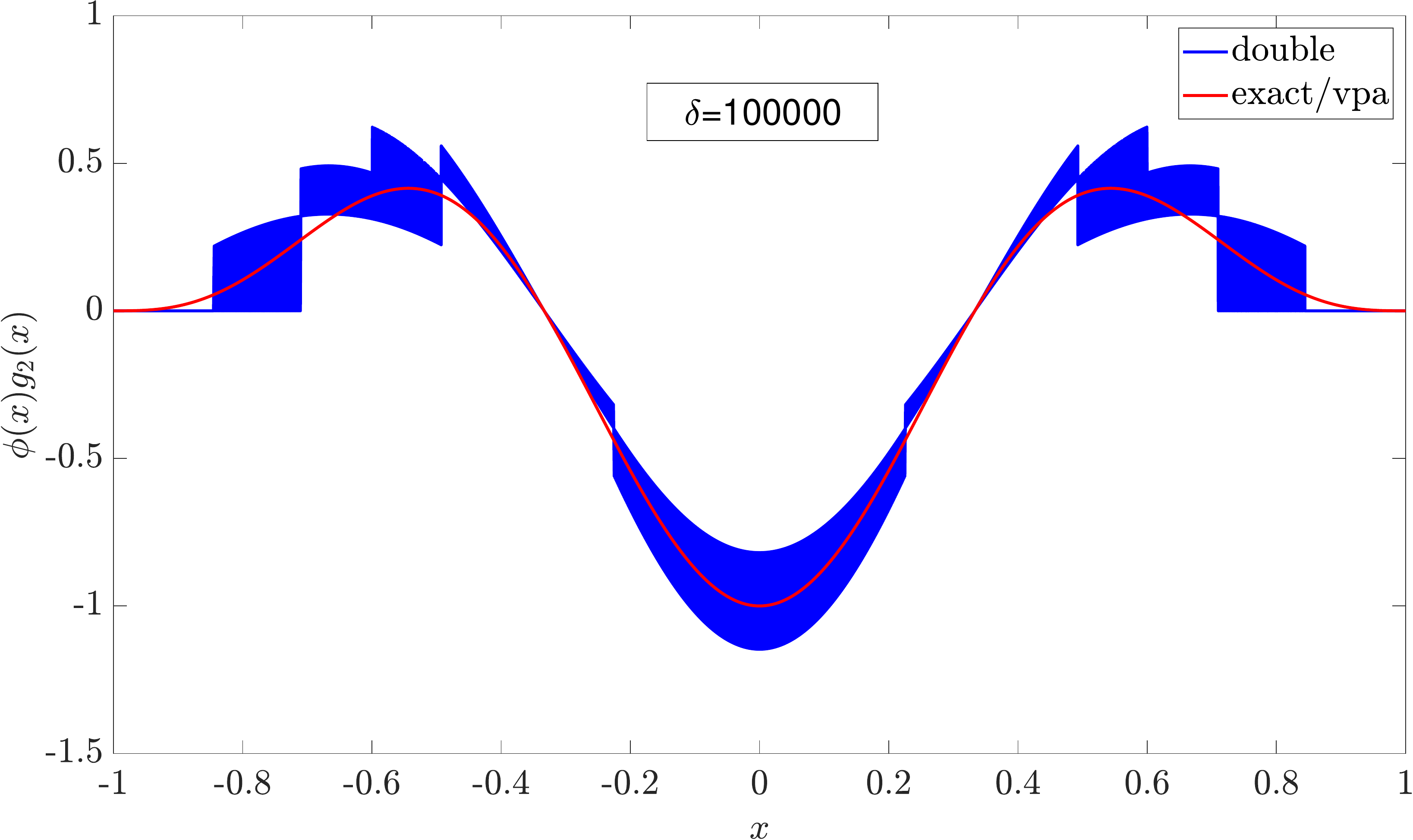}\\
\includegraphics[width=0.49\textwidth]{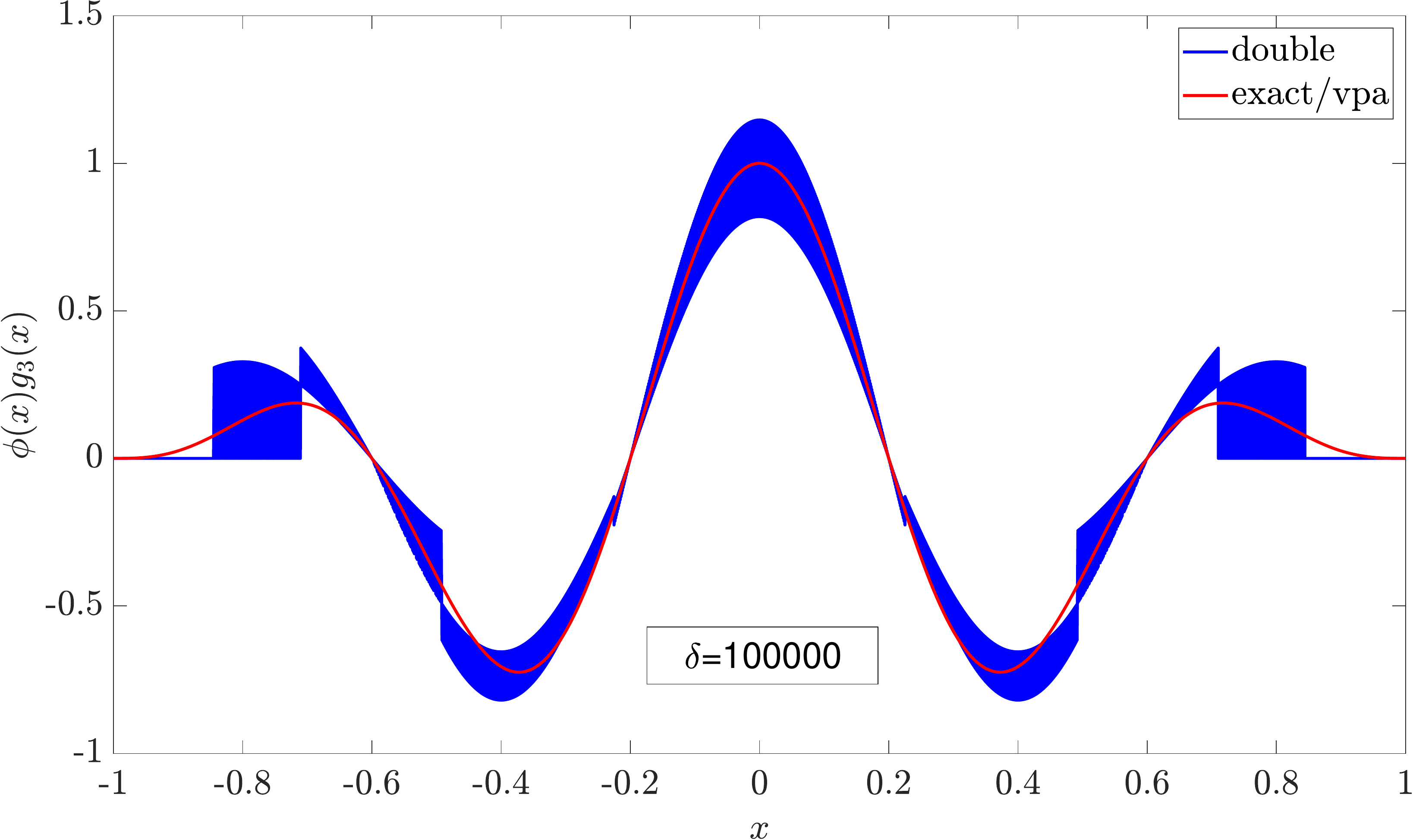}
\includegraphics[width=0.49\textwidth]{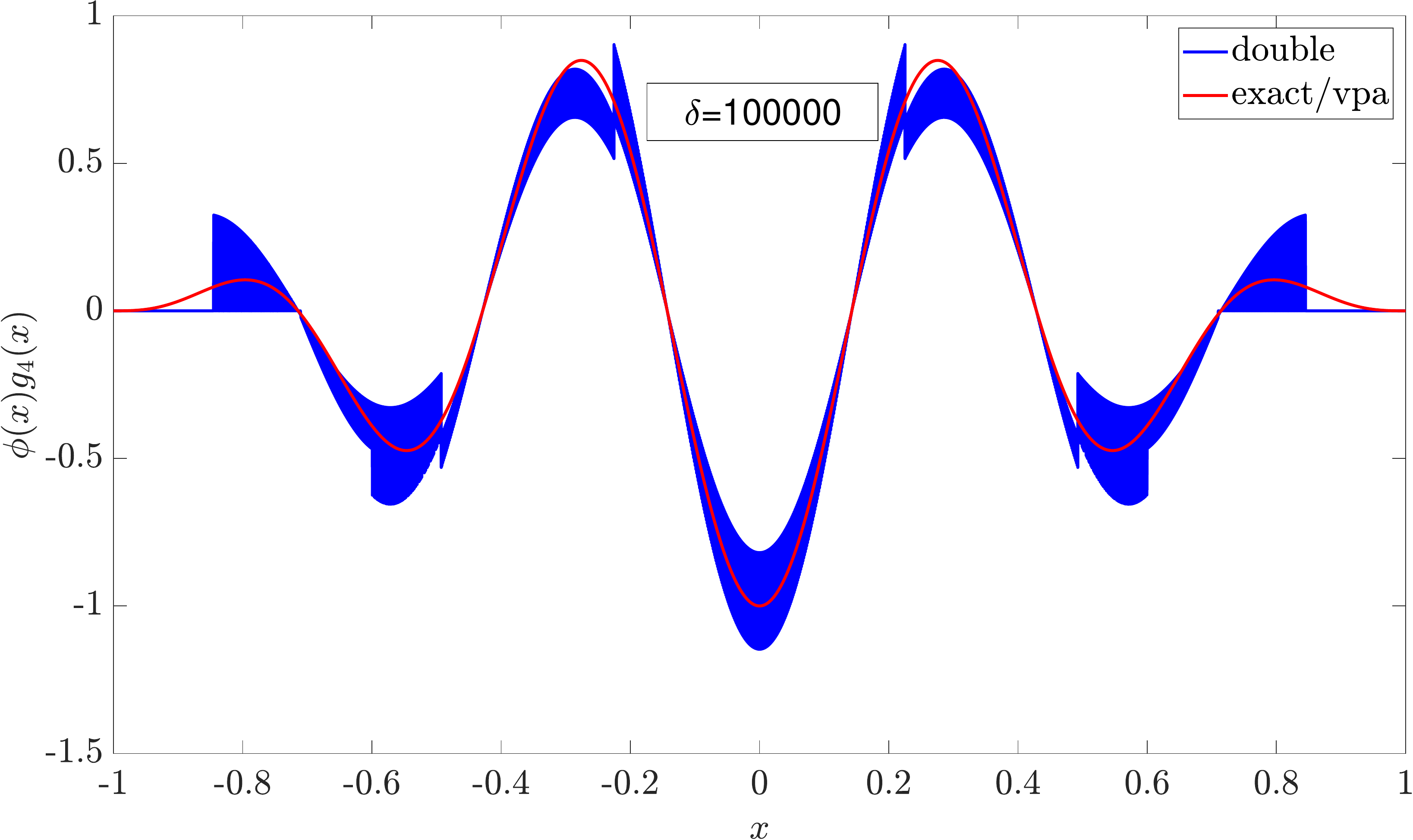}
\caption{Graph of inaccurately evaluated integrand $\phi(x)g_n(x)$ for $\delta=100000$ and different values of $n$.}\label{f:phi_g_n}
\end{figure}

In Figure \ref{f:phi_g_n} we can see how the integrand changes for different values of $n$. The differences in evaluating the integrand values in \texttt{double} (\textcolor{blue}{blue}) and \texttt{vpa} (\textcolor{red}{red}) arithmetic for $\delta=100000$ are visually obvious.

In Tables \ref{t:A1} and \ref{t:A2} we list results for $A_1$ and $A_2$ \eqref{e:A_n} respectively. We can observe a similar behaviour as in Table \ref{t:int_phi} where the integrand did not contain the basis function. As we can see, basis functions do not bring additional difficulties. Due to the boundedness of the basis functions, the numerical difference between the \texttt{double} and \texttt{vpa} can be smaller, but the product $\alert{\phi(x)}g_n(x)$ posses the same bad numerical misbehaviour as the inaccurately evaluated $\alert{\phi(x)}$.

\begin{table}[ht!]
\begin{center}
\caption{Comparison of different quadratures used to calculate the value of Fourier coefficient $A_1$}\label{t:A1}
{\small%
\begin{tabular}{lclcccc}
\hline 
Integral Rule & $\delta$ & Integrand Rule & Value & Time [s] & Fevals & Warning \\ 
\hline 
\hline 
\textcolor{myblue}{\texttt{integral}} & 1000 & \textcolor{myblue}{\texttt{phi\_double}} & 0.891088555909676 & 0.107 & 11 & No \\ 
  &   & \textcolor{myblue}{\texttt{phi\_vpa}} & 0.891088548280465 & 0.076 & 1 & No \\ 
  &   & \textcolor{myblue}{\texttt{phi\_exact}} & 0.891088548280465 & 0.042 & 1 & No \\ 
\hline 
\textcolor{myblue}{\texttt{quad}} & 1000 & \textcolor{myblue}{\texttt{phi\_double}} & 0.891088497788002 & 0.098 & 26 & No \\ 
  &   & \textcolor{myblue}{\texttt{phi\_vpa}} & 0.891088530616134 & 0.182 & 26 & No \\ 
  &   & \textcolor{myblue}{\texttt{phi\_exact}} & 0.891088530616134 & 0.098 & 26 & No \\ 
\hline 
\textcolor{myblue}{\texttt{quadl}} & 1000 & \textcolor{myblue}{\texttt{phi\_double}} & 0.891088558410162 & 0.065 & 8 & No \\ 
  &   & \textcolor{myblue}{\texttt{phi\_vpa}} & 0.891088548278784 & 0.092 & 8 & No \\ 
  &   & \textcolor{myblue}{\texttt{phi\_exact}} & 0.891088548278784 & 0.057 & 8 & No \\ 
\hline 
\textcolor{myblue}{\texttt{trapz}} & 1000 & \textcolor{myblue}{\texttt{phi\_double}} & 0.891088548289315 & 0.036 & 1 & No \\ 
  &   & \textcolor{myblue}{\texttt{phi\_vpa}} & 0.891088548280466 & 2.369 & 1 & No \\ 
  &   & \textcolor{myblue}{\texttt{phi\_exact}} & 0.891088548280466 & 0.033 & 1 & No \\ 
\hline 
\textcolor{myblue}{\texttt{tanh-sinh}} & 1000 & \textcolor{myblue}{\texttt{phi\_double}} & 0.891088505867673 & 0.032 & 1 & No \\ 
  &   & \textcolor{myblue}{\texttt{phi\_vpa}} & 0.891088548280696 & 0.053 & 1 & No \\ 
  &   & \textcolor{myblue}{\texttt{phi\_exact}} & 0.891088548280696 & 0.059 & 1 & No \\ 
\hline 
\hline 
\textcolor{myblue}{\texttt{integral}} & 10000 & \textcolor{myblue}{\texttt{phi\_double}} & 0.891088948393962 & 0.087 & 11 & Yes$^1$ \\ 
  &   & \textcolor{myblue}{\texttt{phi\_vpa}} & 0.891088548280465 & 0.073 & 1 & No \\ 
  &   & \textcolor{myblue}{\texttt{phi\_exact}} & 0.891088548280465 & 0.044 & 1 & No \\ 
\hline 
\textcolor{myblue}{\texttt{quad}} & 10000 & \textcolor{myblue}{\texttt{phi\_double}} & 0.891114412663896 & 0.091 & 24 & No \\ 
  &   & \textcolor{myblue}{\texttt{phi\_vpa}} & 0.891088530616134 & 0.181 & 26 & No \\ 
  &   & \textcolor{myblue}{\texttt{phi\_exact}} & 0.891088530616134 & 0.097 & 26 & No \\ 
\hline 
\textcolor{myblue}{\texttt{quadl}} & 10000 & \textcolor{myblue}{\texttt{phi\_double}} & 0.891069033973833 & 0.056 & 8 & No \\ 
  &   & \textcolor{myblue}{\texttt{phi\_vpa}} & 0.891088548278784 & 0.092 & 8 & No \\ 
  &   & \textcolor{myblue}{\texttt{phi\_exact}} & 0.891088548278784 & 0.056 & 8 & No \\ 
\hline 
\textcolor{myblue}{\texttt{trapz}} & 10000 & \textcolor{myblue}{\texttt{phi\_double}} & 0.891089828004092 & 0.034 & 1 & No \\ 
  &   & \textcolor{myblue}{\texttt{phi\_vpa}} & 0.891088548280466 & 2.385 & 1 & No \\ 
  &   & \textcolor{myblue}{\texttt{phi\_exact}} & 0.891088548280466 & 0.033 & 1 & No \\ 
\hline 
\textcolor{myblue}{\texttt{tanh-sinh}} & 10000 & \textcolor{myblue}{\texttt{phi\_double}} & 0.890985796123688 & 0.032 & 1 & No \\ 
  &   & \textcolor{myblue}{\texttt{phi\_vpa}} & 0.891088548280696 & 0.082 & 1 & No \\ 
  &   & \textcolor{myblue}{\texttt{phi\_exact}} & 0.891088548280696 & 0.035 & 1 & No \\ 
\hline 
\hline 
\textcolor{myblue}{\texttt{integral}} & 100000 & \textcolor{myblue}{\texttt{phi\_double}} & 0.909820722614225 & 0.094 & 12 & Yes$^1$ \\ 
  &   & \textcolor{myblue}{\texttt{phi\_vpa}} & 0.891088548280465 & 0.071 & 1 & No \\ 
  &   & \textcolor{myblue}{\texttt{phi\_exact}} & 0.891088548280465 & 0.043 & 1 & No \\ 
\hline 
\textcolor{myblue}{\texttt{quad}} & 100000 & \textcolor{myblue}{\texttt{phi\_double}} & 0.924049764163183 & 11.322 & 5006 & Yes$^2$ \\ 
  &   & \textcolor{myblue}{\texttt{phi\_vpa}} & 0.891088530616134 & 0.179 & 26 & No \\ 
  &   & \textcolor{myblue}{\texttt{phi\_exact}} & 0.891088530616134 & 0.099 & 26 & No \\ 
\hline 
\textcolor{myblue}{\texttt{quadl}} & 100000 & \textcolor{myblue}{\texttt{phi\_double}} & 0.930107304312503 & 4.773 & 2018 & Yes$^2$ \\ 
  &   & \textcolor{myblue}{\texttt{phi\_vpa}} & 0.891088548278784 & 0.090 & 8 & No \\ 
  &   & \textcolor{myblue}{\texttt{phi\_exact}} & 0.891088548278784 & 0.054 & 8 & No \\ 
\hline 
\textcolor{myblue}{\texttt{trapz}} & 100000 & \textcolor{myblue}{\texttt{phi\_double}} & 0.907718939549998 & 0.036 & 1 & No \\ 
  &   & \textcolor{myblue}{\texttt{phi\_vpa}} & 0.891088548280465 & 2.376 & 1 & No \\ 
  &   & \textcolor{myblue}{\texttt{phi\_exact}} & 0.891088548280466 & 0.036 & 1 & No \\ 
\hline 
\textcolor{myblue}{\texttt{tanh-sinh}} & 100000 & \textcolor{myblue}{\texttt{phi\_double}} & 0.890340539064744 & 0.032 & 1 & No \\ 
  &   & \textcolor{myblue}{\texttt{phi\_vpa}} & 0.891088548280696 & 0.050 & 1 & No \\ 
  &   & \textcolor{myblue}{\texttt{phi\_exact}} & 0.891088548280696 & 0.033 & 1 & No \\ 
\hline 
\hline 
\textcolor{myblue}{\texttt{integral}} & 250000 & \textcolor{myblue}{\texttt{phi\_double}} & 0.000000000000000 & 0.046 & 1 & No \\ 
  &   & \textcolor{myblue}{\texttt{phi\_vpa}} & 0.891088548280465 & 0.101 & 1 & No \\ 
  &   & \textcolor{myblue}{\texttt{phi\_exact}} & 0.891088548280465 & 0.042 & 1 & No \\ 
\hline 
\textcolor{myblue}{\texttt{quad}} & 250000 & \textcolor{myblue}{\texttt{phi\_double}} & 0.000000000000000 & 0.043 & 4 & No \\ 
  &   & \textcolor{myblue}{\texttt{phi\_vpa}} & 0.891088530616134 & 0.177 & 26 & No \\ 
  &   & \textcolor{myblue}{\texttt{phi\_exact}} & 0.891088530616134 & 0.099 & 26 & No \\ 
\hline 
\textcolor{myblue}{\texttt{quadl}} & 250000 & \textcolor{myblue}{\texttt{phi\_double}} & 0.000000000000000 & 0.043 & 2 & No \\ 
  &   & \textcolor{myblue}{\texttt{phi\_vpa}} & 0.891088548278784 & 0.094 & 8 & No \\ 
  &   & \textcolor{myblue}{\texttt{phi\_exact}} & 0.891088548278784 & 0.054 & 8 & No \\ 
\hline 
\textcolor{myblue}{\texttt{trapz}} & 250000 & \textcolor{myblue}{\texttt{phi\_double}} & 0.000000000000000 & 0.035 & 1 & No \\ 
  &   & \textcolor{myblue}{\texttt{phi\_vpa}} & 0.891088548280466 & 2.390 & 1 & No \\ 
  &   & \textcolor{myblue}{\texttt{phi\_exact}} & 0.891088548280466 & 0.037 & 1 & No \\ 
\hline 
\textcolor{myblue}{\texttt{tanh-sinh}} & 250000 & \textcolor{myblue}{\texttt{phi\_double}} & 0.000000000000000 & 0.032 & 1 & No \\ 
  &   & \textcolor{myblue}{\texttt{phi\_vpa}} & 0.891088548280696 & 0.059 & 1 & No \\ 
  &   & \textcolor{myblue}{\texttt{phi\_exact}} & 0.891088548280696 & 0.032 & 1 & No \\ 
\hline 
\hline 
\end{tabular} %
}
\end{center}
\end{table}

\begin{table}[ht!]
\begin{center}
\caption{Comparison of different quadratures used to calculate the value of Fourier coefficient $A_2$}\label{t:A2}
{\small%
\begin{tabular}{lclcccc}
\hline 
Integral Rule & $\delta$ & Integrand Rule & Value & Time [s] & Fevals & Warning \\ 
\hline 
\hline 
\textcolor{myblue}{\texttt{integral}} & 1000 & \textcolor{myblue}{\texttt{phi\_double}} & -0.132240670973319 & 0.117 & 12 & Yes$^1$ \\ 
  &   & \textcolor{myblue}{\texttt{phi\_vpa}} & -0.132240669593892 & 0.075 & 1 & No \\ 
  &   & \textcolor{myblue}{\texttt{phi\_exact}} & -0.132240669593892 & 0.043 & 1 & No \\ 
\hline 
\textcolor{myblue}{\texttt{quad}} & 1000 & \textcolor{myblue}{\texttt{phi\_double}} & -0.132240607546634 & 0.106 & 30 & No \\ 
  &   & \textcolor{myblue}{\texttt{phi\_vpa}} & -0.132240616544284 & 0.203 & 30 & No \\ 
  &   & \textcolor{myblue}{\texttt{phi\_exact}} & -0.132240616544284 & 0.108 & 30 & No \\ 
\hline 
\textcolor{myblue}{\texttt{quadl}} & 1000 & \textcolor{myblue}{\texttt{phi\_double}} & -0.132240665269117 & 0.058 & 8 & No \\ 
  &   & \textcolor{myblue}{\texttt{phi\_vpa}} & -0.132240669544252 & 0.094 & 8 & No \\ 
  &   & \textcolor{myblue}{\texttt{phi\_exact}} & -0.132240669544252 & 0.060 & 8 & No \\ 
\hline 
\textcolor{myblue}{\texttt{trapz}} & 1000 & \textcolor{myblue}{\texttt{phi\_double}} & -0.132240669210701 & 0.034 & 1 & No \\ 
  &   & \textcolor{myblue}{\texttt{phi\_vpa}} & -0.132240669593892 & 2.396 & 1 & No \\ 
  &   & \textcolor{myblue}{\texttt{phi\_exact}} & -0.132240669593892 & 0.035 & 1 & No \\ 
\hline 
\textcolor{myblue}{\texttt{tanh-sinh}} & 1000 & \textcolor{myblue}{\texttt{phi\_double}} & -0.132240685879949 & 0.033 & 1 & No \\ 
  &   & \textcolor{myblue}{\texttt{phi\_vpa}} & -0.132240671735414 & 0.080 & 1 & No \\ 
  &   & \textcolor{myblue}{\texttt{phi\_exact}} & -0.132240671735414 & 0.033 & 1 & No \\ 
\hline 
\hline 
\textcolor{myblue}{\texttt{integral}} & 10000 & \textcolor{myblue}{\texttt{phi\_double}} & -0.132241272405812 & 0.091 & 11 & Yes$^1$ \\ 
  &   & \textcolor{myblue}{\texttt{phi\_vpa}} & -0.132240669593892 & 0.074 & 1 & No \\ 
  &   & \textcolor{myblue}{\texttt{phi\_exact}} & -0.132240669593892 & 0.046 & 1 & No \\ 
\hline 
\textcolor{myblue}{\texttt{quad}} & 10000 & \textcolor{myblue}{\texttt{phi\_double}} & -0.132263860323737 & 0.129 & 38 & No \\ 
  &   & \textcolor{myblue}{\texttt{phi\_vpa}} & -0.132240616544284 & 0.199 & 30 & No \\ 
  &   & \textcolor{myblue}{\texttt{phi\_exact}} & -0.132240616544284 & 0.105 & 30 & No \\ 
\hline 
\textcolor{myblue}{\texttt{quadl}} & 10000 & \textcolor{myblue}{\texttt{phi\_double}} & -0.132218822360428 & 0.053 & 8 & No \\ 
  &   & \textcolor{myblue}{\texttt{phi\_vpa}} & -0.132240669544252 & 0.092 & 8 & No \\ 
  &   & \textcolor{myblue}{\texttt{phi\_exact}} & -0.132240669544252 & 0.054 & 8 & No \\ 
\hline 
\textcolor{myblue}{\texttt{trapz}} & 10000 & \textcolor{myblue}{\texttt{phi\_double}} & -0.132242712125416 & 0.039 & 1 & No \\ 
  &   & \textcolor{myblue}{\texttt{phi\_vpa}} & -0.132240669593892 & 2.380 & 1 & No \\ 
  &   & \textcolor{myblue}{\texttt{phi\_exact}} & -0.132240669593892 & 0.037 & 1 & No \\ 
\hline 
\textcolor{myblue}{\texttt{tanh-sinh}} & 10000 & \textcolor{myblue}{\texttt{phi\_double}} & -0.132181727630471 & 0.033 & 1 & No \\ 
  &   & \textcolor{myblue}{\texttt{phi\_vpa}} & -0.132240671735414 & 0.050 & 1 & No \\ 
  &   & \textcolor{myblue}{\texttt{phi\_exact}} & -0.132240671735414 & 0.032 & 1 & No \\ 
\hline 
\hline 
\textcolor{myblue}{\texttt{integral}} & 100000 & \textcolor{myblue}{\texttt{phi\_double}} & -0.164105328139715 & 0.124 & 12 & Yes$^1$ \\ 
  &   & \textcolor{myblue}{\texttt{phi\_vpa}} & -0.132240669593892 & 0.070 & 1 & No \\ 
  &   & \textcolor{myblue}{\texttt{phi\_exact}} & -0.132240669593892 & 0.043 & 1 & No \\ 
\hline 
\textcolor{myblue}{\texttt{quad}} & 100000 & \textcolor{myblue}{\texttt{phi\_double}} & -0.217677705429780 & 11.343 & 5006 & Yes$^2$ \\ 
  &   & \textcolor{myblue}{\texttt{phi\_vpa}} & -0.132240616544284 & 0.196 & 30 & No \\ 
  &   & \textcolor{myblue}{\texttt{phi\_exact}} & -0.132240616544284 & 0.108 & 30 & No \\ 
\hline 
\textcolor{myblue}{\texttt{quadl}} & 100000 & \textcolor{myblue}{\texttt{phi\_double}} & -0.175317505741447 & 4.736 & 2012 & Yes$^2$ \\ 
  &   & \textcolor{myblue}{\texttt{phi\_vpa}} & -0.132240669544252 & 0.094 & 8 & No \\ 
  &   & \textcolor{myblue}{\texttt{phi\_exact}} & -0.132240669544252 & 0.056 & 8 & No \\ 
\hline 
\textcolor{myblue}{\texttt{trapz}} & 100000 & \textcolor{myblue}{\texttt{phi\_double}} & -0.162243580664938 & 0.037 & 1 & No \\ 
  &   & \textcolor{myblue}{\texttt{phi\_vpa}} & -0.132240669593892 & 2.380 & 1 & No \\ 
  &   & \textcolor{myblue}{\texttt{phi\_exact}} & -0.132240669593892 & 0.033 & 1 & No \\ 
\hline 
\textcolor{myblue}{\texttt{tanh-sinh}} & 100000 & \textcolor{myblue}{\texttt{phi\_double}} & -0.114446017197740 & 0.032 & 1 & No \\ 
  &   & \textcolor{myblue}{\texttt{phi\_vpa}} & -0.132240671735414 & 0.053 & 1 & No \\ 
  &   & \textcolor{myblue}{\texttt{phi\_exact}} & -0.132240671735414 & 0.060 & 1 & No \\ 
\hline 
\hline 
\textcolor{myblue}{\texttt{integral}} & 250000 & \textcolor{myblue}{\texttt{phi\_double}} & 0.000000000000000 & 0.043 & 1 & No \\ 
  &   & \textcolor{myblue}{\texttt{phi\_vpa}} & -0.132240669593892 & 0.071 & 1 & No \\ 
  &   & \textcolor{myblue}{\texttt{phi\_exact}} & -0.132240669593892 & 0.045 & 1 & No \\ 
\hline 
\textcolor{myblue}{\texttt{quad}} & 250000 & \textcolor{myblue}{\texttt{phi\_double}} & 0.000000000000000 & 0.043 & 4 & No \\ 
  &   & \textcolor{myblue}{\texttt{phi\_vpa}} & -0.132240616544284 & 0.204 & 30 & No \\ 
  &   & \textcolor{myblue}{\texttt{phi\_exact}} & -0.132240616544284 & 0.105 & 30 & No \\ 
\hline 
\textcolor{myblue}{\texttt{quadl}} & 250000 & \textcolor{myblue}{\texttt{phi\_double}} & 0.000000000000000 & 0.039 & 2 & No \\ 
  &   & \textcolor{myblue}{\texttt{phi\_vpa}} & -0.132240669544252 & 0.094 & 8 & No \\ 
  &   & \textcolor{myblue}{\texttt{phi\_exact}} & -0.132240669544252 & 0.054 & 8 & No \\ 
\hline 
\textcolor{myblue}{\texttt{trapz}} & 250000 & \textcolor{myblue}{\texttt{phi\_double}} & 0.000000000000000 & 0.035 & 1 & No \\ 
  &   & \textcolor{myblue}{\texttt{phi\_vpa}} & -0.132240669593892 & 2.393 & 1 & No \\ 
  &   & \textcolor{myblue}{\texttt{phi\_exact}} & -0.132240669593892 & 0.033 & 1 & No \\ 
\hline 
\textcolor{myblue}{\texttt{tanh-sinh}} & 250000 & \textcolor{myblue}{\texttt{phi\_double}} & 0.000000000000000 & 0.033 & 1 & No \\ 
  &   & \textcolor{myblue}{\texttt{phi\_vpa}} & -0.132240671735414 & 0.054 & 1 & No \\ 
  &   & \textcolor{myblue}{\texttt{phi\_exact}} & -0.132240671735414 & 0.033 & 1 & No \\ 
\hline 
\hline 
\end{tabular} %
}
\end{center}
\end{table}

 \clearpage
\section{Problematic use-cases}\label{sec:models}

In this section we introduce three problematic use-cases, one linear and one non-linear ordinary differential equation (ODE) and one linear partial differential equation (PDE). 

{\let\addcontentsline\nocontentsline%
\subsection{Linear ODE}\label{ss:4:linode}
}

We consider an initial problem for linear ODE of the second order
\begin{equation}\label{e:linode}
x''(t) + k^2 x(t) = 0, \quad x(0)=0,\,x'(0)=k,
\end{equation}
where $k\in\R$ is a parameter. Initial value problem \eqref{e:linode} has the unique exact solution $x(t) = \sin(k\cdot t)$. 

\begin{figure}[!ht]
\begin{center}
\includegraphics[width=0.9\textwidth]{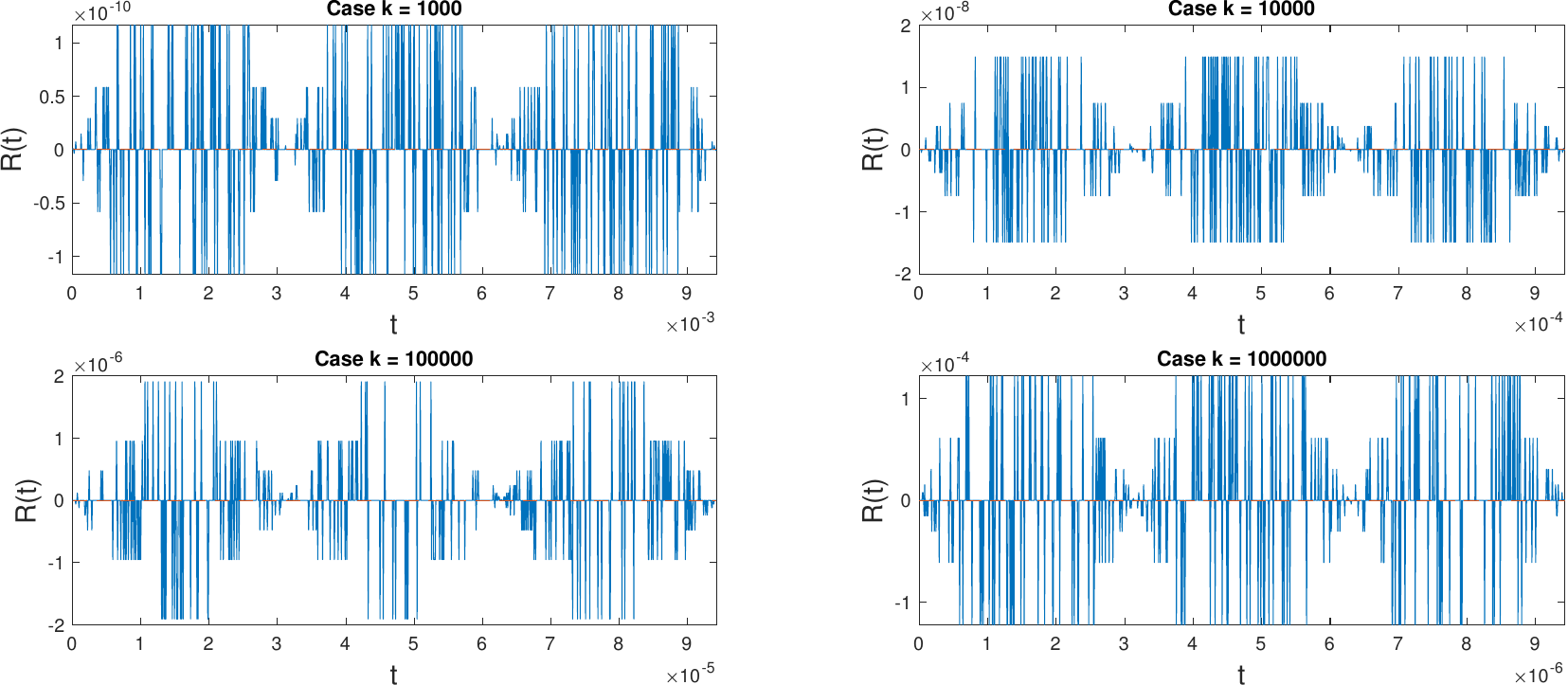}
\end{center}
\caption{Residuum $R(t)$ calculated using AD for different values of $k$.}\label{f:residuals}
\end{figure}

Let $R(t) = x''(t) + k^2 x(t)$ be the residuum for problem \eqref{e:linode} and let us calculate it using the automatic differentiation over the first three periods of the solution. In Figure \ref{f:residuals} we can see that the increasing value of $k$ causes also the higher value of the residua that should be theoretically zero. This problem belongs to the category of highly oscillatory functions mentioned in Section \ref{ss:list}.

{\let\addcontentsline\nocontentsline%
\subsection{Non-linear ODE}\label{ss:4:nonlinode}
}

We consider an initial problem for the \alert{van der Pol oscillator} modelled as the non-linear ODE of the second order
\begin{equation}\label{e:nonlinode}
x''(t) - \mu(1-x^2(t))x'(t) + x(t) = 0, \quad x(0)=2,\,x'(0)=0,
\end{equation}
where $\mu\in\R$ is a dumping parameter. Although the theorem of Picard–Lindel\"of guarantees a unique analytic solution, it does not have an exact solution in closed form, it can only be written as a Taylor (Maclaurin) series. We can rewrite \eqref{e:nonlinode} as the system
\begin{align}\label{e:nonlinodesys}
\begin{split}
    x'(t) &= y(t),  \\
    y'(t) &= \mu(1-x^2(t))y(t) - x(t), \\
    x(0) &= 2, \\
    y(0) &= 0.
\end{split}
\end{align}
To solve the system \eqref{e:nonlinodesys} numerically, we may use some suitable solver for stiff differential equations such as \texttt{ode15s} in MATLAB. In Section \ref{ss:5:nonlinode} we take \texttt{ode15s} solution as a reference when comparing it to a solution using PINNs.

{\let\addcontentsline\nocontentsline%
\subsection{Heat equation problem}\label{sec:heat}
}

We consider the heat equation in one spatial variable with Dirichlet boundary conditions:
\begin{align}\label{e:heat}
\begin{split}
\dfrac{\partial u}{\partial t}(t,x) - \alpha \dfrac{\partial^2 u}{\partial x^2}(t,x) &=0, \quad t\in(0,1], x\in(-1,1), \\
u(t,-1) = u(t,1) &= 0, \quad t\in[0,1], \\
u(0,x) &= \alert{\phi(x)}, \quad x\in[-1,1]. 
\end{split}
\end{align}
A solution by separation of variables with 
$f_n(t) = \e^{-\alpha t ((2n-1)\pi/2)^2}$, base functions $g_n(x)$ given by \eqref{e:g_n} and Fourier coefficients $A_n$ given by \eqref{e:A_n} takes the form
\begin{align}
u(t,x) &= \sum_{n=1}^{+\infty} A_n f_n(t) g_n(x). \label{e:u}
\intertext{We denote}
u_n(t,x) &= A_n f_n(t) g_n(x) \label{e:u_n}
\intertext{and the truncated solution after the $N$-th term by}
u^N(t,x) &= \sum_{n=1}^{N} u_n(t,x). \label{e:u_N}
\end{align}

\begin{figure}[ht!]
\includegraphics[width=\textwidth]{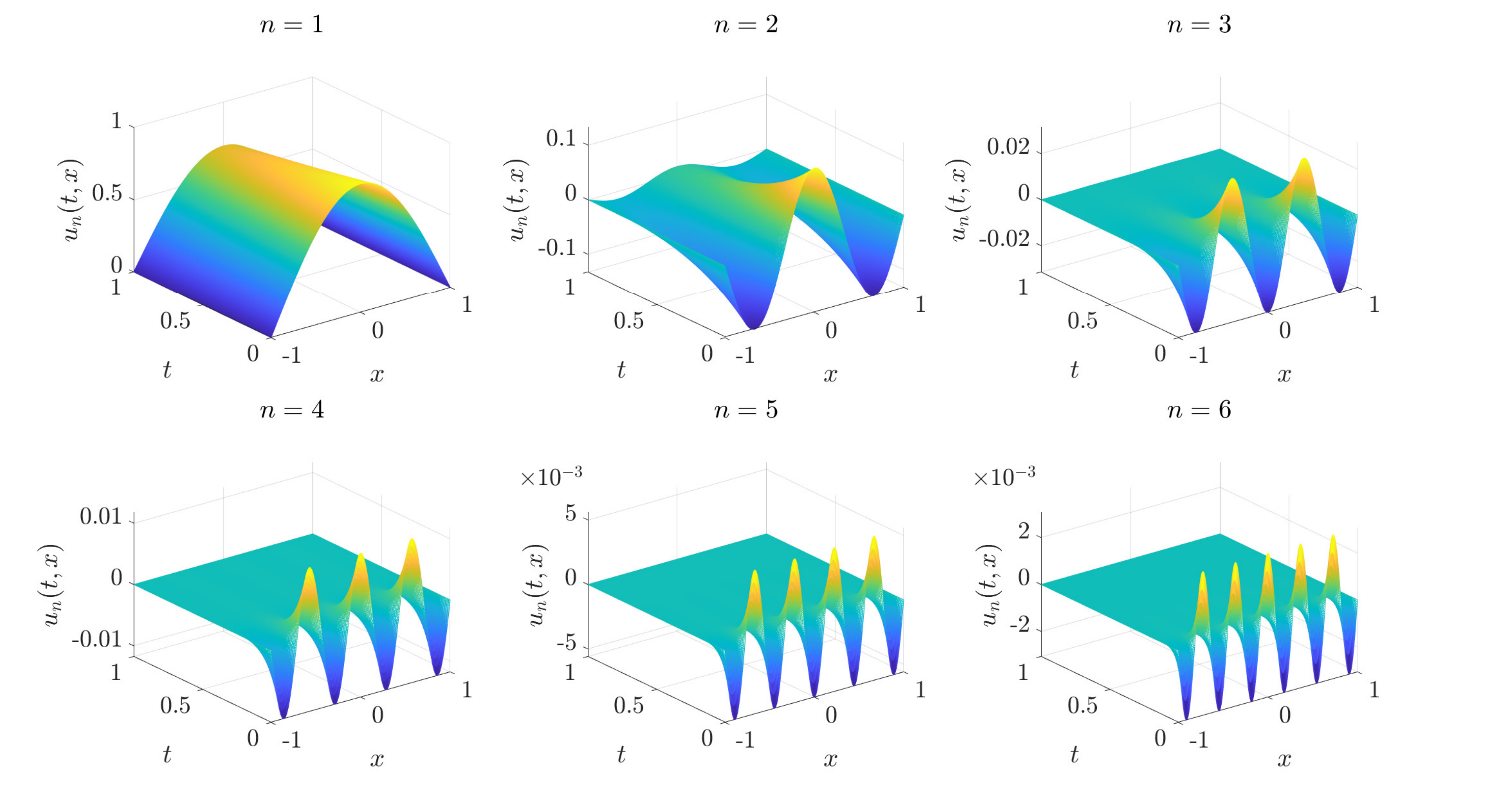}\\
\includegraphics[width=\textwidth]{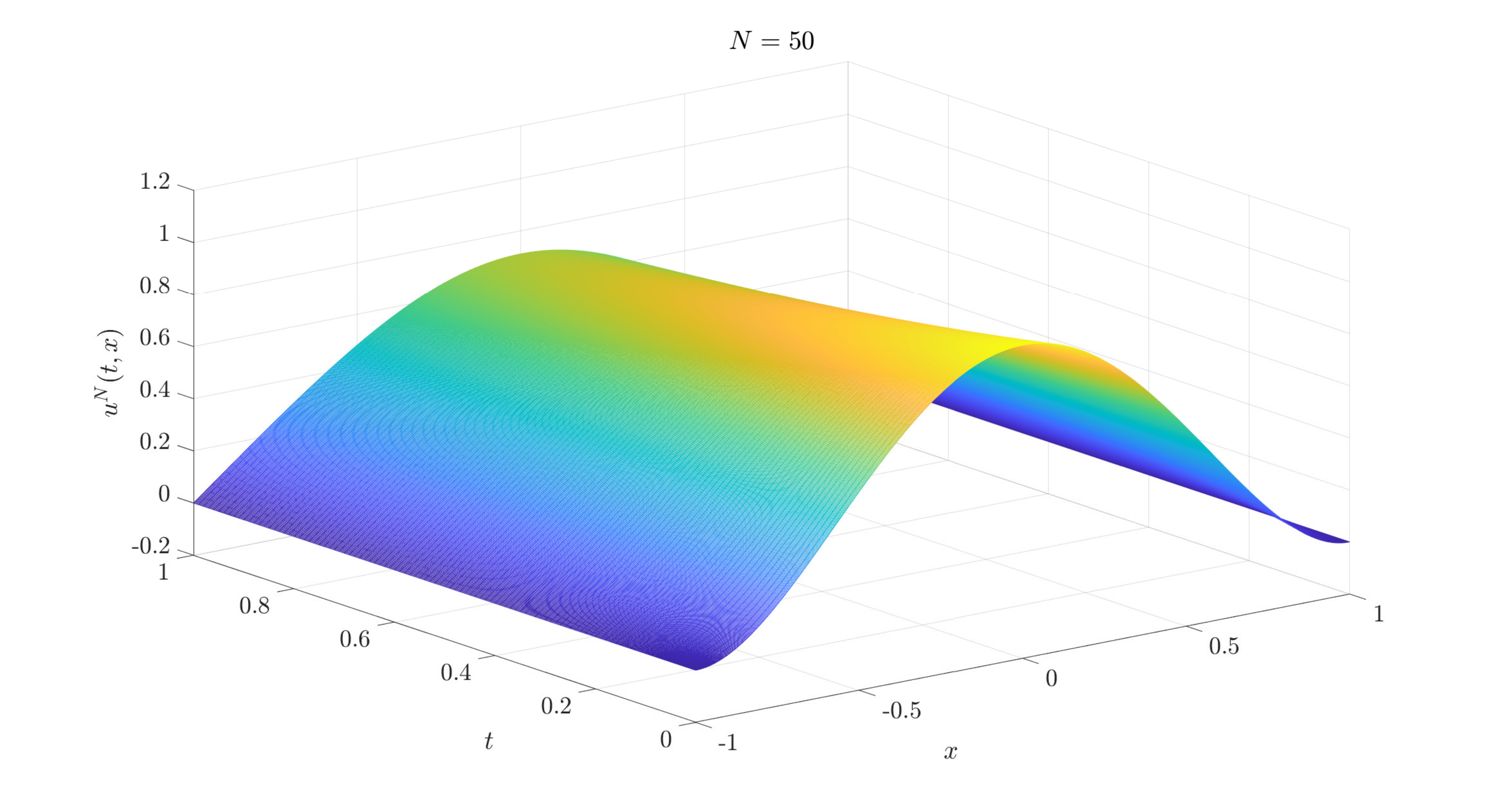}
\caption{Graphs of the first six $u_n(t,x)$, $n=1,\dots,6$, (top) and $u^N(t,x)$ for $N=50$ (bottom) in the case when the initial condition $\alert{\phi(x)}$ was taken as the \texttt{phi\_vpa()}.}\label{f:u_exact_N}
\end{figure}

In Figure~\ref{f:u_exact_N} we can see graphs of the first six $u_n(t,x)$, $n=1,\dots,6$, (top) given by \eqref{e:u_n} and $u^N(t,x)$ for $N=50$ (bottom) given by \eqref{e:u_N} in the case when the initial condition $\alert{\phi(x)}$ was taken as the \texttt{phi\_exact()} or \texttt{phi\_vpa()} (the difference between these two cases is again at the order of machine precision and hence only one figure is provided).

\begin{figure}[ht!]
\includegraphics[width=\textwidth]{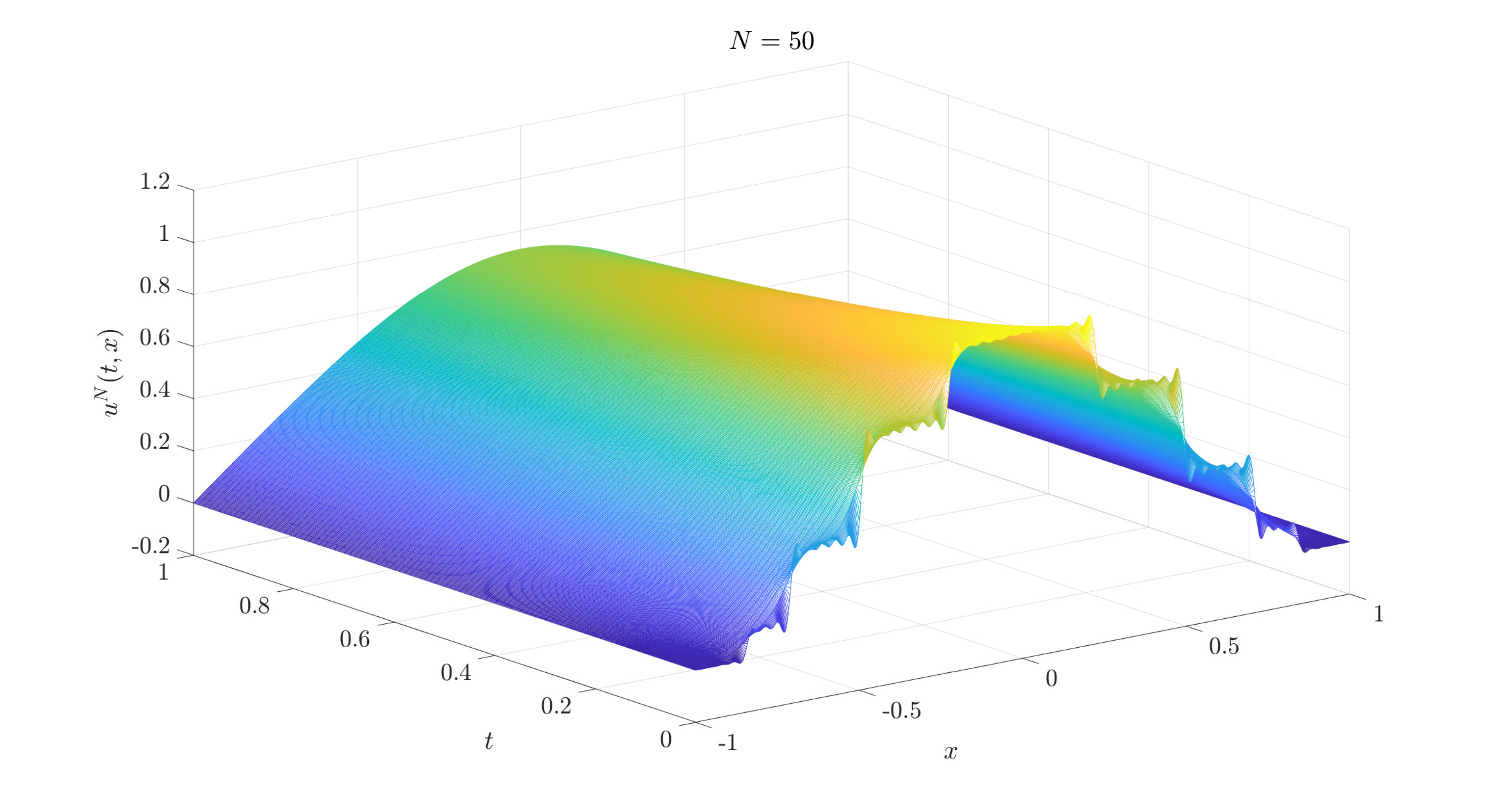}
\includegraphics[width=\textwidth]{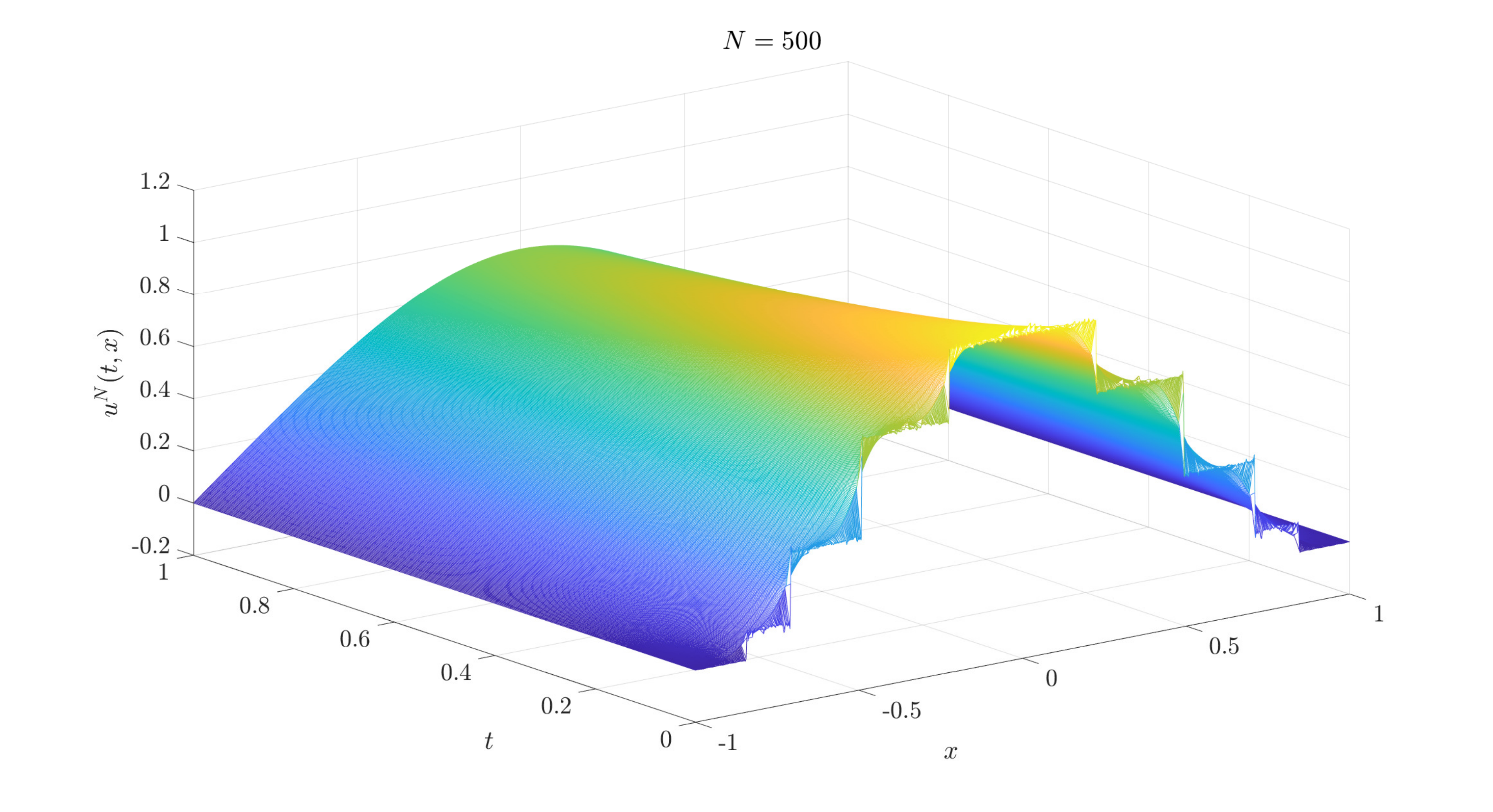}
\caption{Graphs of $u^N(t,x)$ for $N=50$ (top) and $N=500$ (bottom) in the case when the initial condition $\alert{\phi(x)}$ was taken as the \texttt{phi\_double()}.}\label{f:u_double_N}
\end{figure}

If $\alert{\phi(x)}$ is taken as \texttt{phi\_double()}, the graphs of $u_n(t,x)$ are visually indistinguishable from the exact version presented in Figure~\ref{f:u_exact_N}. However,  the graph of $u^N(t,x)$ now differs significantly as we can see in Figure~\ref{f:u_double_N}. Inaccurately evaluated initial condition $\alert{\phi(x)}$ and hence the inaccurately evaluated Fourier coefficients $A_n$ causes the serious error at the beginning of the time domain. Figure~\ref{f:u_double_N} also shows the influence of the truncation $N$, namely that increasing $N$ does not help to solve the problem caused by the inaccurate evaluation described above.

 \clearpage
\section{Physics-informed neural network models}\label{sec:pinns}

In this section we now solve all three problematic cases from Section \ref{sec:models} using PINNs. 
\subsection{Linear ODE}\label{ss:5:linode}

Let us now consider a PINNs solution of the problem \eqref{ss:4:linode}. Training data are generated on the go in $N_f=128$ random collocation points $t_i$ covering the first period, i.e. the interval $[0,2\pi/k]$. We use a multilayer perceptron NN architecture with 2 fully connected hidden layers with 50 neurons each. Input layer has one input channel corresponding to the input $t$, output layer has one output channel $x(t)$. As activation functions, we take hyperbolic tangents. Total number of learnable parameters is 2700, all of type \texttt{double}. During the NN training, the \alert{Adam optimizer} is used with the default setting (\texttt{beta1=0.9}, \texttt{beta2=0.999} and \texttt{epsilon=1e-8}). In the following experiments, number of epochs and the value of the learning rate (LR) differs.  

We consider the model loss function: 
$$\text{loss} = \text{MSE}_R + \lambda\cdot\text{MSE}_I,\text{ where } \text{MSE}_R = \frac{1}{N_f}\sum\limits_{i=1}^{N_f} \left\vert R(t_i) \right\vert ^2 \text{ and }\text{MSE}_I = x^2(0) + (x'(0)-k)^2.$$
Parameter $\lambda$ is used to put more stress on satisfying the initial conditions, in particular the value $\lambda=k^2$ is used. Calculating $\text{MSE}_R$ and $\text{MSE}_I$ requires the derivatives $x'$ and $x''$ of the output $x$ of the model using automatic differentiation.

\begin{figure}[!ht]
\includegraphics[width=0.49\textwidth]{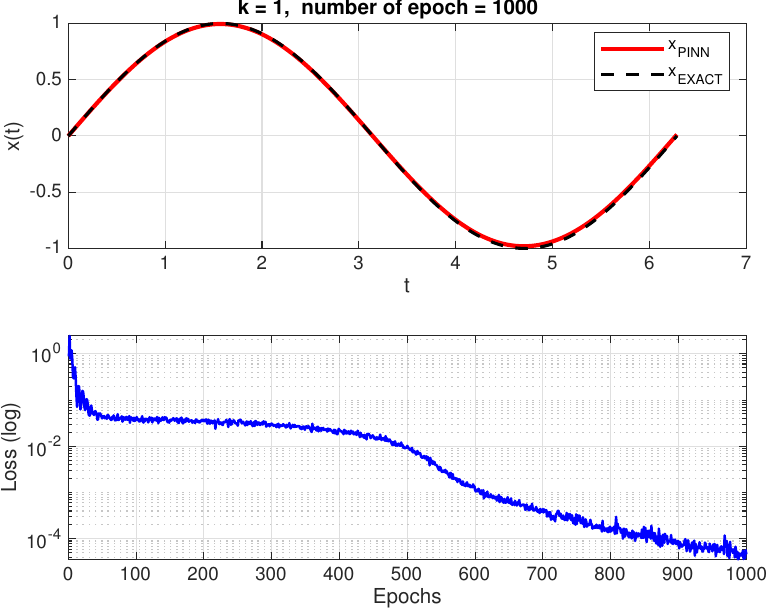}
\includegraphics[width=0.49\textwidth]{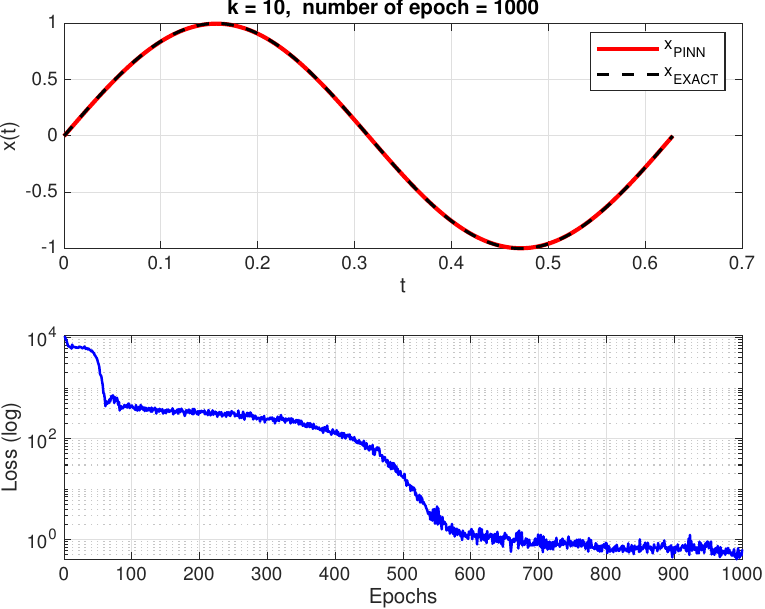}\\[4mm]
\includegraphics[width=0.49\textwidth]{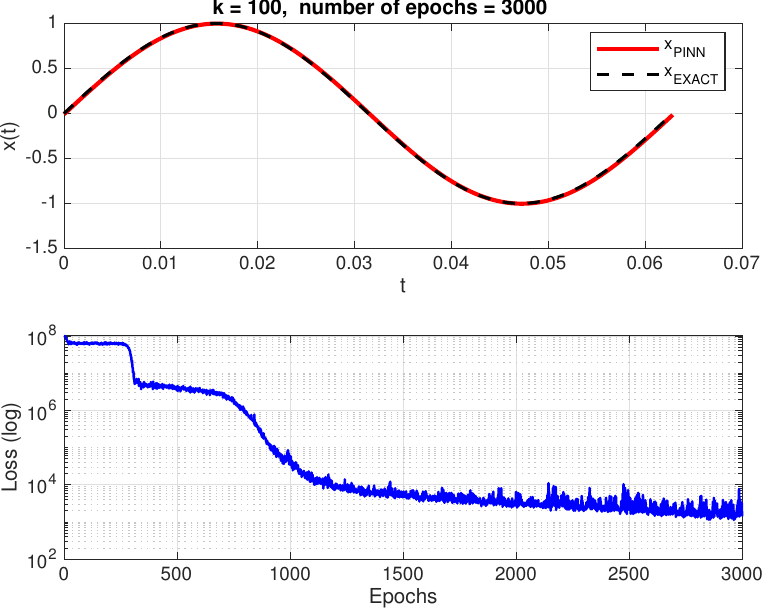}
\includegraphics[width=0.49\textwidth]{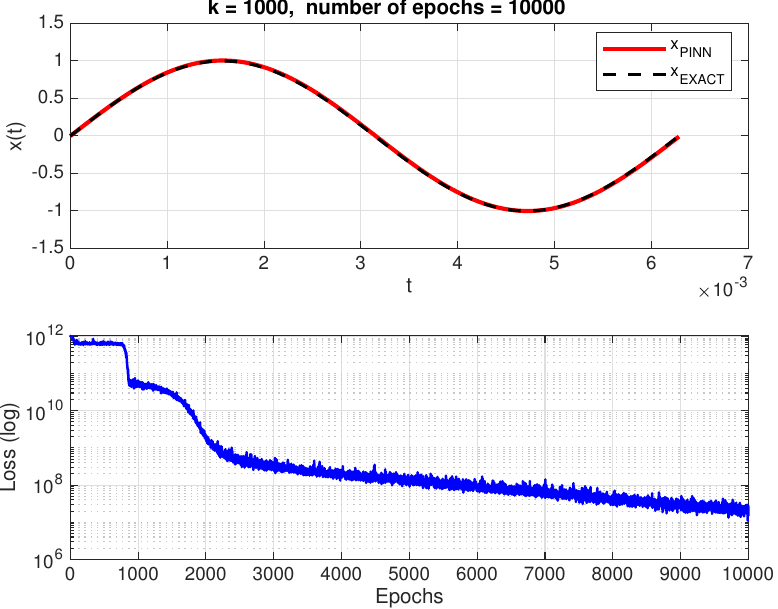}
\caption{Comparison of the PINNs solution to the exact solution of problem \eqref{e:linode}, four experiments considered, see Table \ref{t:linode}.}\label{f:linode}
\end{figure}

In Figure \ref{f:linode} we present a comparison of the PINNs solution to the exact solution for relatively small values of the parameter $k\in\{1,10,100,1000\}$. An increasing $k$ requires an increased number of epochs. The learning rate is changed adaptively (in a prescribed schedule) according to the Table \ref{t:linode} below. In all four experiments, the PINNs solution is satisfactory, revealing no numerical problems. 

\begin{table}[!ht]
\caption{Settings used in experiments depicted in Figure \ref{f:linode}}\label{t:linode}
\begin{center}
\begin{tabular}{lcccc}
\hline 
$k$  & 1 & 10 & 100 & 1000 \\ 
\hline 
\hline 
number of epochs & 1000 & 1000 & 3000 & 10000 \\
\texttt{LRinit}  & 1e-2 & 1e-2 & 1e-2 & 1.3e-2 \\
adaptive \texttt{LR} factor  & 0.99 & 0.99 & 0.98 & 0.97 \\
after each \# epoch  & 10 & 10 & 100 & 100 \\
\hline
\end{tabular} \end{center}
\end{table}

\begin{figure}[!ht]
\includegraphics[width=0.49\textwidth]{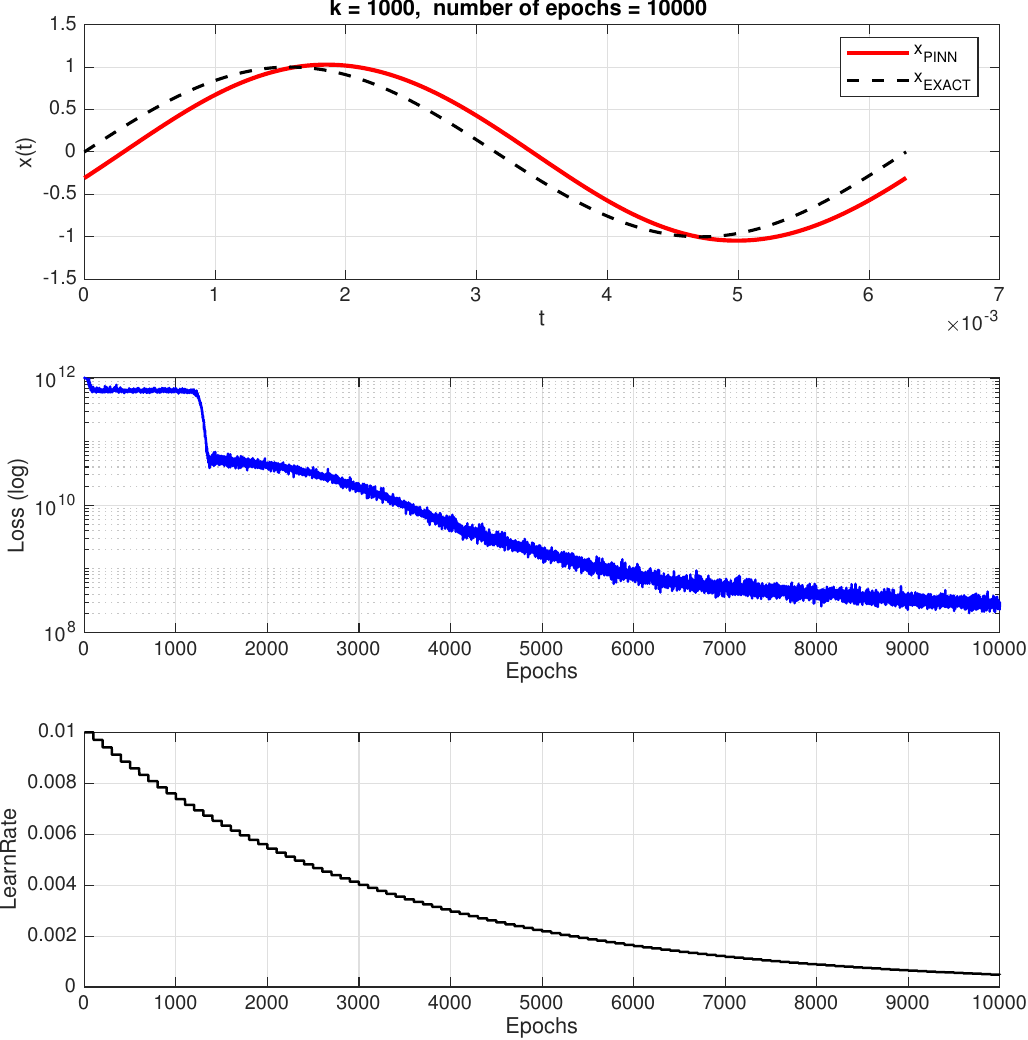}
\includegraphics[width=0.49\textwidth]{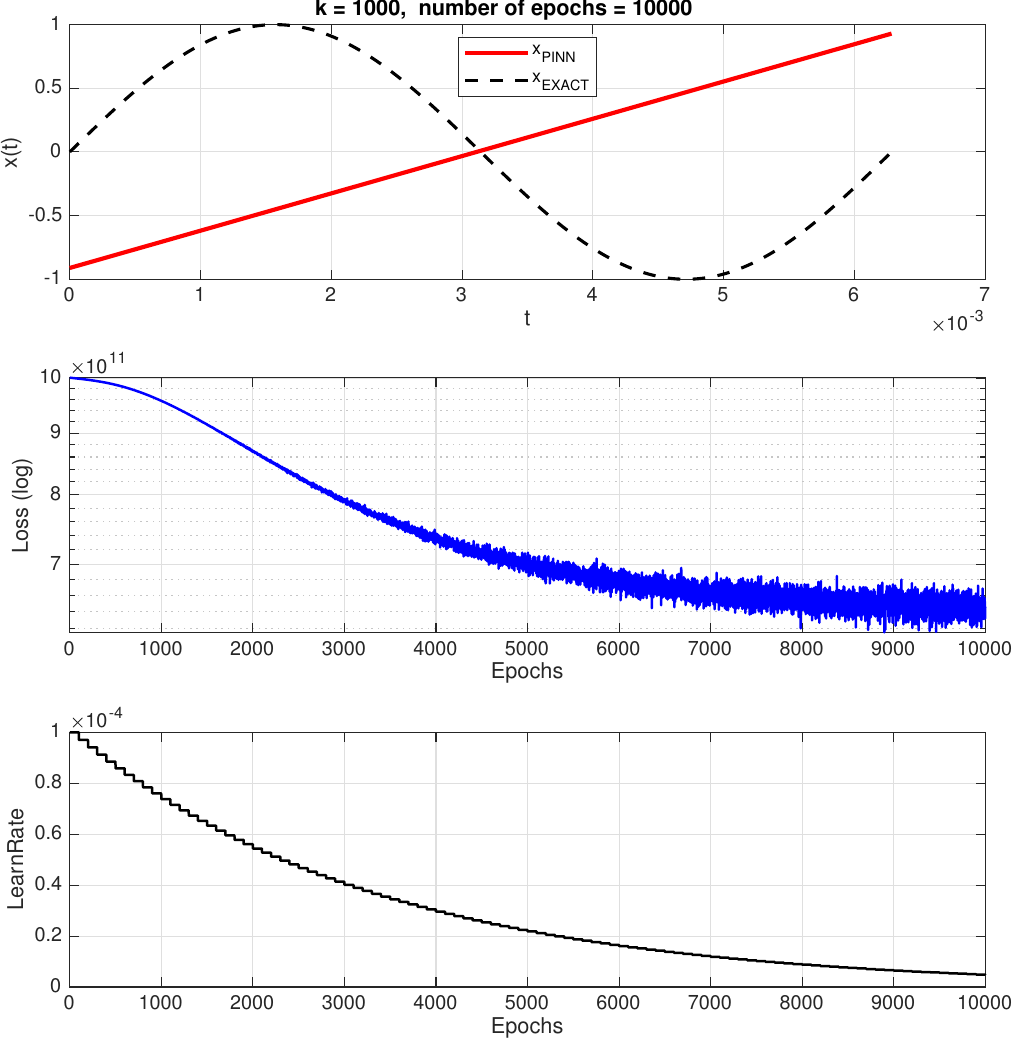}
\caption{Comparison of the PINNs solution to the exact solution of problem \eqref{e:linode}, with $k=1000$ and adaptive LR with \texttt{LRinit=1e-2} (left) and \texttt{LRinit=1e-4} (right).}\label{f:linode2}
\end{figure}

\begin{figure}[!ht]
\includegraphics[width=0.49\textwidth]{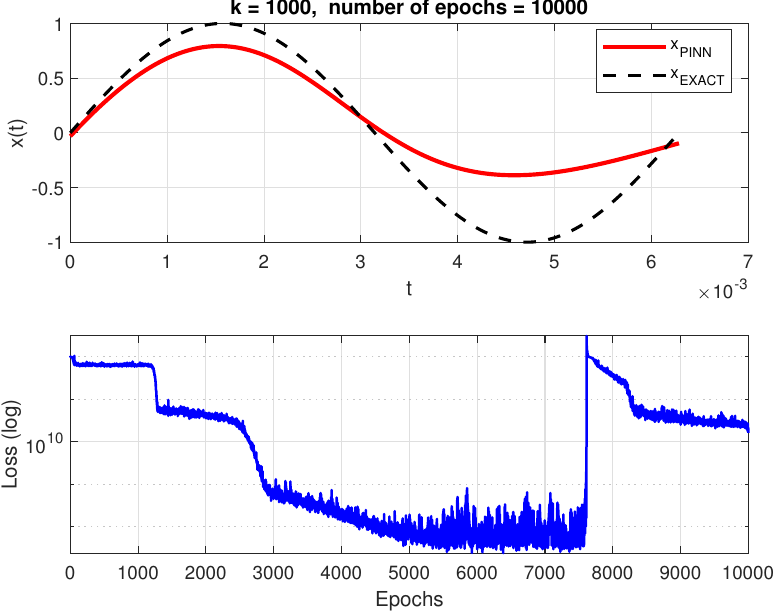}
\includegraphics[width=0.49\textwidth]{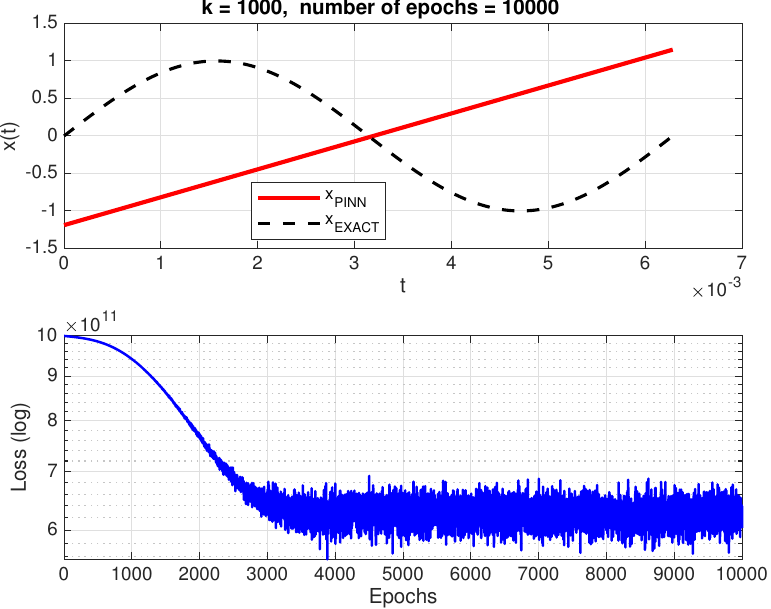}
\caption{Comparison of the PINNs solution to the exact solution of problem \eqref{e:linode}, with $k=1000$ and fixed \texttt{LR=1e-2} (left) and \texttt{LR=1e-4} (right).}\label{f:linode3}
\end{figure}

In Section \ref{ss:4:linode} we showed that error of automatic differentiation applied to \eqref{e:linode} increases with increasing value of $k$. In Figures \ref{f:linode2} and \ref{f:linode3} we consider only the value $k=1000$ that can lead to severe numerical problems of the PINNs solution when the setting is not carefully and suitably chosen. In particular, in Figure \ref{f:linode2}, only the initial learning rate is different compared to the settings $k=1000$ in Table \ref{t:linode} (the last column). Whereas on the left, the PINNs solution is biased with smaller error in initial conditions, the PINNs solution on the right is not satisfactory at all, in both cases after 10000 epochs. For convenience, we depict here also the learning rate decay as described in Table \ref{t:linode}. Similar misbehaviour can be observed in the case, where the learning rate is fixed in all epochs, see  Figure \ref{f:linode3}. Constant learning rate causes higher oscillations in the loss function (even with a big ``jump'' after ca 7500 epoch in the left case) and consequently worse results. 

\clearpage
\subsection{Nonlinear ODE}\label{ss:5:nonlinode}

Let us now consider a PINNs solution of the problem \eqref{ss:4:nonlinode}. Training data are generated on the go in $N_f=128$ random collocation points $t_i$ covering the interval $[0,7]$. We use a multilayer perceptron NN architecture as in the previous section, i.e. with 2 fully connected hidden layers with 50 neurons each. Input layer has one input channel corresponding to the input $t$, output layer has one output channel $x(t)$. As activation functions, we take hyperbolic tangents. Total number of learnable parameters is again 2700, all of type \texttt{double}. During the NN training, the \alert{Adam optimizer} is used with the default setting (\texttt{beta1=0.9}, \texttt{beta2=0.999} and \texttt{epsilon=1e-8}). The learning rate is changed adaptively (in a prescribed schedule) according to the Table \ref{t:nonlinode} below.
\begin{table}[!ht]
\caption{Settings used in experiments depicted in Figures \ref{f:nonlinode} and \ref{f:nonlinodebad}}\label{t:nonlinode}
\begin{center}
\begin{tabular}{lcccccc}
\hline 
$\mu$  & 1 & 2 & 2 & 2 & 2 & 2\\ 
\hline 
\hline 
figure & \ref{f:nonlinode}(a) &  \ref{f:nonlinode}(b) &
         \ref{f:nonlinodebad}(a) & \ref{f:nonlinodebad}(b) & 
         \ref{f:nonlinodebad}(c) & \ref{f:nonlinodebad}(d) \\  
number of epochs & 3000 & 25000 & 50000 & 70000  & 50000  & 50000 \\
\texttt{LRinit}  & 7e-3 & 8e-3 & 2.5e-2 & 6e-3 & 3e-3 & 5e-3 \\
adaptive \texttt{LR} factor  & 0.97 & 0.99 & 0.99 & 0.99 & 0.99 & 0.99 \\
after each \# epoch  & 100 & 100 & 100 & 100 & 200 & 200 \\
\hline
\end{tabular} \end{center}
\end{table}

We consider the model loss function: $\text{loss} = \text{MSE}_R + \lambda\cdot\text{MSE}_I$, where
\begin{align*}
\text{MSE}_R &= \frac{1}{N_f}\sum\limits_{i=1}^{N_f} \left\vert x''(t_i) - \mu(1-x^2(t_i))x'(t_i) + x(t_i) \right\vert ^2, \\
\text{MSE}_I &= (x(0)-2)^2 + (x'(0))^2.
\end{align*}
Parameter $\lambda$ can be used to put more stress on satisfying the initial conditions, however, in the following experiments when $\mu\in\{1,2\}$, the value of $\lambda$ does not bring any improvements, so that the value $\lambda=\mu$ is used for convenience. Calculating $\text{MSE}_R$ and $\text{MSE}_I$ requires the derivatives $x'$ and $x''$ of the output $x$ of the model using automatic differentiation.

\begin{figure}[!ht]
\includegraphics[width=0.49\textwidth]{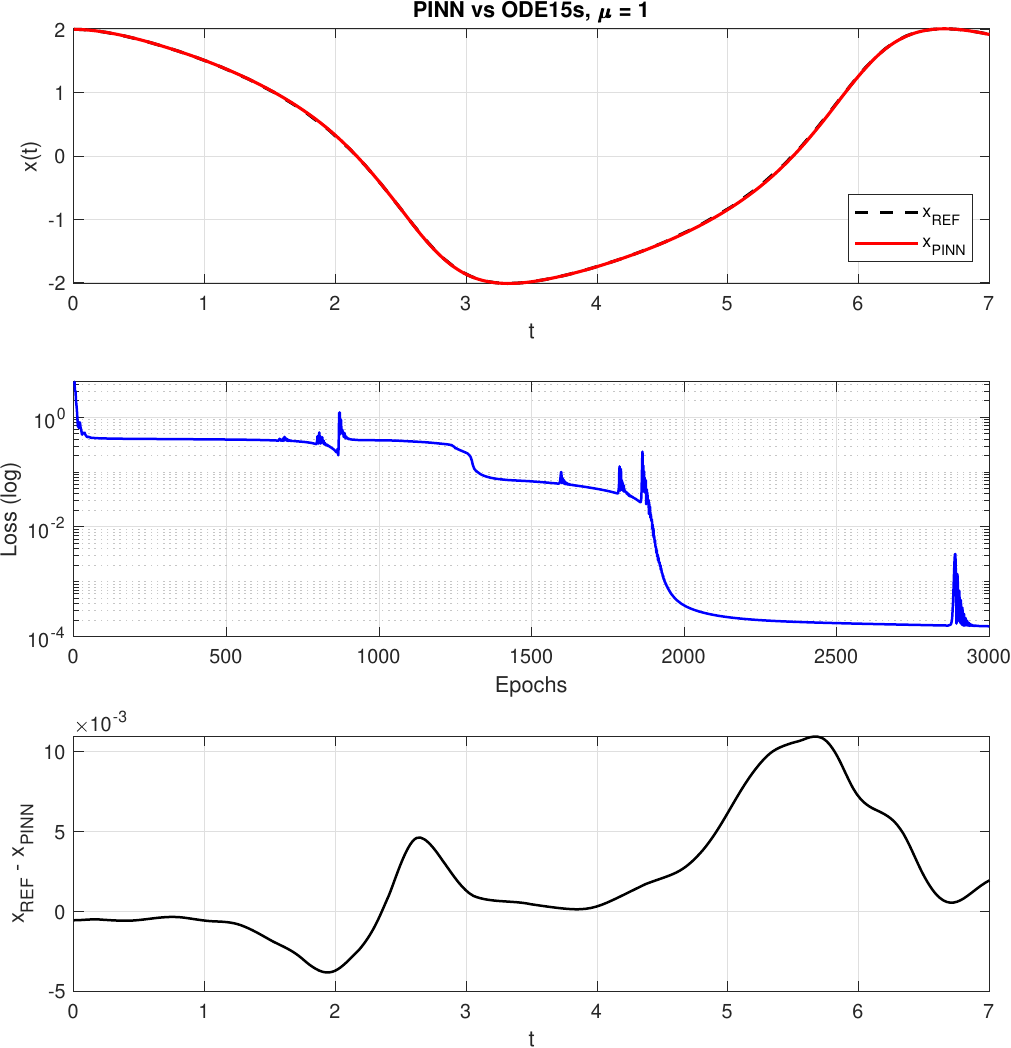}
\includegraphics[width=0.49\textwidth]{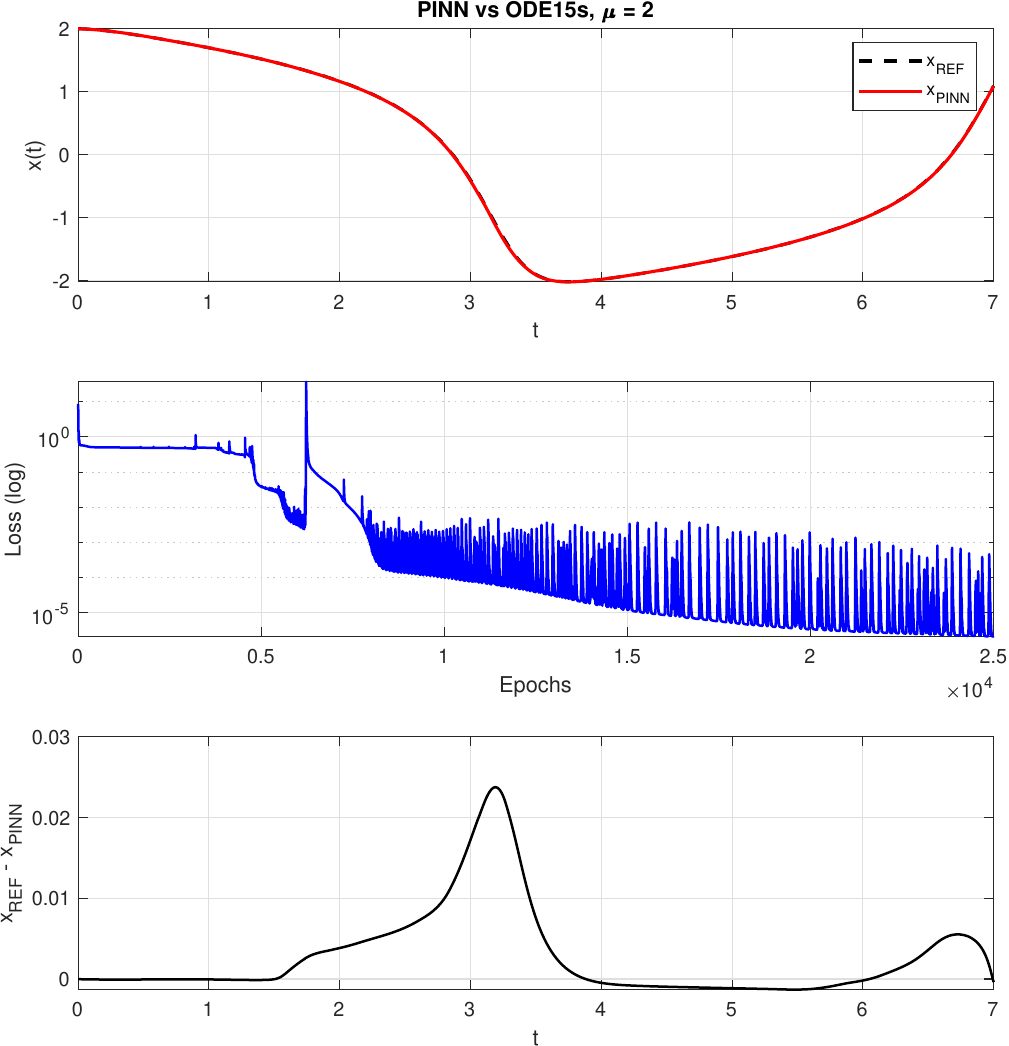}\\[-.1mm]
\mbox{ \hspace*{1mm} } \hfill (a) \hfill \mbox{} \hfill (b) \hfill \mbox{} 
\caption{Comparison of the PINNs solution (correct) to the reference numerical solution (\texttt{ode15s}) of problem \eqref{e:nonlinode} for (a)~$\mu=1$ and (b)~$\mu=2$.}\label{f:nonlinode}
\end{figure}

In Figure \ref{f:nonlinode} we present a comparison of the PINNs solution to the reference numerical solution obtained by \texttt{ode15s} for different values of the parameter $\mu\in\{1,2\}$. For $\mu=1$ the PINNs solution is satisfactory almost regardless of the LR settings. For $\mu=2$, the PINNs solution becomes highly sensitive not only to LR settings, but also to a particular realization of the random number generator used for the stochastic gradient descent algorithm. In  case (b), the PINNs solution is satisfactory, but it was obtained as one of the many realizations, although using the same LR settings.  

\begin{figure}[!ht]
\includegraphics[width=0.49\textwidth]{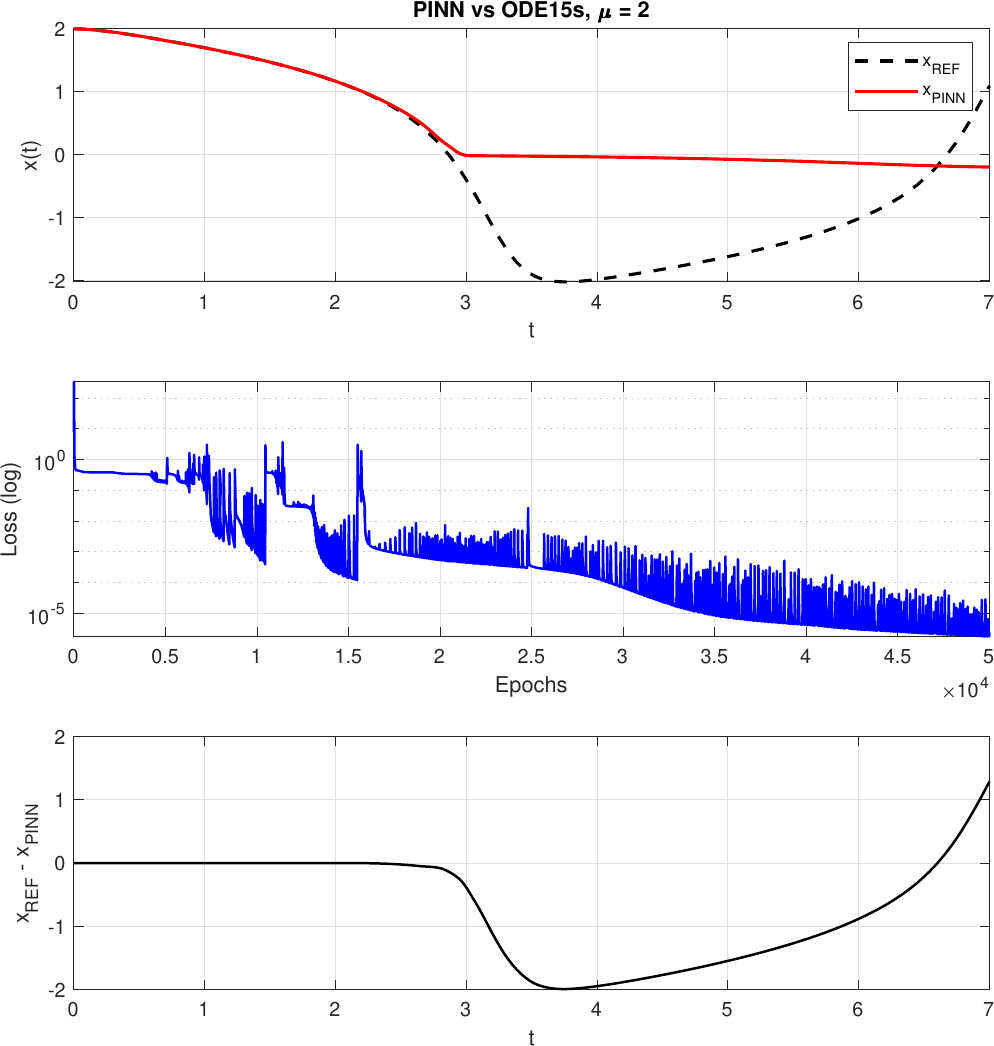}
\includegraphics[width=0.49\textwidth]{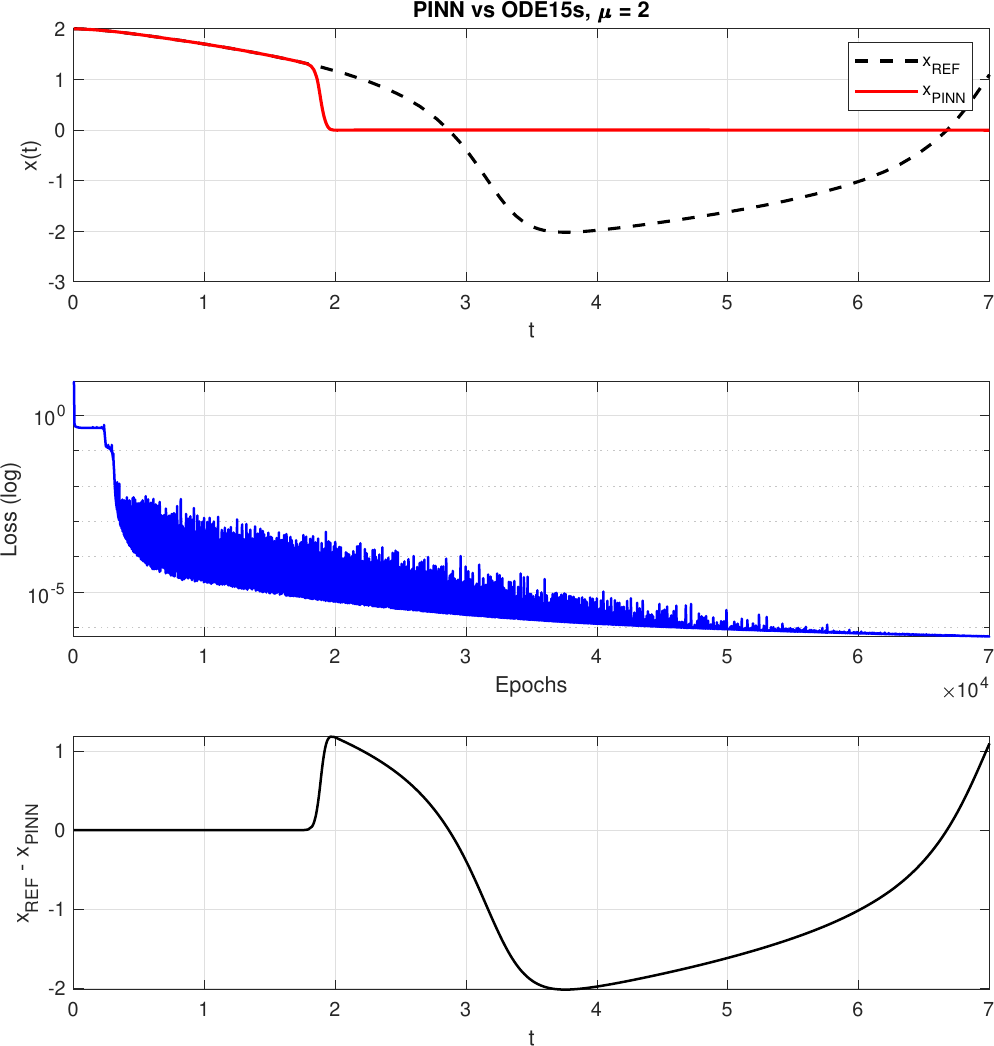}\\[-.1mm]
\mbox{ \hspace*{1mm} } \hfill (a) \hfill \mbox{} \hfill (b) \hfill \mbox{}  \\[5mm]
\includegraphics[width=0.49\textwidth]{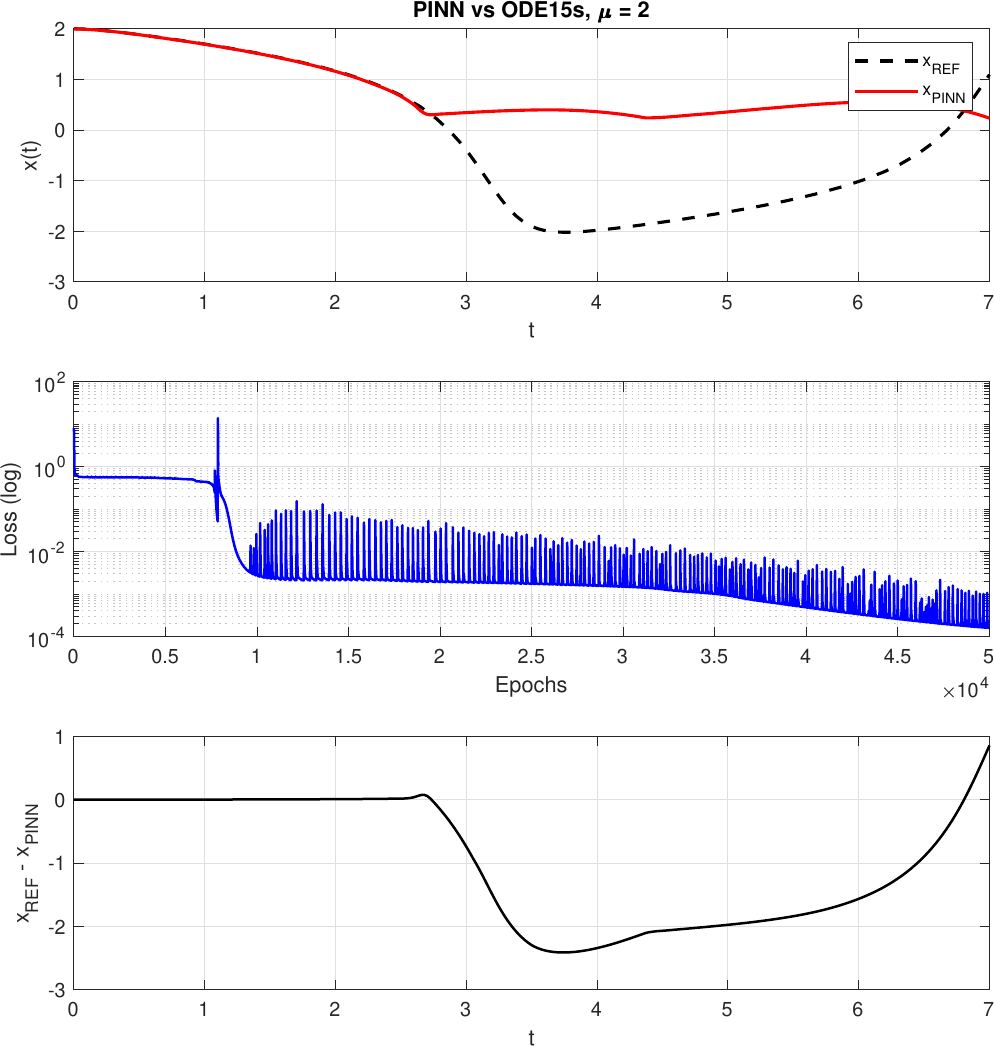}
\includegraphics[width=0.49\textwidth]{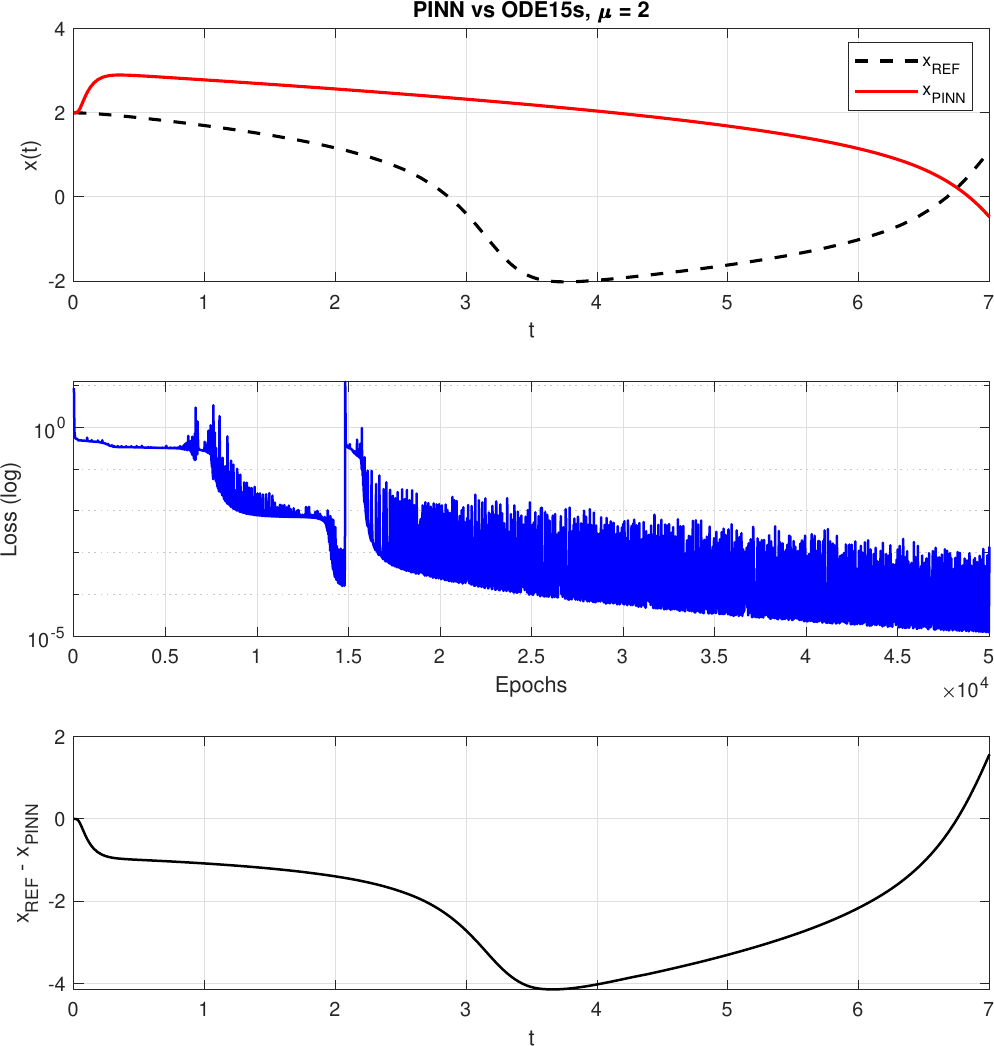}\\[-.1mm]
\mbox{ \hspace*{1mm} } \hfill (c) \hfill \mbox{} \hfill (d) \hfill \mbox{}
\caption{The examples of the PINNs solution (non-correct) to the reference numerical solution (\texttt{ode15s}) of problem \eqref{e:nonlinode} for $\mu=2$.}\label{f:nonlinodebad}
\end{figure}

In Figure \ref{f:nonlinodebad} we give examples of non-correct PINNs solutions in comparison to the reference numerical solution obtained by \texttt{ode15s} for parameter $\mu=2$. In all four cases, the loss function is very small, but the solutions are far from being satisfactory. In cases (a) and (b), although we have a satisfactory overlap of the PINNs and the reference solution at the beginning of the time interval, the PINNs solution suddenly switches to the almost zero singular solution. Moreover, in case (c), the solution does not continue as a zero solution, and in case (d) the solution is completely wrong from the beginning of the considered time interval even after 50 thousand epochs. 

\clearpage
\subsection{Heat equation problem}

We now train a PINN to numerically compute the \alert{solution} of the heat equation \eqref{e:heat} by using the \alert{limited-memory Broyden-Fletcher-Goldfarb-Shanno (L-BFGS)} algorithm.
\begin{itemize}
\item Training data generated on the go:
\begin{itemize}
\item 25 equally spaced time points to enforce each of the boundary conditions $u(t,-1)=0$ and $ u(t,1) = 0$.
\item 50 equally spaced spatial points to enforce the initial condition $u(0,x) = \alert{\phi(x)}$.
\item 10,000 interior points to enforce the output of the network to fulfil the heat equation, \alert{QRNG} (low discrepancy) Sobol sequence.
\end{itemize}
\item Deep neural network (DNN) architecture:
\begin{itemize}
\item Multilayer perceptron DNN architecture with 9 fully connected hidden layers with 20 neurons each.
\item Input layer has two input channels corresponding to the inputs $x$ and $t$, output layer has one output channel $u(x,t)$.
\item As activation functions we take hyperbolic tangents.
\item Total 3021 learnable parameters, all of type \texttt{double}.
\end{itemize}
\item Training options (custom train loop):
\begin{itemize}
\item Train for max \texttt{3e3} iterations.
\item Stop training when the norm of the gradients or steps are smaller than \texttt{1e-12} (never).
\item Use the default options for the L-BFGS solver state.
\end{itemize}
\begin{lstlisting}[firstline=142,lastline=145]
maxIterations = 3e3;
gradientTolerance = 1e-12;
stepTolerance = 1e-12; 
solverState = lbfgsState;
\end{lstlisting}
\item Model loss function: $\text{loss} = \text{MSE}_f + \text{MSE}_u$, where
\begin{itemize}
\item $\text{MSE}_f = \frac{1}{N_f}\sum\limits_{i=1}^{N_f} \left\vert f(t_f^i, x_f^i) \right\vert ^2$ 
and $\text{MSE}_u = \frac{1}{N_u}\sum\limits_{i=1}^{N_u} \left\vert u(t_u^i, x_u^i) - u^i \right\vert ^2$.
\item Here, $\{t_u^i, x_u^i \}_{i=1}^{N_u}$ correspond to collocation points on the boundary of the computational domain and account for both boundary and initial condition and $\{t_f^i, x_f^i \}_{i=1}^{N_f}$ are points in the interior of the domain.
\item Calculating $\text{MSE}_f$ requires the derivatives $\frac{\partial u}{\partial t}, \frac{\partial^2 u}{\partial x^2}$ of the output $u$ of the model -- \alert{automatic differentiation}.
\item For MATLAB implementation see Figure \ref{f:modelLoss}.
\end{itemize}
\end{itemize}

\begin{figure}[ht!]
\begin{lstlisting}[basicstyle=\scriptsize,lastline=19]
function [loss,gradients] = modelLoss(net,T,X,T0,X0,U0,opt)
TX = cat(1,X,T);
U = forward(net,TX);
gradientsU = dlgradient(sum(U,"all"),{T,X},EnableHigherDerivatives=true);
Ut = gradientsU{1}; Ux = gradientsU{2};
Uxx = dlgradient(sum(Ux,"all"),X,EnableHigherDerivatives=true);
% Calculate mseF. Enforce heat equation.
f = Ut - (opt.alpha).*Uxx;
zeroTarget = zeros(size(f),"like",f);
mseF = l2loss(f,zeroTarget);
% Calculate mseU. Enforce initial and boundary conditions.
TX0 = cat(1,X0,T0);
U0Pred = forward(net,TX0);
mseU = l2loss(U0Pred,U0);
% Calculated loss to be minimized by combining errors.
loss = mseF + mseU;
% Calculate gradients with respect to the learnable parameters.
gradients = dlgradient(loss,net.Learnables);
end
\end{lstlisting}
\begin{lstlisting}[basicstyle=\scriptsize,firstline=20]
% MATLAB (since R2021a) automatically optimizes, caches, and reuses the traces
accfun = dlaccelerate(@modelLoss); 
\end{lstlisting}
\caption{MATLAB implementation of the \texttt{modelLoss()} function.}\label{f:modelLoss}
\end{figure}

\begin{figure}[ht!]
\includegraphics[width=\textwidth]{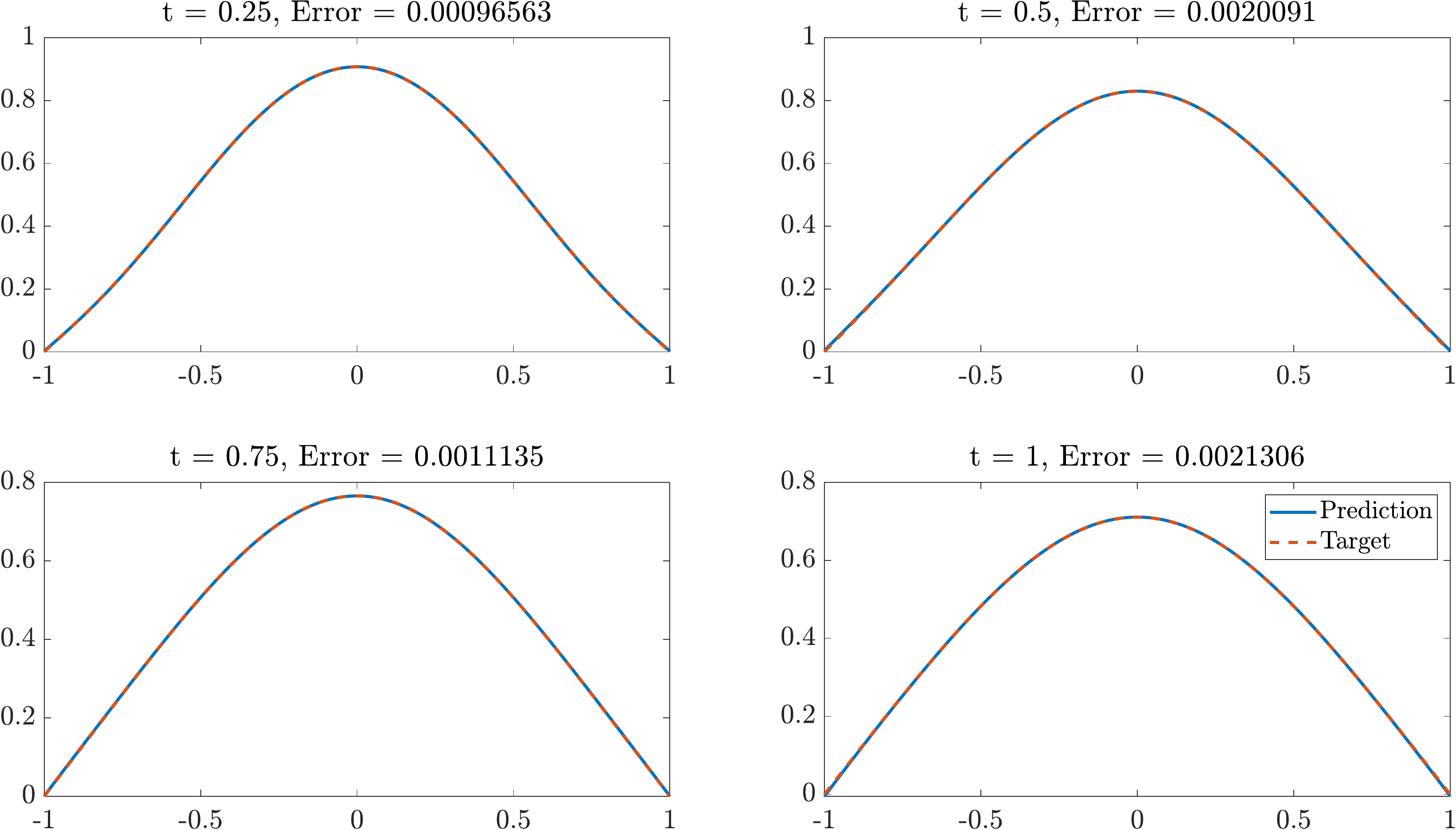}
\includegraphics[width=\textwidth]{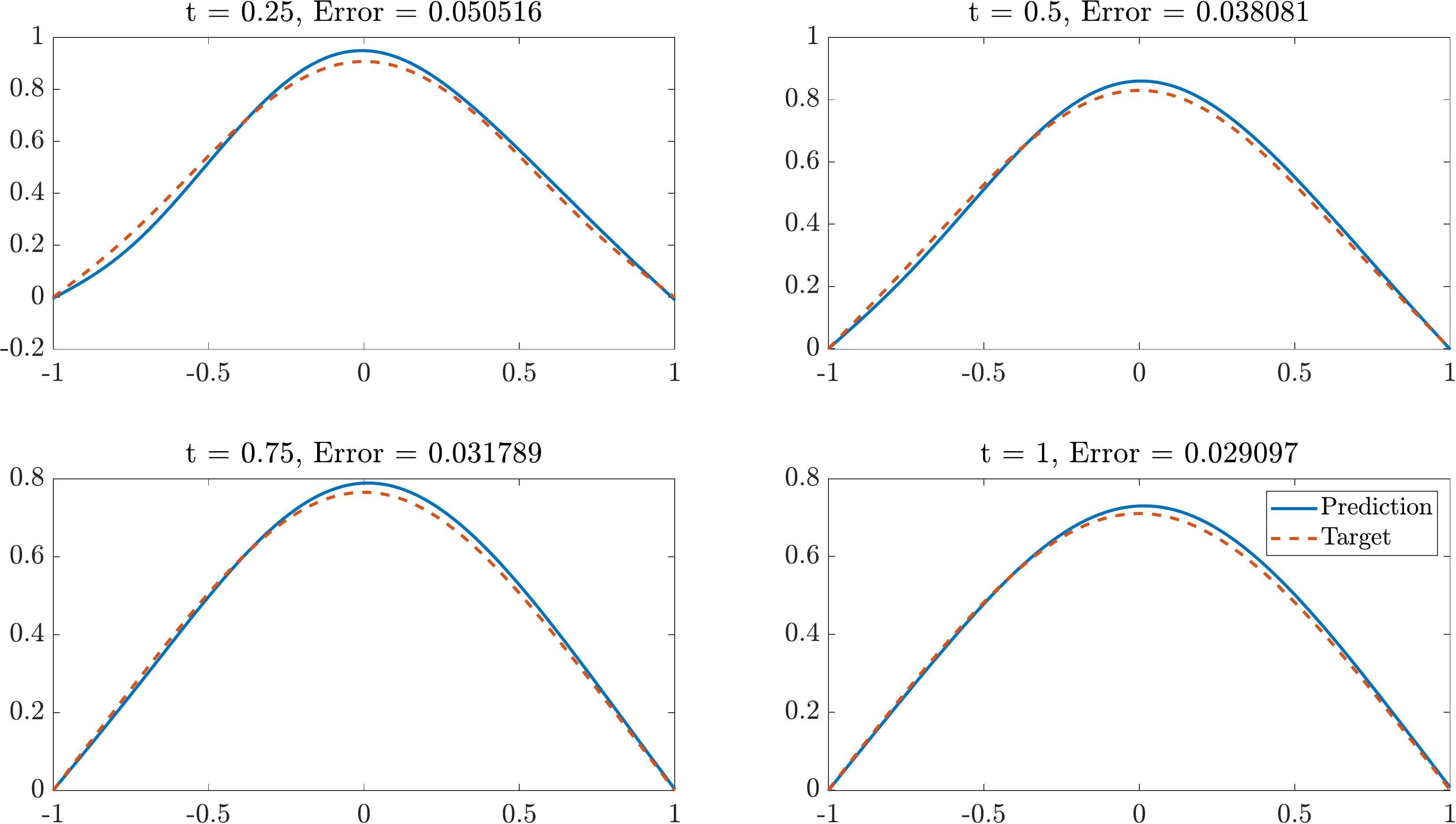}
\caption{Cut of the solution at four different timestamps. Top four pictures show the solution when \texttt{phi\_vpa()} is used as the initial condition, bottom four are for the \texttt{phi\_double()} case.}\label{f:intpinn_test}
\end{figure}

In Figure \ref{f:intpinn_test}, top four pictures show a very good fit of our DNN model to the (exact) solution, with \texttt{phi\_vpa()} being used as the initial condition. Bottom four pictures show the influence of inaccurately evaluated initial condition \texttt{phi\_double()}. The solution is highly biased. It is worth to mention that in many real world examples we usually do not have an exact solution to compare our fitted DNN to and hence to realize that there is a huge bias is far from being straightforward.

 \clearpage
\section{Conclusion}\label{sec:conclusion}

Calculating the value of a function together with its derivatives or integrals are basic numerical operations in most numerical algorithms. Often, we incorrectly automatically assume that the value calculated is within machine precision in a given arithmetic, not realizing that this may not be the case. Inaccurate function values can subsequently cause very unreliable inputs to algorithms working with that function. Yet we may not even notice that the obtained results are incorrect. 

Automatic differentiation is a core tool to automate the simultaneous computation of the numerical values of arbitrarily complex functions and their derivatives. It computes exact, machine‑precision derivatives by threading the chain rule through each operation as code runs. Unlike finite differences it needs no step‑size tuning, and unlike symbolic math it acts on program instructions, avoiding complicated expression manipulations. Gradients and Jacobians emerge alongside function evaluation in forward or reverse mode, yet they still inherit floating‑point round‑off issues. It is worth to mention that neural network back‑propagation is simply the reverse‑mode automatic differentiation \citep{Baydin2018AD} and hence many neural networks trainings can suffer from numerical precision if not treated correctly. 

In this article, we have analysed three problematic use-cases of using PINNs to solve differential equations. In particular, a linear and a non-linear stiff ODE were considered to show how changing one model parameter or NN settings produces non-satisfactory results. To demonstrate the detection difficulty, we presented several experiments with very low loss function, but with the solution being highly biased. For a linear partial differential equation we show how inaccurate function values can negatively influence the results of numerical integration (quadrature) and subsequently the PINNs solution. 

For quadrature calculations, we are interested in the accuracy of the result and the computational time. Requested optimal combination is to have an accurate result in minimum computational time. We have shown in the examples that a longer computational time does not automatically mean a more accurate result, on the contrary, it can cause problems. On the other hand, a sufficiently accurate result may not always be time consuming. It is worth mentioning that one has to be careful with low computational times that can sometimes occur for inaccurately evaluated integrand. Moreover, in practice, an integrand could be available only as a black box and hence it might be difficult to conclude that the problems are caused only by the inaccurately evaluated integrand. Therefore, for each computation, it is important to analyse the entire algorithm including all the data. In contrast to the results of \cite{DanekPospisil20ijcm}, where authors propose a regime switching algorithm to change the integrand evaluation to \texttt{vpa} only in problematic cases, here the PDE solution method is significantly improved only by evaluating the initial condition in \texttt{vpa}. 

In this paper we focused on failures of numerical calculations that are based on using the IEEE-754 Standard floating-point formats. In the last decade we can observe also an increase of usage of the IEEE-1788 Standard interval arithmetic floating-point formats or a completely new generalizing approach using the so called universal number format (unum) described for example in the book by \cite{Gustafson2015unum}. Since these approaches are still not widely used in practice we omitted them in the current study, however, it will be interesting to further research their usage not only in the PINNs modelling context.

Numerical integration challenges, as discussed in this paper, arise due to inaccuracies in evaluating integrands, leading to severe precision problems in many models, and as demonstrated in particular in physics-informed neural networks. These issues stem from floating-point arithmetic limitations and improper numerical quadrature methods, especially when model parameters exacerbate these problems. Similar challenges can occur for example in the peridynamic models, which rely on integral formulations to capture long-range interactions and avoid spatial derivatives. The evaluation of integro-differential operators in peridynamics, as seen in its bond-based models, requires precise numerical integration for accurate representation of material behaviour, including crack propagation and heterogeneous material response \citep{Haghighat2021nonlocal, Ning2023peridynamic}.

In PINNs applied to peridynamics, numerical integration problems may emerge when enforcing physical constraints, such as integro-differential operators, within the loss function. For example, in peridynamic-informed PINNs, inaccurate integration of forces across horizons or improperly handled kernel functions can degrade model accuracy, particularly under complex boundary or material conditions \citep{Ning2023peridynamic, Difonzo2024Physics}. These problems echo those highlighted in this paper, where detection of precision loss can be non-trivial without detailed numerical diagnostics. Thus, addressing these integration challenges in peridynamic models often necessitates adopting adaptive quadrature methods or higher-precision arithmetic to maintain robustness and accuracy \citep{Ning2023peridynamic}.

\section*{Acknowledgements}

Computational resources were provided by the \href{https://www.e-infra.cz/en}{e-INFRA CZ} project (ID:90254), supported by the Ministry of Education, Youth and Sports of the Czech Republic.

We thank our colleagues Kamil Ek\v{s}tein and Ond\v{r}ej Pra\v{z}\'{a}k from the Department of Computer Science and Engineering for their comments related to the implementation of PINNs.

\begin{switch}{3}
\case{1}{%
\bibliographystyle{references/styles/jp+doi+mr+zbl}
\bibliography{new-intpinns,references/references,references/references-books,references/references-own,references/preprints,references/preprints-own}
}
\case{2}{%
\bibliographystyle{jp+doi+mr+zbl}
\bibliography{new-intpinns,references-export}
}
\case{3}{%

 }
\default{}
\end{switch}

\end{document}